# Characteristic cycles and wave front cycles of representations of reductive Lie groups

By Wilfried Schmid and Kari Vilonen*

## 1. Introduction

In the papers [V1] and [BV], Vogan and Barbasch-Vogan attach two similar invariants to representations of a reductive Lie group, one by an algebraic process, the other analytic. They conjectured that the two invariants determine each other in a definite manner. Here we prove the conjecture. Our arguments involve two finer invariants – the characteristic cycles of representations – which are interesting in their own right.

To describe the invariants, we consider a linear, reductive Lie group $G_{\mathbb{R}}$ and fix a maximal compact subgroup $K_{\mathbb{R}} \subset G_{\mathbb{R}}$. We denote their Lie algebras by $\mathfrak{g}_{\mathbb{R}}$ and $\mathfrak{k}_{\mathbb{R}}$, and the complexified Lie algebras by $\mathfrak{g}$, $\mathfrak{k}$. An element $\zeta$ of the dual space $\mathfrak{g}^*$ is said to be nilpotent if it corresponds to a nilpotent element of $[\mathfrak{g},\mathfrak{g}]$ via the isomorphism $[\mathfrak{g},\mathfrak{g}] \cong [\mathfrak{g},\mathfrak{g}]^* \subset \mathfrak{g}^*$ induced by the Killing form. Via the adjoint action, the complexification $G$ of $G_{\mathbb{R}}$ acts with finitely many orbits on $\mathcal{N}^*$, the cone of all nilpotents in $\mathfrak{g}^*$. Like all coadjoint orbits, each $G$-orbit $\mathcal{O} \subset \mathcal{N}^*$ carries a distinguished (complex algebraic) symplectic structure; the intersection of $\mathcal{O}$ with $i\mathfrak{g}_{\mathbb{R}}^* = \{\, \zeta \in \mathfrak{g}^* \mid \langle \zeta, \mathfrak{g}_{\mathbb{R}} \rangle \subset i\mathbb{R} \,\}$ consists of finitely many $G_{\mathbb{R}}$-orbits, each of which is Lagrangian in $\mathcal{O}$. The choice of a maximal compact subgroup determines a Cartan decomposition $\mathfrak{g} = \mathfrak{k} \oplus \mathfrak{p}$, and dually $\mathfrak{g}^* = \mathfrak{k}^* \oplus \mathfrak{p}^*$. The complexification $K$ of the group $K_{\mathbb{R}}$ operates on $\mathcal{N}^* \cap \mathfrak{p}^*$ with finitely many orbits, and each of these orbits is Lagrangian in the $G$-orbit which contains it [KR].

Now let $\pi$ be an irreducible, admissible representation of $G_{\mathbb{R}}$ – for example an irreducible unitary representation. To such a representation, one can associate its Harish-Chandra module $V$, which is simultaneously and compatibly a module for the Lie algebra $\mathfrak{g}$ and a locally finite module for the algebraic group $K$. Then $V$ admits $K$-invariant "good filtrations," as module over the universal enveloping algebra $U(\mathfrak{g})$, relative to its canonical filtered structure. Vogan [V1] shows that the annihilator of the graded module defines an equidi-

*W. Schmid was partially supported by the NSF. K. Vilonen was partially supported by the NSA, the NSF and the Guggenheim Foundation.

mensional, $K$-invariant algebraic cycle, independently of the choice of good filtration, whose support is contained in $\mathcal{N}^*\cap\mathfrak{p}^*$. Since $K$ acts on $\mathcal{N}^*\cap\mathfrak{p}^*$ with finitely many orbits, this "associated cycle" becomes a linear combination

$$\operatorname{Ass}(\pi) \;=\; \sum a_j\,[\mathcal{O}_{\mathfrak{p},j}] \qquad (\,a_j\in\mathbb{Z}_{\geq 0}\,) \tag{1.1}$$

of fundamental cycles of $K$-orbits $\mathcal{O}_{\mathfrak{p},j}$ in $\mathcal{N}^*\cap\mathfrak{p}^*$, all of the same dimension. This is the first of the two invariants which we relate. The other is constructed in [BV], from the (Harish-Chandra) character $\Theta_\pi$. When $\Theta_\pi$ is pulled back to $\mathfrak{g}_\mathbb{R}$ via the exponential map and the argument is scaled by a multiplicative parameter $t\in\mathbb{R}_{>0}$, the resulting family of distributions has an asymptotic expansion as $t\to 0$. The Fourier transform of the leading term can be thought of as a complex linear combination of fundamental cycles of $G_\mathbb{R}$-orbits in $i\mathcal{N}^*_\mathbb{R}=\mathcal{N}^*\cap i\mathfrak{g}^*_\mathbb{R}$,

$$\operatorname{WF}(\pi) \;=\; \sum b_j\,[\mathcal{O}_{\mathfrak{g}_\mathbb{R},j}] \qquad (\,b_j\in\mathbb{C}\,)\,, \tag{1.2}$$

in the following sense: up to the multiplication by $i$, each $\mathcal{O}_{\mathfrak{g}_\mathbb{R},j}$ is a coadjoint orbit, hence carries a canonical measure, whose Fourier transform defines a distribution on $\mathfrak{g}_\mathbb{R}$. We shall call $\operatorname{WF}(\pi)$ the "wave front cycle" of $\pi$. Its support coincides with the wave front set of the distribution $\Theta_\pi$ at the identity, as was proved by Rossmann [R3], [R4].

The similar nature of the two invariants (1.1–1.2) led Barbasch and Vogan to suggest the existence of a natural bijection between the $K$-orbits in $\mathcal{N}^*\cap\mathfrak{p}^*$ and the $G_\mathbb{R}$-orbits in $i\mathcal{N}^*_\mathbb{R}$. This conjectured correspondence

$$K\backslash(\mathcal{N}^*\cap\mathfrak{p}^*) \;\longleftrightarrow\; G_\mathbb{R}\backslash i\mathcal{N}^*_\mathbb{R} \tag{1.3}$$

was established by Sekiguchi [Se] and Kostant (unpublished). We shall show:

1.4 Theorem. *The associated cycle* $\operatorname{Ass}(\pi)$ *coincides with the wave front cycle* $\operatorname{WF}(\pi)$ *via the correspondence* (1.3).

This result settles a conjecture of Barbasch-Vogan [V2]. In particular, it implies that the coefficients $b_j$ of the wave front cycle are nonnegative integers.

The so-called orbit method suggests that certain irreducible unitary representations of the reductive group $G_\mathbb{R}$ should be attached to nilpotent orbits. Ideally one would like to realize these "unipotent representations" geometrically, as spaces of sections, or perhaps cohomology groups, of line bundles on the nilpotent orbits in question. Such direct geometric constructions have been carried out only in isolated cases. On the other hand, the associated cycle and the wave front cycle attach nilpotent orbits to representations, and these nilpotent invariants can be used in the process of labeling some representations as unipotent. The affirmative answer to the Barbasch-Vogan conjecture thus settles a natural question: the two types of nilpotent invariants, arising from the associated cycle and the wave front cycle respectively, give absolutely equivalent information. Vogan's paper [V3] contains a broad survey of the notion of unipotent representation and various related matters. It is also a convenient reference for a good part of the material used by us.

Theorem 1.4 fits into a general pattern. There are several other invariants and constructions – such as $\mathfrak{n}$-homology, induction, geometric realizations – that can be carried out alternatively for Harish-Chandra modules or $G_{\mathbb{R}}$-representations; the former emphasizes the role of the $K$-action, and the latter, the role of $G_{\mathbb{R}}$. It sometimes happens that a calculation is doable on one side, but has implications on the other side. In this spirit, J.-T. Chang [C2] has used Theorem 1.4 to give a simple, conceptual proof of a theorem of Vogan [V1].

Our proof of the theorem relates the associated cycle to the wave front cycle via two geometric invariants. On the side of Harish-Chandra modules, the Beilinson-Bernstein construction attaches representations to $K$-equivariant $\mathcal{D}$-modules on the flag variety $X$, and these in turn correspond to $K$-equivariant sheaves on $X$ via the Riemann-Hilbert correspondence. The passage from Harish-Chandra modules to $K$-equivariant sheaves has an analogue on the side of $G_{\mathbb{R}}$-representations. This gives three sides of a square, with vertices "Harish-Chandra modules," "$G_{\mathbb{R}}$-representations," "$K$-equivariant sheaves," and "$G_{\mathbb{R}}$-equivariant sheaves." The fourth side, the "Matsuki correspondence for sheaves," makes the square commute; this is the commutative square (2.9) below. Kashiwara's characteristic cycle construction applies, in particular, to the $K$-equivariant and $G_{\mathbb{R}}$-equivariant sheaves arising from representations, producing Lagrangian cycles in the cotangent bundle $T^*X$. The main point of our proof is a microlocalization of the Matsuki correspondence – an explicit, geometric passage from the characteristic cycle on the $K$-equivariant side to that of the corresponding $G_{\mathbb{R}}$-equivariant sheaf.

The moment map of the $G$-action on the cotangent bundle sends $T^*X$ to the nilpotent cone $\mathcal{N}^*$. It turns out that the associated cycle is the image, in an appropriate sense, of the characteristic cycle of the $K$-equivariant sheaf under the moment map. Similarly, the characteristic cycle of the $G_{\mathbb{R}}$-equivariant sheaf determines the wave front cycle. In both cases the characteristic cycles carry more information than $\mathrm{Ass}(\pi)$ and $\mathrm{WF}(\pi)$. The characteristic cycles merit further study, we think; they may turn out to be more interesting invariants than the associated cycle and the wave front cycle.

The two types of characteristic cycles, i.e., the associated cycle and the wave front cycle, take values in four abelian groups which form the vertices of a commutative square. We already mentioned three of the arrows: the microlocalization of the Matsuki correspondence and the cycle maps induced by the moment map. The fourth arrow, the "push-down" of the microlocalized

Matsuki correspondence, is an explicit geometric passage from $K$-orbits in $\mathcal{N}^* \cap \mathfrak{p}^*$, viewed as cycles in $\mathcal{N}^*$, to $G_{\mathbb{R}}$-orbits in $\mathcal{N}^* \cap i\mathfrak{g}^*_{\mathbb{R}}$, again viewed as cycles. This, we prove, coincides with the Sekiguchi correspondence (1.3). The final ingredient of the argument identifies the push-downs of the characteristic cycles with the associated cycle and the wave front cycle, respectively.

A two-column commutative diagram, (7.1) in Section 7, encapsulates the entire argument. Here we give it in heuristic form:

$$\{\text{H-C-modules}\}@>>>\{G_{\mathbb{R}}\text{-representations}\}$$
$$@VL_KVV@VVL_{G_{\mathbb{R}}}V$$
$$\{K\text{-equivariant sheaves}\}@>\ \gamma\ >>\{G_{\mathbb{R}}\text{-equivariant sheaves}\}$$
$$@V\mathrm{CC}VV@VV\mathrm{CC}V \tag{1.5}$$
$$\{\text{characteristic cycles}\}@>\ \Phi\ >>\{\text{characteristic cycles}\}$$
$$@V\mu_*VV@VV\mu_*V$$
$$\{K\text{-orbits in } \mathcal{N}^* \cap \mathfrak{p}^*\}@>\ \phi\ >>\{G_{\mathbb{R}}\text{-orbits in } \mathcal{N}^* \cap i\mathfrak{g}^*_{\mathbb{R}}\}\,.$$

The top arrow represents some right inverse of Harish-Chandra's passage from representations to Harish-Chandra modules, such as the maximal globalization functor [S]. Beilinson-Bernstein's localization functor [B], [BB1], [BB2], followed by the Riemann-Hilbert correspondence [K1], [Me], is $L_K$. It has a $G_{\mathbb{R}}$-analogue $L_{G_{\mathbb{R}}}$, whose inverse is constructed in [KSd]. Kashiwara [K3] conjectured the Matsuki correspondence for sheaves $\gamma$, an elaboration of Matsuki's correspondence between $K$-orbits and $G_{\mathbb{R}}$-orbits in the flag variety [Ma]; the paper [MUV] establishes Kashiwara's conjectured description of $\gamma$. The arrows CC refer to Kashiwara's characteristic cycle construction [K2], [KSa]. Theorem 3.7 below describes $\Phi$, the microlocalization of the Matsuki correspondence. Its proof, in Section 4, depends heavily on the "open embedding theorem" [SV3], which describes the effect on characteristic cycles of push-forward under an open embedding. Our construction of the functor $\Phi$ leads us outside the customary real analytic context – i.e., outside the subanalytic setting; instead, we need to work inside one of the "analytic-geometric categories" of van den Dries-Miller [DM]. On the $G_{\mathbb{R}}$-side, the push-down map $\mu_*$ is based on ideas of Rossmann [R2], and the $K$-version first appears in [C1]. Both of these maps are discussed in Section 5, where we also deduce a description of $\phi$ from that of $\Phi$. Theorem 6.3 below identifies $\phi$ with the Sekiguchi correspondence. Our proof of this theorem again uses analytic-geometric categories; it also depends on certain geometric properties of nilpotent orbits which are established in [SV5]. We complete the proof of our main Theorem 1.4 in Section 7, by identifying the composition of three vertical arrows $L_K$, CC, $\mu_*$ with the associated cycle construction, and the composition of $L_{G_{\mathbb{R}}}$, CC, $\mu_*$ with the construction of the wave front cycle. The former amounts to a rephrasing of a result of J.-T. Chang [C1]. On the $G_{\mathbb{R}}$-side, we crucially use our integral formula for characters [SV4], which is based on ideas of Rossmann [R1].

Our proof of the Barbasch-Vogan conjecture was announced and sketched in [SV2]. Earlier, Chang [C1] had deduced the conjecture, in the case of complex groups, from results of Rossmann [R2].

## 2. Geometric parametrization of representations

Our hypotheses and notation are those established in [SV4]. We recall some of the results – not due to us – collected in Section 2 of that paper, which will serve as general reference. In particular, we suppose that $G_{\mathbb{R}}$ is a real form of a connected, complex, linear, reductive group $G$. We choose a maximal compact subgroup $K_{\mathbb{R}} \subset G_{\mathbb{R}}$; its complexification $K$ is a subgroup of $G$. We write $\mathfrak{g}_{\mathbb{R}}$, $\mathfrak{g}$, $\mathfrak{k}_{\mathbb{R}}$, $\mathfrak{k}$ for the Lie algebras of $G_{\mathbb{R}}$, $G$, $K_{\mathbb{R}}$, $K$.

The group $G$ acts transitively and algebraically on the flag variety $X$ of $\mathfrak{g}$. The two subgroups $K$, $G_{\mathbb{R}}$ act with finitely many orbits. According to Matsuki [Ma], the two types of orbits are in one-to-one correspondence

$$K\backslash X \;\longleftrightarrow\; G_{\mathbb{R}}\backslash X\,, \tag{2.1}$$

with a $K$-orbit $S_K$ matched to a $G_{\mathbb{R}}$-orbit $S_{G_{\mathbb{R}}}$ if and only if the two orbits intersect along exactly one $K_{\mathbb{R}}$-orbit. As in [SV4, §2], we consider the "universal Cartan algebra" $\mathfrak{h}$ for $\mathfrak{g}$. Its dual space $\mathfrak{h}^*$ contains the universal root system $\Phi$ and the universal system of positive roots $\Phi^+$, as well as the universal weight lattice $\Lambda$. We fix a "localization parameter" $\lambda \in \mathfrak{h}^*$ and introduce $\mathrm{D}_{G_{\mathbb{R}}}(X)_\lambda$, the $G_{\mathbb{R}}$-equivariant derived category with twist $(\lambda - \rho)$ [SV4, §2], and totally analogously $\mathrm{D}_K(X)_\lambda$, the $K$-equivariant derived category with the same twist. We recall that the objects of these derived categories are represented by complexes of equivariant monodromic sheaves on the enhanced flag variety. Here $\rho$ denotes the half-sum of the positive roots; thus $\mathrm{D}_{G_{\mathbb{R}}}(X)_\rho$, $\mathrm{D}_K(X)_\rho$ reduce to the usual (untwisted) equivariant derived categories.

To each $\mathcal{F} \in \mathrm{D}_{G_{\mathbb{R}}}(X)_\lambda$, one can associate a family of admissible representations of $G_{\mathbb{R}}$, as follows. We let $\mathcal{O}_X^{\mathrm{hol}}(\lambda)$ denote the twisted sheaf of holomorphic functions on $X$ with shift $(\lambda - \rho)$, i.e., the same twist as in the definition of $\mathrm{D}_{G_{\mathbb{R}}}(X)_\lambda$. Thus one can introduce the groups $\mathrm{Ext}^p(\mathcal{F}, \mathcal{O}_X^{\mathrm{hol}}(\lambda))$ by deriving the functor Hom on the category of twisted sheaves with twist $(\lambda - \rho)$. These Ext groups carry a natural action of $G_{\mathbb{R}}$ and, less obviously, a natural Fréchet topology. The resulting representations are continuous, admissible, of finite length, with infinitesimal character $\chi_\lambda$, in Harish-Chandra's notation [HC]. The groups $\mathrm{Ext}^p(\mathcal{F}, \mathcal{O}_X^{\mathrm{hol}}(\lambda))$ depend contravariantly on $\mathcal{F}$. For technical reasons, we want to make the dependence covariant, by composing it with the Verdier duality operator $\mathbb{D} : \mathrm{D}_{G_\mathbb{R}}(X)_{-\lambda} \to \mathrm{D}_{G_\mathbb{R}}(X)_\lambda$, as in [SV4]. Taking the alternating sum with respect to $p$, we obtain a map

$$\begin{aligned} \beta \ :\ \mathrm{D}_{G_\mathbb{R}}(X)_{-\lambda} \ &\longrightarrow\ \{\text{virtual } G_\mathbb{R}\text{-representations}\}_\lambda\,, \\ \beta(\mathcal{F}) \ &=\ \textstyle\sum_p (-1)^p \ \mathrm{Ext}^p(\mathbb{D}\mathcal{F}, \mathcal{O}_X^{\mathrm{hol}}(\lambda))\,. \end{aligned} \tag{2.2}$$

Here "virtual representations" is shorthand for "integral linear combination of irreducible, admissible representations," which we take up to infinitesimal equivalence. The subscript $\lambda$ refers to the infinitesimal character of the summands, namely $\chi_\lambda$.

The differential operators acting on the twisted sheaf $\mathcal{O}_X^{\mathrm{alg}}(\lambda)$ constitute a sheaf of "twisted differential operators" $\mathcal{D}_{X,\lambda}$, a sheaf (untwisted!) relative to the Zariski topology on X. We let $\mathrm{D}_K(\mathrm{Mod}_{\mathrm{coh}}(\mathcal{D}_{X,\lambda}))$ denote the bounded $K$-equivariant derived category of coherent sheaves of $\mathcal{D}_{X,\lambda}$-modules [BL], [KSd]. Objects in this category are regular holonomic because $K$ acts on $X$ with finitely many orbits; see [Ma], for example. According to Beilinson-Bernstein [BB2], the cohomology groups $\mathrm{H}^p(X,\mathfrak{M})$ of any $\mathfrak{M} \in \mathrm{D}_K(\mathrm{Mod}_{\mathrm{coh}}(\mathcal{D}_{X,\lambda}))$ are Harish-Chandra modules with infinitesimal character $\chi_\lambda$. Thus we can take the alternating sum over $p$,

$$\sum (-1)^p \ \mathrm{H}^p(X,\mathfrak{M}) \ \in\ \{\text{virtual Harish-Chandra modules }\}_\lambda\,. \tag{2.3}$$

Here, as before, the subscript $\lambda$ refers to the infinitesimal character. The covariant de Rham functor

$$\begin{aligned} \mathrm{DR} \ :\ \mathrm{D}_K(\mathrm{Mod}_{\mathrm{coh}}(\mathcal{D}_{X,\lambda})) \ &\longrightarrow\ \mathrm{D}_K(X)_{-\lambda}\,, \\ \mathrm{DR}(\mathfrak{M}) \ &=\ \mathrm{R\,Hom}_{\mathcal{D}_{X,\lambda}^{\mathrm{hol}}}(\,\mathcal{O}_X^{\mathrm{hol}}(\lambda)\,,\, \mathcal{D}_{X,\lambda}^{\mathrm{hol}} \otimes_{\mathcal{D}_{X,\lambda}} \mathfrak{M}\,) \end{aligned} \tag{2.4}$$

can be defined in the $K$-equivariant setting just as in the absolute case. Since the twisted sheaf $\mathcal{O}_X^{\mathrm{hol}}(\lambda)$ "lives" on the enhanced flag variety, the operation R Hom must be performed there, and produces a twisted sheaf with the opposite twist since R Hom is contravariant in the first variable. The de Rham functor implements the Riemann-Hilbert correspondence, which is an equivalence of categories [K1], [Me], also in the equivariant case [BL]. Thus we can take its inverse and compose it with the operation (2.3), to produce a map

$$\alpha \ :\ \mathrm{D}_K(X)_{-\lambda} \ \longrightarrow\ \{\text{virtual Harish-Chandra modules}\}_\lambda\,, \tag{2.5}$$

in complete analogy to (2.2). Results of Beilinson-Bernstein [BB1], [BB2] imply that $\alpha$ is surjective. We shall recall the relevant statements in more precise form later, in Section 7.

In the definition of the group of virtual $G_\mathbb{R}$-representations with infinitesimal character $\chi_\lambda$, we have taken representations up to infinitesimal equivalence; so formally each representation is completely determined by its Harish-Chandra module. Conversely, each irreducible Harish-Chandra module can be lifted to a $G_{\mathbb{R}}$-representation; hence

$$\{\text{virtual Harish-Chandra modules}\}_\lambda \longleftrightarrow \{\text{virtual } G_{\mathbb{R}}\text{-representations }\}_\lambda \tag{2.6}$$

is a natural bijection. Kashiwara conjectured the existence of an equivalence of categories

$$\gamma : \mathrm{D}_K(X)_{-\lambda} @> \sim >> \mathrm{D}_{G_{\mathbb{R}}}(X)_{-\lambda}\,, \tag{2.7}$$

the "Matsuki correspondence for sheaves," which was established in [MUV]. Define maps

$$X @< a << G_{\mathbb{R}} \times X @> q >> G_{\mathbb{R}}/K_{\mathbb{R}} \times X @> p >> X \tag{2.8}$$

by $a(g,x) = g^{-1}x$, $q(g,x) = (gK_{\mathbb{R}}, x)$, $p(gK_{\mathbb{R}}, x) = x$. They become $G_{\mathbb{R}} \times K_{\mathbb{R}}$-equivariant with respect to the following actions on the four spaces in (2.8), going from left to right: $(g,k) \cdot x = k \cdot x$, $(g,k) \cdot (g', x) = (gg'k^{-1}, g \cdot x)$, $(g,k) \cdot (g'K_{\mathbb{R}}, x) = (gg'K_{\mathbb{R}}, g \cdot x)$, $(g,k) \cdot x = g \cdot x$. Then any $\mathcal{F} \in \mathrm{D}_K(X)_{-\lambda}$ can be regarded as an object in $\mathrm{D}_{G_{\mathbb{R}} \times K_{\mathbb{R}}}(X)_{-\lambda}$. Thus, by equivariance, $a^!(\mathcal{F}) \in \mathrm{D}_{G_{\mathbb{R}} \times K_{\mathbb{R}}}(G_{\mathbb{R}} \times X)_{-\lambda}$. Now $K_{\mathbb{R}}$ acts freely on $G_{\mathbb{R}} \times X$, so there exists a distinguished $\tilde{\mathcal{F}} \in \mathrm{D}_{G_{\mathbb{R}}}(G_{\mathbb{R}}/K_{\mathbb{R}} \times X)_{-\lambda}$ such that $a^!\mathcal{F} \cong q^!\tilde{\mathcal{F}}$. Then $Rp_!\tilde{\mathcal{F}}$ is an object in $\mathrm{D}_{G_{\mathbb{R}}}(X)_{-\lambda}$ – in the last two steps, we have dropped the subscript $K_{\mathbb{R}}$ since $K_{\mathbb{R}}$ acts trivially. By definition, $\mathcal{F} \mapsto Rp_!\tilde{\mathcal{F}}$ is the map (2.7). The four morphisms (2.2), (2.5–2.7) fit into a diagram

$$\begin{CD} \{\text{virtual H-C-modules}\}_\lambda @> \sim >> \{\text{virtual } G_{\mathbb{R}}\text{-representations}\}_\lambda \\ @A\alpha AA @AA\beta A \\ \mathrm{D}_K(X)_{-\lambda} @> \gamma >> \mathrm{D}_{G_{\mathbb{R}}}(X)_{-\lambda}\,, \end{CD} \tag{2.9}$$

whose commutativity is implicit in Kashiwara's conjectures [K3]. In the remainder of this section, we deduce this commutativity from known facts.

The paper [KSd], which constructs the representations $\mathrm{Ext}^p(\mathcal{F}, \mathcal{O}_X^{\mathrm{hol}}(\lambda))$, also relates them to the Beilinson-Bernstein modules $\mathrm{H}^p(X, \mathfrak{M})$. Let $\mathcal{M} = D\,\mathrm{DR}(\mathfrak{M})$ be the image of $\mathfrak{M} \in \mathrm{D}_K(\mathrm{Mod}_{\mathrm{coh}}(\mathcal{D}_{X,\lambda}))$ in $\mathrm{D}_K(X)_{-\lambda}$ under the de Rham functor. Then, for all $p$, and with $n = \dim_{\mathbb{C}} X$,

(2.10) $\mathrm{H}^p(X, \mathfrak{M})$ is the dual, in the category of Harish-Chandra modules, of the Harish-Chandra module of $\mathrm{Ext}^{n-p}(\gamma(\mathcal{M}), \mathcal{O}_X^{\mathrm{hol}}(-\lambda))$

[KSd, 1.1f]. We should remark that the functor $\gamma$ of [MUV] and its analogue in [KSd] – where it is denoted by $\Phi$ – are defined differently on the surface. To see that the two functors are actually the same, we note that the operations of restricting from $G$ to $G_{\mathbb{R}}$ and from $S = G/K$ to $S_{\mathbb{R}} = G_{\mathbb{R}}/K_{\mathbb{R}}$ in [KSd, (5.8)] amount to replacing the induction from $K$ to $G$ by induction from $K_{\mathbb{R}}$ to $G_{\mathbb{R}}$. In the language of [KSd], this means identifying

$$\mathrm{D}_{K_{\mathbb{R}}}(X)_{-\lambda} \cong \mathrm{D}_{G_{\mathbb{R}} \times K_{\mathbb{R}}}(X \times G_{\mathbb{R}})_{-\lambda} \ \text{ and } \ \mathrm{D}_{G_{\mathbb{R}} \times K_{\mathbb{R}}}(X \times G_{\mathbb{R}})_{-\lambda} \cong \mathrm{D}_{G_{\mathbb{R}}}(X \times S_{\mathbb{R}})_{-\lambda}\,.$$

The first of these operations is $a^*$ in our previous notation, and the second coincides with $q^*$. The definition of $\gamma$ above involves $a^!$ and $q^!$, which accounts for a shift in degree by the complex dimension of $S$, since $a$ and $q$ are smooth morphisms. The restriction from $S$ to $S_{\mathbb{R}}$ in [KSd] involves $i^!$ rather than $i^*$, which accounts for another shift in degree; the combination of these two shifts agrees, finally, with the shift in the definition of the functor $\Phi$ in [KSd].

At this point, to deduce the commutativity of (2.9) from the statement (2.10), we need to know:

2.11 PROPOSITION. *The virtual representation* $\sum(-1)^p \operatorname{Ext}^p(\mathbb{D}\mathcal{F}, \mathcal{O}_X^{\mathrm{hol}}(\lambda))$ *is* (*up to infinitesimal equivalence*) *the dual of* $\sum(-1)^{n+p} \operatorname{Ext}^p(\mathcal{F}, \mathcal{O}_X^{\mathrm{hol}}(-\lambda))$, *for every* $\mathcal{F} \in \mathrm{D}_{G_{\mathbb{R}}}(X)_{-\lambda}$.

A more precise statement, which puts individual summands into duality with each other, was conjectured by Kashiwara [K3]. Although the conjecture is surely accessible with present techniques, no proof has appeared in the literature. Here we shall deduce the weaker statement (2.11) from our proof [SV4] of Kashiwara's conjecture on characters and fixed point formulas [K4]. We should remark that the proposition is more than a purely formal assertion; in particular, the analogous statement about the association (2.5) is incorrect. As far as we know, there does not exist an explicit interpretation of the duality of Harish-Chandra modules in terms of $K$-equivariant sheaves.

*Proof of* (2.11). Let $\tilde{G}$ denote the set of pairs $(g,x) \in G \times X$ such that $gx = x$ and $\tilde{G}_{\mathbb{R}}$ the inverse image of $G_{\mathbb{R}}$ in $\tilde{G}$, as in our paper [SV4]. Following Kashiwara, we assign a cycle $c(\mathcal{F}) \in \mathrm{H}_d^{\mathrm{inf}}(\tilde{G}_{\mathbb{R}}, \mathbb{C}_{-\lambda})$ to each $\mathcal{F} \in \mathrm{D}_{G_{\mathbb{R}}}(X)_{-\lambda}$; here $d$ is the real dimension of $G_{\mathbb{R}}$ and $\mathbb{C}_{-\lambda}$ a certain local system on the universal Cartan, pulled back to $\tilde{G}_{\mathbb{R}}$ in a natural way. We write $\Theta(\mathcal{F})$ for the character of the virtual representation $\sum(-1)^p \operatorname{Ext}^p(\mathbb{D}\mathcal{F}, \mathcal{O}_X^{\mathrm{hol}}(\lambda))$. Then

$$\int_{G_{\mathbb{R}}} \Theta(\mathcal{F})\, \phi\, dg \;=\; \int_{c(\mathcal{F})} (q^*\phi)\, \tilde{\omega} \qquad (\ \phi \in C_c^\infty(G_{\mathbb{R}})\ ) \tag{2.12}$$

[SV4, Theorem 5.12]. In this formula, $q : \tilde{G}_{\mathbb{R}} \to G_{\mathbb{R}}$ refers to the projection and $\tilde{\omega}$ denotes a differential form on $\tilde{G}$ derived from the Haar measure $dg$ and a choice of orientation of $G_{\mathbb{R}}$. The same choice of orientation is used to produce the cycle $c(\mathcal{F})$: reversing the orientation of $G_{\mathbb{R}}$ affects both $c(\mathcal{F})$ and $\tilde{\omega}$ by the factor $(-1)^d$.

The anti-involution $g \mapsto g^{-1}$ induces an anti-involution $u : \tilde{G}_{\mathbb{R}} \to \tilde{G}_{\mathbb{R}}$. Examining the fixed point formalism which produces $c(\mathcal{F})$ from $\mathcal{F}$, one finds

$$c(\mathbb{D}\mathcal{F}) \;=\; e^{-2\rho}\, u_* c(\mathcal{F})\,, \tag{2.13}$$

as was already observed in [SV1, Prop. 5.3]. Several remarks are in order. First, the twisted derived category $\mathrm{D}_{G_{\mathbb{R}}}(X)_{-\lambda}$ involves the twist $(-\lambda-\rho)$, so Verdier duality for twisted sheaves maps $\mathrm{D}_{G_{\mathbb{R}}}(X)_{-\lambda}$ to $\mathrm{D}_{G_{\mathbb{R}}}(X)_{\lambda+2\rho}$. On the other hand, $\mathrm{D}_{G_{\mathbb{R}}}(X)_{\lambda+2\rho} \cong \mathrm{D}_{G_{\mathbb{R}}}(X)_{\lambda}$ because $2\rho$ is an integral weight. The factor $e^{-2\rho}$ in (2.13) is forced by the identification $\mathrm{D}_{G_{\mathbb{R}}}(X)_{\lambda+2\rho} \cong \mathrm{D}_{G_{\mathbb{R}}}(X)_{\lambda}$ – recall that the cycle $c(\mathcal{F})$ can be viewed as a linear combination of $d$-dimensional simplices with coefficients that are sections, over the simplices in question, of the local system generated by $e^{\lambda-\rho}$ [SV4]. In our paper [SV1] we use a different convention: there $\mathrm{D}_{G_{\mathbb{R}}}(X)_{\lambda}$ refers to the derived category with twist $\lambda$; consequently the version of (2.13) stated in [SV1] does not involve the factor $e^{-2\rho}$. The formula (2.13) is ambiguous until we specify the orientations on $G_{\mathbb{R}}$ used to construct the two cycles. Since $u$ is an anti-involution, we take orientations related by $u_*$ on the two sides. With this convention, there is no sign change in (2.13), as can be checked by tracing through the construction of the cycle $c(\mathcal{F})$ in [SV4].

The character $\Theta_\pi$ of an irreducible admissible representation $\pi$ is obtained by summing the diagonal matrix coefficients of $\pi$ in the sense of distributions. Thus, for entirely formal reasons,

$$\Theta_\pi(g) \;=\; \Theta_{\pi^*}(g^{-1}) \qquad (\ \pi^* = \text{dual of } \pi\ )\,. \tag{2.14}$$

Let $\Theta^*(\mathcal{F})$ denote the virtual character dual to $\Theta(\mathcal{F})$. Then, for any test function $\phi \in C_c^\infty(G_{\mathbb{R}})$,

$$\int_{G_{\mathbb{R}}} \Theta^*(\mathcal{F})\,\phi(g)\,dg \;=\; \int_{G_{\mathbb{R}}} \Theta(\mathcal{F})\,\phi(g^{-1})\,dg \;=\; \int_{c(\mathcal{F})} (u^*q^*\phi)\,\tilde\omega\,. \tag{2.15}$$

In the first step, we have used (2.14), and the second equality follows from (2.12). The definition of $\tilde\omega$ in [SV4] implies

$$u^*\tilde\omega \;=\; (-1)^{d+n} e^{-2\rho}\tilde\omega\,. \tag{2.16}$$

We combine (2.13) with (2.16) and recall our choice of orientations of $G_{\mathbb{R}}$ on the two sides of (2.12), to conclude

$$\begin{aligned}
\int_{c(\mathcal{F})} u^*(q^*\phi)\,\tilde\omega &= (-1)^{d+n}\int_{c(\mathcal{F})} e^{2\rho}\,u^*(q^*\phi\,\tilde\omega)\\
&= (-1)^{d+n}\int_{c(\mathcal{F})} u^*(e^{-2\rho}\,q^*\phi\,\tilde\omega) \;=\; (-1)^n\int_{u_*c(\mathcal{F})} e^{-2\rho}\,q^*\phi\,\tilde\omega\\
&= (-1)^n\int_{e^{-2\rho}u_*c(\mathcal{F})} q^*\phi\,\tilde\omega \;=\; (-1)^n\int_{c(\mathbb{D}\mathcal{F})} q^*\phi\,\tilde\omega\,.
\end{aligned} \tag{2.17}$$

At this point (2.15), (2.17), and another application of (2.12) give

$$\int_{G_{\mathbb{R}}} \Theta^*(\mathcal{F})\,\phi(g)\,dg \;=\; (-1)^n\int_{G_{\mathbb{R}}} \Theta(\mathbb{D}\mathcal{F})\,\phi(g)\,dg\,, \tag{2.18}$$

in other words, the assertion of the proposition.

## 3. Microlocalization of the Matsuki correspondence

In the previous section, we described $\gamma$, the Matsuki correspondence for sheaves, as a composition of certain geometrically induced morphisms. Objects in the derived categories $\mathrm{D}_K(X)_{-\lambda}$ and $\mathrm{D}_{G_{\mathbb{R}}}(X)_{-\lambda}$ can be regarded as complexes of (semi-algebraically) constructible sheaves. As such, they have characteristic cycles in the sense of Kashiwara [K2], [KSa]. We shall now determine the effect of $\gamma$ on these characteristic cycles.

It will be convenient for us to adopt the geometric view [SV3] of characteristic cycles, which was written with the present application in mind. In particular, the characteristic cycle $\mathrm{CC}(\mathcal{F})$ of a complex of sheaves $\mathcal{F}$, constructible with respect to a particular semi-algebraic (Whitney) stratification $\mathcal{S}$, is a top dimensional cycle with infinite support on $T^*_{\mathcal{S}}X$, the union of the conormal bundles of the strata $S \in \mathcal{S}$. In regarding $\mathrm{CC}(\mathcal{F})$ as a cycle in $T^*_{\mathcal{S}}X \subset T^*X$, we treat $X$ as a real algebraic manifold. Thus it would be notationally consistent to work in $T^*(X^{\mathbb{R}})$, the real cotangent bundle of $X$, considered as a manifold without complex structure. On the other hand, at various points we do use the complex structure of $X$, for example, in putting canonical orientation on the complex manifold $X$. We therefore identify the real cotangent bundle $T^*(X^{\mathbb{R}})$ with the (holomorphic!) cotangent bundle of the complex manifold $X$,

$$T^*(X^{\mathbb{R}}) \;\cong\; T^*X\,, \tag{3.1a}$$

using the convention of [KSa, (11.1.2)]. Concretely, in terms of local holomorphic coordinates $z_j = x_j + iy_j$, $1 \leq j \leq n$,

$$dx_j \;\mapsto\; \frac{1}{2}\,dz_j\;, \qquad dy_j \;\mapsto\; \frac{-i}{2}\,dz_j\;. \tag{3.1b}$$

This convention will remain in force throughout the paper. As the cotangent bundle of the complex manifold $X$, $T^*X$ carries a canonical holomorphic, nondegenerate, closed 2-form $\sigma$. It is related to the canonical 2-form $\sigma^{\mathbb{R}}$ on $T^*(X^{\mathbb{R}})$ by the formula

$$\sigma^{\mathbb{R}} \;=\; 2\ \mathrm{Re}\,\sigma \tag{3.2}$$

via the identification (3.1) [KSa, (11.1.3)].

The actions of $K$ and $G_{\mathbb{R}}$ on $X$ are, respectively, complex and real algebraic, and both groups act with finitely many orbits. It follows that the orbit stratifications are semi-algebraic and satisfy the Whitney condition. We let $T^*_K X$ and $T^*_{G_{\mathbb{R}}} X$ denote the unions of the conormal bundles of the orbits of the two groups. They are complex or real algebraic, Lagrangian subvarieties of $T^*(X^{\mathbb{R}}) \cong T^*X$ – in the case of $K$, Lagrangian even with respect to the complex algebraic symplectic form on $T^*X$. Thus both have real dimension $2n$, with $n = \dim_{\mathbb{C}} X$ as before. Objects in $\mathrm{D}_K(X)$ and $\mathrm{D}_{G_\mathbb{R}}(X)$ are constructible with respect to the orbit stratifications, so the characteristic cycle construction defines maps CC from these two derived categories to top dimensional cycles on $T^*_K X$ and $T^*_{G_\mathbb{R}} X$, respectively. The characteristic cycle construction is local with respect to the base $X$, so CC makes sense also in the twisted case:

$$\text{(3.3)}\qquad \begin{aligned} \mathrm{CC} \ &: \ \mathrm{D}_K(X)_{-\lambda} \ \longrightarrow \ \mathrm{H}^{\inf}_{2n}(T^*_K X, \mathbb{Z})\,, \\ \mathrm{CC} \ &: \ \mathrm{D}_{G_\mathbb{R}}(X)_{-\lambda} \ \longrightarrow \ \mathrm{H}^{\inf}_{2n}(T^*_{G_\mathbb{R}} X, \mathbb{Z})\,. \end{aligned}$$

Alternatively but equivalently, we can represent objects in $\mathrm{D}_K(X)$ and $\mathrm{D}_{G_\mathbb{R}}(X)$ by complexes of constructible sheaves on the enhanced flag variety $\hat{X}$, with monodromic behavior along the fibers of $\hat{X} \to X$; as such they have characteristic cycles on $T^*\hat{X}$. Because of the monodromicity, the characteristic cycles in $T^*\hat{X}$ descend to cycles in $T^*X$, and this construction coincides with (3.1). This will be made explicit later, below the proof of Lemma 4.2. We should remark that Kashiwara defines the characteristic cycle of a sheaf as a cycle with values in the orientation sheaf of the base $X$ without choosing an orientation of $X$ first. We use the complex structure to put a definite orientation on $X$, and thus may think of the characteristic cycles as absolute cycles.

The action of $G$ on $X$ induces a Hamiltonian, complex algebraic action on the cotangent bundle $T^*X$. We let $m : T^*X \to \mathfrak{g}^*$ denote the moment map. It is $G$-equivariant, complex algebraic. In our particular situation, $m$ is easy to describe: the fiber $T^*_x X$ at any $x \in X$ is naturally isomorphic to $(\mathfrak{g}/\mathfrak{b}_x)^* = \mathfrak{b}_x^\perp \subset \mathfrak{g}^*$; here $\mathfrak{b}_x$ denotes the Borel subalgebra which fixes $x$. In terms of this identification, the moment map $m$ is the identity when we regard the fibers of $T^*X$ as subspaces of $\mathfrak{g}^*$. We choose a nondegenerate, symmetric, $G$-invariant bilinear form $B$ on $\mathfrak{g}$, which is defined over $\mathbb{R}$ and agrees with the Killing form on $[\mathfrak{g}, \mathfrak{g}]$. Then $B$ induces an isomorphism $\mathfrak{g}^* \cong \mathfrak{g}$. We write

$$\text{(3.4a)}\qquad \mu \ : \ T^*X \longrightarrow \ \mathfrak{g}$$

for the composition of this identification with the moment map $m$. We shall refer to $\mu$ as the moment map from now on. By construction, it is holomorphic and $G$-invariant. Note that

$$\text{(3.4b)}\qquad \mu(T^*X) \ = \ \mathcal{N} \ = \ \text{nilpotent cone in } \mathfrak{g}\,,$$

since $\mathfrak{b}_x^\perp \cong \mathfrak{n}_x =_{\mathrm{def}} [\mathfrak{b}_x, \mathfrak{b}_x]$ via the isomorphism $\mathfrak{g}^* \cong \mathfrak{g}$. The moment map provides a useful characterization of the subvarieties $T^*_K X$ , $T^*_{G_\mathbb{R}} X$ of $T^*X$ :

$$\text{(3.5)}\qquad T^*_K X \ = \ \mu^{-1}(\mathfrak{p})\,, \qquad T^*_{G_\mathbb{R}} X \ = \ \mu^{-1}(i\mathfrak{g}_\mathbb{R})\,.$$

To see this, we observe that a cotangent vector $(x, \xi)$ with $\xi \in T^*_x X \cong \mathfrak{n}_x$ is normal to the $G_\mathbb{R}$-orbit through $x$ precisely when $\operatorname{Re} B(\xi, i\mathfrak{g}_\mathbb{R}) = 0$, in other words, when $\xi \in \mathfrak{n}_x \cap i\mathfrak{g}_{\mathbb{R}}$. In the case of the first identity, we argue the same way.

In the following, $\operatorname{Re}\mu : T^*X \to \mathfrak{g}_{\mathbb{R}}$ refers to the real part of the moment map relative to the real form $\mathfrak{g}_{\mathbb{R}} \subset \mathfrak{g}$. We define a one parameter family of bianalytic maps:

$$(3.6)\qquad \begin{gathered} F_s : T^*X \to T^*X, \qquad s \in \mathbb{R}_{>0}, \\ F_s(\xi) = \ell(\exp(-s^{-1}\mathrm{Re}\mu(\xi)))^*\xi \qquad (\xi \in T^*X), \end{gathered}$$

where $\ell(g^{-1}) : X \to X$, for $g \in G$, denotes translation by $g^{-1}$, and $\ell(g^{-1})^*$ the induced map from $T_x^*X$ to $T_{gx}^*X$. Since $F_s$ preserves the symplectic structure on $T^*X$, $(F_s)_*(C)$ is a Lagrangian cycle, for each $C \in \mathrm{H}_{2n}^{\mathrm{inf}}(T_K^*X, \mathbb{Z})$ and each $s > 0$. We recall the notion of the limit of a family of cycles and the various equivalent ways of defining it; for a detailed discussion of these matters, see [SV3, §3].

3.7 Theorem. *For* $C \in \mathrm{H}_{2n}^{\mathrm{inf}}(T_K^*X, \mathbb{Z})$, *the limit of cycles* $\lim_{s\to 0^+}(F_s)_*(C)$ *exists and is supported on* $T_{G_{\mathbb{R}}}^*X$. *The resulting homomorphism*

$$\Phi \ : \ \mathrm{H}_{2n}^{\mathrm{inf}}(T_K^*X, \mathbb{Z}) \ \longrightarrow \ \mathrm{H}_{2n}^{\mathrm{inf}}(T_{G_{\mathbb{R}}}^*X, \mathbb{Z}), \qquad \Phi(C) \ = \ \lim_{s\to 0^+}(F_s)_*(C),$$

*coincides with the map on characteristic cycles induced by* $\gamma$. *In other words,*

$$\mathrm{D}_K(X)_{-\lambda}@> \ \gamma \ >> \mathrm{D}_{G_{\mathbb{R}}}(X)_{-\lambda}$$

$$@V\mathrm{CC}VV@VV\mathrm{CC}V$$

$$\mathrm{H}_{2n}^{\mathrm{inf}}(T_K^*X, \mathbb{Z})@> \ \Phi \ >> \mathrm{H}_{2n}^{\mathrm{inf}}(T_{G_{\mathbb{R}}}^*X, \mathbb{Z})$$

*is a commutative diagram.*

The existence of the limit is not entirely obvious. However, once it is known to exist, it must have support in $T_{G_{\mathbb{R}}}^*X$ for elementary reasons. Indeed, let $(x, \xi)$ be a cotangent vector in the boundary of $\cup_{s>0}F_s(T_K^*X)$. Then there must exist sequences $\{(x_k, \xi_k)\}$ and $\{s_k\}$, such that $\mu(x_k, \xi_k) \in \mathfrak{p}$, $s_k \to 0^+$, $F_{s_k}(x_k, \xi_k) \to (x, \xi)$. We regard the cotangent spaces of $X$ as subspaces of $\mathfrak{g}$ via $\mu$. Thus we consider $\xi$ and the $\xi_k$ as lying in $\mathfrak{g}$, and more specifically, in $\mathcal{N}$; further, $\mu(x_k, \xi_k) = \xi_k \in \mathfrak{p}$ by assumption. For $g \in G$, the induced map $\ell(g^{-1})^* : T_x^*X \to T_{gx}^*X$ reduces to $\mathrm{Ad}(g) : \mathfrak{n}_x \to \mathfrak{n}_{gx}$ when we identify the cotangent spaces with subspaces of $\mathfrak{g}$ via $\mu$. The assumption $F_{s_k}(x_k, \xi_k) \to (x, \xi)$ implies

$$(3.8)\qquad \mathrm{Ad}(\exp(s_k^{-1}\ \mathrm{Re}\ \xi_k))(\xi_k) \ \longrightarrow \ \xi\,.$$

In particular, $\mathrm{Re}\ \xi_k \to \mathrm{Re}\ \xi \in \mathfrak{p}_{\mathbb{R}}$. We choose a maximal abelian subspace $\mathfrak{a}_{\mathbb{R}}$ in $\mathfrak{p}_{\mathbb{R}}$ so that $\mathrm{Re}\ \xi \in \mathfrak{a}_{\mathbb{R}}$. Replacing the $(x_k, \xi_k)$ by appropriate $K_{\mathbb{R}}$-conjugates, we can arrange Re $\xi_k \in \mathfrak{a}_{\mathbb{R}}$ for all $k$ without destroying any of our hypotheses. Let us write

$$\xi_k \;=\; \text{Re}\,\xi_k \;+\; i \textstyle\sum_{\alpha \in R} (\text{Im}\,\xi_k)^\alpha \,, \tag{3.9}$$

where $R$ denotes the restricted root system of $(\mathfrak{g}, \mathfrak{a})$ and $(\text{Im}\,\xi_k)^\alpha$ the component of Im $\xi_k$ in the $\alpha$-root space. Then

$$\begin{aligned} \text{Ad}(\exp(s_k^{-1}\,\text{Re}\,\xi_k))(\xi_k) \;&=\; \text{Re}\,\xi_k \;+\; i \textstyle\sum_{\alpha \in R} e^{s_k^{-1}\langle \alpha, \text{Re}\,\xi_k\rangle} (\text{Im}\,\xi_k)^\alpha \,, \\ \text{and} \quad i\, e^{s_k^{-1}\langle \alpha, \text{Re}\,\xi_k\rangle} &(\text{Im}\,\xi_k)^\alpha \in i\mathfrak{g}_{\mathbb{R}} \quad \text{for each } \alpha \in R\,. \end{aligned} \tag{3.10}$$

The Cartan involution $\theta$ maps the $\alpha$-root space to the $-\alpha$-root space, and acts as multiplication by $-1$ on $\mathfrak{p}$. It follows that $\theta(\text{Im}\,\xi_k)^\alpha \;=\; -(\text{Im}\,\xi_k)^{-\alpha}$; in particular,

$$\|(\text{Im}\,\xi_k)^\alpha\| \;=\; \|(\text{Im}\,\xi_k)^{-\alpha}\|\,. \tag{3.11}$$

Since $s_k^{-1} \to +\infty$, $\langle \alpha, \text{Re}\,\xi \rangle > 0$ would imply $(\text{Im}\,\xi_k)^\alpha \to 0$; hence $(\text{Im}\,\xi_k)^{-\alpha} \to 0$ and $e^{s_k^{-1}\langle -\alpha, \text{Re}\,\xi_k\rangle}(\text{Im}\,\xi_k)^{-\alpha} \to 0$. Thus, for $\alpha \in R$,

$$\langle \alpha, \text{Re}\,\xi \rangle \;>\; 0 \quad \Longrightarrow \quad (\text{Im}\,\xi)^{-\alpha} \;=\; 0\,. \tag{3.12}$$

Let $\mathfrak{m}$ denote the centralizer of Re$\,\xi$ and $\mathfrak{u}$ the linear span of the root spaces corresponding to roots $\alpha$ with $\langle \alpha, \text{Re}\,\xi \rangle \;>\; 0$. Then $\mathfrak{m} \oplus \mathfrak{u}$ is a parabolic subalgebra of $\mathfrak{g}$, defined over $\mathbb{R}$, with nilpotent radical $\mathfrak{u}$. Because of (3.12) and the definition of $\mathfrak{m}$,

$$\text{Re}\,\xi \in \mathfrak{m} \quad \text{and} \quad \text{Im}\,\xi \in \mathfrak{m} \oplus \mathfrak{u}\,. \tag{3.13}$$

As the limit of nilpotents, $\xi$ is nilpotent. Thus (3.13) forces the nilpotence also of $\xi_{\mathfrak{m}}$, the $\mathfrak{m}$-component of $\xi$. But Re $\xi_{\mathfrak{m}} = \text{Re}\,\xi$ is semisimple and commutes with Im $\xi_{\mathfrak{m}}$. We conclude that Re $\xi = 0$; hence $\xi \in i\mathfrak{g}_{\mathbb{R}}$ and $(x, \xi) \in T^*_{G_{\mathbb{R}}} X$, as was asserted.

The somewhat lengthy proof of the remaining (and deeper) parts of the theorem occupies the next section.

## 4. Proof of Theorem 3.7

In Section 2, we had described the functor $\gamma$ in terms of the operations $a^!$, $q^!$, and $Rp_!$ induced by the maps $a, q, p$ in the diagram (2.8). We begin by identifying the effect on characteristic cycles of the passage, via $a$ and $q$, from $X$ to $G_{\mathbb{R}}/K_{\mathbb{R}} \times X$. We mentioned already that we regard characteristic cycles as geometric cycles – not as cycles with values in an orientation sheaf as in [KSa] – by putting a definite orientation on the ambient space. In the case of $X$, we use the complex structure, and in the case of $G_{\mathbb{R}}/K_{\mathbb{R}}$, some as yet unspecified orientation; the particular choice will not matter. We follow the conventions of [SV3] in all orientation questions. In particular, we orient products of manifolds by choosing forms of top degree on the two factors which are positive with respect to the orientations; we then orient the product by the wedge product of the two forms, in the given order of the factors. Also recall the rule [SV3, (2.3)] for orienting the conormal bundle of a submanifold of an oriented manifold.

We consider a particular $\mathcal{F} \in \mathrm{D}_K(X)_{-\lambda}$ and its characteristic cycle $\mathrm{CC}(\mathcal{F})$. As a top dimensional cycle in $T^*_K X$, we can express it as an integral linear combination of conormal bundles of $K$-orbits in $X$,

$$\mathrm{CC}(\mathcal{F}) \;=\; \textstyle\sum_\ell c_\ell \, [T^*_{K\cdot x_\ell} X]\,, \tag{4.1}$$

with $x_\ell$ running over a complete set of coset representatives. We had argued in Section 2 that there exists a canonical $\tilde{\mathcal{F}} \in \mathrm{D}_{G_{\mathbb{R}}}(G_{\mathbb{R}}/K_{\mathbb{R}} \times X)_{-\lambda}$ such that $q^! \tilde{\mathcal{F}} = a^! \mathcal{F}$.

4.2 Lemma. *For each* $\ell$, $M_\ell = \{(gK_{\mathbb{R}}, gkx_\ell) \in G_{\mathbb{R}}/K_{\mathbb{R}} \times X \mid g \in G_{\mathbb{R}}, k \in K\}$ *is a real algebraic submanifold of* $G_{\mathbb{R}}/K_{\mathbb{R}} \times X$ , *and*

$$\mathrm{CC}(\tilde{\mathcal{F}}) \;=\; \textstyle\sum_\ell c_\ell \,(-1)^{\dim \mathfrak{p}_{\mathbb{R}}} \, [T^*_{M_\ell}(G_{\mathbb{R}}/K_{\mathbb{R}} \times X)]\,.$$

*Proof.* We observe first of all that dim $\mathfrak{p}_{\mathbb{R}}$ is the difference of the dimensions of the fibers of the maps $a$ and $q$, i.e., of the dimensions of $G_{\mathbb{R}}$ and $G_{\mathbb{R}}/K_{\mathbb{R}}$. Both $a$ and $q$ are smooth fibrations; hence $q^*$ agrees with $q^!$ except for a shift in degree equal to the fiber dimension, and similarly in the case of $a$; see [KSa, 3.3.2], for example. Note that the preceding statements use the explicit isomorphism $\mathbb{D}_M = or_M[\dim_{\mathbb{R}} M] \cong \mathbb{C}_M[\dim_{\mathbb{R}} M]$ given by the orientation of any smooth manifold $M$. Thus we can drop the sign factor in the statement of the lemma when we replace $a^!, q^!$ by $a^*, q^*$. The characteristic cycle of a sheaf is a local invariant, and locally $a, q$ are products. Since $a^{-1}(K \cdot x_\ell) = q^{-1}(M_\ell)$ by definition, our assertion now follows from these three facts:

$$\mathrm{CC}(\mathbb{C}_S) \;=\; [T^*_S M] \tag{4.3a}$$

when $S$ is a closed submanifold of an oriented manifold $M$;

$$\mathrm{CC}(\mathbb{C}_{M_1} \boxtimes \mathcal{F}) \;=\; [M_1] \times \mathrm{CC}(\mathcal{F}) \tag{4.3b}$$

when $M_1, M_2$ are oriented manifolds and $\mathcal{F}$ a semi-algebraically constructible sheaf on $M_2$; and

$$T^*_{M_1 \times S}(M_1 \times M_2) = M_1 \times T^*_S M_2\,, \qquad \text{as oriented manifolds,} \tag{4.3c}$$

with $M_1, M_2$ as before and $S \subset M_2$ a submanifold. These three statements are contained in the formalism and conventions of [KSa], but can also be deduced directly from the geometric definition of characteristic cycles and the sign conventions in [SV3], specifically (2.3–2.7) and the convention for orienting a product.

The question of sign is the only subtle matter in the proof of Lemma 4.2. At first glance, it might appear that the signs of characteristic cycles depends on whether we regard twisted sheaves on $X$ as monodromic sheaves $\hat{X}$ or – as we have chosen to do – locally as sheaves on $X$. Not so: the fiber $H$ of $\hat{X} \to X$ is complex, hence even dimensional; when we treat this fibration locally as a product, it does not matter where we place the factor $H$.

With $\mathcal{F}$ as in (4.1–4.2), we need to determine $\mathrm{CC}(Rp_!\tilde{\mathcal{F}})$. We shall do so using Theorem 6.9 in [SV3]. This requires a compactification $\overline{G_{\mathbb{R}}/K_{\mathbb{R}}}$ of $G_{\mathbb{R}}/K_{\mathbb{R}}$, as well as a function $f : \overline{G_{\mathbb{R}}/K_{\mathbb{R}}} \to \mathbb{R}$ which vanishes on $\partial(G_{\mathbb{R}}/K_{\mathbb{R}})$ and is strictly positive on $G_{\mathbb{R}}/K_{\mathbb{R}}$. The statement of [SV3, (6.9)] requires the compactification and the function $f$ to be semi-algebraic. To make our computation of $\mathrm{CC}(Rp_!\tilde{\mathcal{F}})$ manageable, we shall need to work with a certain specific function $f$ which is not even subanalytic. However, with the present application in mind, we showed in [SV3, §10] how to extend the validity of Theorem 6.9 and its generalization 6.10 beyond the semi-algebraic and subanalytic contexts. Specifically, the two theorems apply in the setting of any analytic-geometric category $\mathcal{C}$ as defined in [DM] – the sheaf $\tilde{\mathcal{F}}$ must be constructible with respect to a $\mathcal{C}$-stratification of $G_{\mathbb{R}}/K_{\mathbb{R}} \times X$ which is extendable to a $\mathcal{C}$-stratification of the compactification, and the function $f$ must be a $\mathcal{C}$-function of class $C^1$.

In the argument below, we shall use the analytic geometric category [DM] which corresponds to $\mathbb{R}_{\mathrm{an,exp}}$, in the notation of [DMM]. As in Section 3, we choose a nondegenerate, symmetric, $G$-invariant bilinear form $B$ on $\mathfrak{g}$, which is defined over $\mathbb{R}$ and agrees with the Killing form on $[\mathfrak{g},\mathfrak{g}]$. We compactify $G_{\mathbb{R}}/K_{\mathbb{R}} \cong \mathfrak{p}_{\mathbb{R}}$ by adding a single point,

$$\overline{G_{\mathbb{R}}/K_{\mathbb{R}}} \;\cong\; \overline{\mathfrak{p}_{\mathbb{R}}} \;=\; \mathfrak{p}_{\mathbb{R}} \cup \{\infty\}\,, \tag{4.4}$$

with the real analytic structure of the standard sphere containing the Euclidean space $\mathfrak{p}_{\mathbb{R}}$ as the complement of $\{\infty\}$. The function

$$\begin{gathered} f : \overline{\mathfrak{p}_{\mathbb{R}}} \;\to\; \mathbb{R}\,, \quad \text{defined by} \\ f(\infty) = 0\,, \qquad f(\zeta) = e^{-\frac{1}{2}B(\zeta,\zeta)} \;\; \text{for } \zeta \in \mathfrak{p}_{\mathbb{R}}\,, \end{gathered} \tag{4.5}$$

takes real values since $B$ is defined over $\mathbb{R}$. It is $C^\infty$ because $B > 0$ on $\mathfrak{p}_{\mathbb{R}}$, and visibly a $\mathcal{C}$-function: the map $\zeta \mapsto \|\zeta\|^{-2}\zeta$ gives a coordinate system at infinity, and the graph of $\zeta \mapsto e^{-\frac{1}{2}\|\zeta\|^{-2}}$ is a $\mathcal{C}$-set.

Let us argue that $\tilde{\mathcal{F}}$ is indeed constructible with respect to a $\mathcal{C}$-stratification of $\mathfrak{p}_{\mathbb{R}} \times X$ which extends to a $\mathcal{C}$-stratification of $\overline{\mathfrak{p}_{\mathbb{R}}} \times X$. Since $\tilde{\mathcal{F}}$ is constructible with respect to a semi-algebraic stratification of $G_{\mathbb{R}}/K_{\mathbb{R}} \times X$ , it suffices to show

$$\text{if } S \subset \mathfrak{p}_{\mathbb{R}} \times X \text{ is semi-algebraic with respect to the real algebraic structure coming from } G_{\mathbb{R}}/K_{\mathbb{R}} \times X, \text{ then } S \text{ is a } \mathcal{C}\text{-set in } \overline{\mathfrak{p}_{\mathbb{R}}} \times X\,. \tag{4.6}$$

Let $\mathfrak{a}_{\mathbb{R}} \subset \mathfrak{p}_{\mathbb{R}}$ be a maximal abelian subspace and $A_{\mathbb{R}} \cong (\mathbb{R}_{>0})^r$ the connected subgroup of $G_{\mathbb{R}}$ generated by $\mathfrak{a}_{\mathbb{R}}$. Every $K_{\mathbb{R}}$-orbit in $\mathfrak{p}_{\mathbb{R}}$ meets $\mathfrak{a}_{\mathbb{R}}$; hence

$$K_{\mathbb{R}} \times A_{\mathbb{R}} \longrightarrow G_{\mathbb{R}}/K_{\mathbb{R}}\,, \qquad (k,a) \mapsto kaK_{\mathbb{R}}\,, \tag{4.7}$$

is a surjective algebraic map. It follows that $\tilde{S} \subset K_{\mathbb{R}} \times A_{\mathbb{R}} \times X$ , the inverse image of the semi-algebraic subset $S$ of $G_{\mathbb{R}}/K_{\mathbb{R}} \times X$ , is semi-algebraic. Note that $S$ is the image in $\mathfrak{p}_{\mathbb{R}} \times X$ of $\tilde{S}$ under the map $K_{\mathbb{R}} \times A_{\mathbb{R}} \times X \to \mathfrak{p}_{\mathbb{R}} \times X$, which is induced by

$$K_{\mathbb{R}} \times A_{\mathbb{R}} \longrightarrow \mathfrak{p}_{\mathbb{R}}\,, \qquad (k, \exp\zeta) \mapsto \mathrm{Ad}(k)\zeta\,. \tag{4.8}$$

To conclude (4.6), it now suffices to extend the map (4.8) to a $\mathcal{C}$-map from an algebraic compactification of $K_{\mathbb{R}} \times A_{\mathbb{R}}$ to $\overline{\mathfrak{p}_{\mathbb{R}}}$. What matters are the following two general properties of $\mathcal{C}$-maps and $\mathcal{C}$-sets: a) the product of two $\mathcal{C}$-maps is a $\mathcal{C}$-map, and b) the image of a $\mathcal{C}$-set under a proper $\mathcal{C}$-map is a $\mathcal{C}$-set.

The map (4.8) is $K_{\mathbb{R}}$-equivariant with respect to the action by left translation on itself, the trivial action on $A_{\mathbb{R}}$, and the adjoint action on $\mathfrak{p}_{\mathbb{R}}$; moreover, $K_{\mathbb{R}}$ acts algebraically on $\mathfrak{p}_{\mathbb{R}}$, relative to the algebraic structure coming from $G_{\mathbb{R}}/K_{\mathbb{R}}$. The factor $K_{\mathbb{R}}$ in (4.8) is therefore innocuous: it suffices to complete the map $A_{\mathbb{R}} \to \mathfrak{p}_{\mathbb{R}}$, $\exp\zeta \mapsto \zeta$ to a $\mathcal{C}$-map between their compactifications. We choose coordinates $y_1, \dots, y_r, y_{r+1}, \dots, y_q$ in the compactification of $\mathfrak{p}_{\mathbb{R}}$, centered at $\infty$, so that $(\sum_k y_k^2)^{-1} y_j$ , $1 \le j \le r$, are linear coordinates on $\mathfrak{a}_{\mathbb{R}}$, and $y_j = 0$ on $\mathfrak{a}_{\mathbb{R}}$ for $j > r$. The linear coordinates on $\mathfrak{a}_{\mathbb{R}}$ induce isomorphisms $\mathfrak{a}_{\mathbb{R}} \cong \mathbb{R}^r$ and $A_{\mathbb{R}} \cong (\mathbb{R}_{>0})^r$. We compactify $A_{\mathbb{R}}$ algebraically by viewing $x_j$ and $x_j^{-1}$ as algebraic coordinates on the $j^{\text{th}}$ factor $\mathbb{R}_{>0}$ of $A_{\mathbb{R}} \cong (\mathbb{R}_{>0})^r$, near 0 and $\infty$, respectively. With these choices of coordinates, the graph of the map $A_{\mathbb{R}} \to \mathfrak{p}_{\mathbb{R}}$ is given by the equations

$$y_j \;=\; \varepsilon_j \,(\textstyle\sum_k y_k^2)\, \log x_j\,, \quad 1 \le j \le r\,; \qquad y_j \;=\; 0\,, \quad r < j \le q\,, \tag{4.9}$$

with $\varepsilon_j = \pm 1$ specifying a particular component in the boundary of $A_{\mathbb{R}} \cong (\mathbb{R}_{>0})^r$. The equations (4.9) characterize the graph as a $\mathcal{C}$-set. It is clear that the map extends at least continuously to the boundary. The graph of the continuous extension is the closure of the graph, hence a $\mathcal{C}$-set. By definition, a $\mathcal{C}$-map is a continuous map whose graph is a $\mathcal{C}$-set. Thus $A_{\mathbb{R}} \to \mathfrak{p}_{\mathbb{R}}$ can indeed be compactified as a $\mathcal{C}$-map. This completes the verification of (4.6). Note that the particular nature of the compactification of $\mathfrak{p}_{\mathbb{R}}$ plays no role in Theorems 6.9 and 6.10 of [SV3] – all that matters is the existence of some compactification with the right properties.

For the statement of our next lemma, we fix a particular $K$-orbit $Q = K \cdot x_\ell$ in $X$, and let $M = M_\ell$ denote the corresponding submanifold of $G_\mathbb{R}/K_\mathbb{R} \times X$ defined in the statement of Lemma 4.2. Recall the definition (3.6) of the family of bianalytic maps $F_s$.

4.10 LEMMA. *For* $0 < s < \infty$, *the submanifold* $T^*_M(G_\mathbb{R}/K_\mathbb{R} \times X) - s\, d\log f$ *of* $T^*(G_\mathbb{R}/K_\mathbb{R} \times X)$ *intersects the submanifold* $G_\mathbb{R}/K_\mathbb{R} \times T^*X$ *transversely along a submanifold* $N_s$. *The projection* $G_\mathbb{R}/K_\mathbb{R} \times T^*X \to T^*X$ *maps* $N_s$ *isomorphically onto* $F_s(T^*_Q X)$.

The transversality statement in this lemma will allow us to apply the results of [SV3] in the present context. Let $\mathcal{F} \in \mathrm{D}_K(X)_{-\lambda}$ be given, and let $\tilde{\mathcal{F}} \in \mathrm{D}_{G_\mathbb{R}}(G_\mathbb{R}/K_\mathbb{R} \times X)_{-\lambda}$ be the distinguished sheaf such that $q^!\tilde{\mathcal{F}} = Da^!\mathcal{F}$. To calculate the characteristic cycle of $\gamma(\mathcal{F}) = Rp_!\tilde{\mathcal{F}}$ in terms of $\mathrm{CC}(\mathcal{F})$, we appeal to Theorem 6.9 of [SV3], or more precisely, to its generalization for $\mathcal{C}$-constructible sheaves as described in [SV3, §10]. The first ingredient, namely the constructibility of $\tilde{\mathcal{F}}$ with respect to a $\mathcal{C}$-stratification of $\overline{G_\mathbb{R}/K_\mathbb{R}} \times X$, has already been established in (4.6). We have also produced a $\mathcal{C}$-function $f$ on $\overline{G_\mathbb{R}/K_\mathbb{R}}$ which vanishes precisely on the boundary. The transversality hypothesis of Theorem 6.9, finally, follows from Lemma 4.2 and the transversality statement in Lemma 4.10. Theorem 6.9 is stated for $Rp_*$ rather than for $Rp_!$. We need to apply the version for $Rp_!$, which is completely analogous to the other version, both in statement and in proof, with one exception: one uses the open embedding theorem for $Rj_*$, which involves addition of the term $s d\log f$, and the other the open embedding theorem for $Rj_!$, which involves subtracting $s d\log f$ – see [SV3, Theorem 4.2]. Alternatively, one can deduce the version of Theorem 6.9 for $Rp_!$ from the stated version by appealing to the following two facts. The operation of Verdier duality relates $Rp_!$ to $Rp_*$, $Rp_! = \mathbb{D}Rp_*\mathbb{D}$. Also, for any constructible sheaf $\mathcal{G}$, $\mathrm{CC}(\mathbb{D}\mathcal{G}) = A_*(\mathrm{CC}(\mathcal{G}))$, with $A$ denoting the antipodal map on the cotangent bundle, i.e., the bundle map which acts as multiplication by $-1$ on the fibers.

*Proof of* 4.10. The one form $d\log f$ on $\mathfrak{p}_\mathbb{R}$ defines a section of $\mathfrak{p}_\mathbb{R} \times \mathfrak{p}_\mathbb{R} \cong T^*\mathfrak{p}_\mathbb{R} \to \mathfrak{p}_\mathbb{R}$, namely "minus the identity" on the fiber,

$$d\log f \;=\; \{\, (\zeta, -\zeta) \mid \zeta \in \mathfrak{p}_\mathbb{R} \,\} \;\subset\; \mathfrak{p}_\mathbb{R} \times \mathfrak{p}_\mathbb{R} \;\cong\; T^*\mathfrak{p}_\mathbb{R}\,; \tag{4.11}$$

indeed, the function $f$ was chosen to put $d\log f$ into this particularly simple form. We want to lift this description of the one form $d\log f$ to $G_\mathbb{R}/K_\mathbb{R}$ via the identification $G_\mathbb{R}/K_\mathbb{R} \cong \mathfrak{p}_\mathbb{R}$ induced by the exponential map. For this purpose, we identify

$$T^*_{eK_\mathbb{R}}(G_\mathbb{R}/K_\mathbb{R}) \;\cong\; (\mathfrak{g}_\mathbb{R}/\mathfrak{k}_\mathbb{R})^* \;\cong\; {\mathfrak{p}_\mathbb{R}}^* \;\cong\; \mathfrak{p}_\mathbb{R}\,; \tag{4.12}$$

the last step refers to the isomorphism determined by the symmetric form $B$. We claim:

$$(4.13) \qquad d\log f \;=\; \{\, (\exp\zeta\cdot K_{\mathbb{R}}\,,\, -\ell^*_{\exp(-\zeta)}\zeta\,) \mid \zeta\in\mathfrak{p}_{\mathbb{R}}\,\} \;\subset\; T^*(G_{\mathbb{R}}/K_{\mathbb{R}})\,.$$

To deduce this from (4.11), recall the formula

$$(\exp_*)|_\zeta \;=\; (\ell_{\exp\zeta})_*\circ\frac{1-e^{-\operatorname{ad}\zeta}}{\operatorname{ad}\zeta}$$

(see, for example, [He, §II, Theorem 1.7]). Dually,

$$(4.14a) \qquad \begin{gathered}\exp^*\;:\;T^*_{\exp\zeta\cdot K_{\mathbb{R}}}(G_{\mathbb{R}}/K_{\mathbb{R}})\;\longrightarrow\;{\mathfrak{p}_{\mathbb{R}}}^*\;\cong\;\mathfrak{p}_{\mathbb{R}}\\ \text{coincides with}\quad\left(\frac{e^{\operatorname{ad}\zeta}-1}{\operatorname{ad}\zeta}\right)^*\circ\ell^*_{\exp(-\zeta)}\,.\end{gathered}$$

But $\operatorname{ad}\zeta$, for $\zeta\in\mathfrak{p}_{\mathbb{R}}$, is self-adjoint with respect to $B$, so

$$(4.14b) \qquad \left(\frac{e^{\operatorname{ad}\zeta}-1}{\operatorname{ad}\zeta}\right)^*(\zeta)\;=\;(\,1+\frac{1}{2}\operatorname{ad}^*\zeta+\dots)(\zeta)\;=\;\zeta\,.$$

Thus (4.11) and (4.14) do imply (4.13).

We need to describe the conormal bundle of the submanifold $M\subset(G_{\mathbb{R}}/K_{\mathbb{R}})\times X$. At the typical point $x=k\cdot x_\ell$ in $Q=K\cdot x_\ell$, we identify

$$(4.15) \qquad T^*_xX\;=\;(\mathfrak{g}/\mathfrak{b}_x)^*\;\cong\;\mathfrak{n}_x\quad\text{and}\quad(T^*_QX)_x\;\cong\;\mathfrak{p}\cap\mathfrak{n}_x\,.$$

Note that $M=\{(\exp\zeta\cdot K_{\mathbb{R}},\exp\zeta\cdot x)\mid\zeta\in\mathfrak{p}_{\mathbb{R}},x\in Q\}$. With convention (4.15), for $x\in Q$ and $\zeta\in\mathfrak{p}_{\mathbb{R}}$, we have

$$(4.16) \qquad \begin{aligned}&(T^*_M(G_{\mathbb{R}}/K_{\mathbb{R}}\times X))_{(\exp\zeta\cdot K_{\mathbb{R}},\exp\zeta\cdot x)}\\ &\quad=\;\{\,(\,-(\ell_{\exp(-\zeta)})^*\operatorname{Re}\eta\,,\,(\operatorname{Ad}\exp\zeta)\eta\,)\mid\eta\in\mathfrak{p}\cap\mathfrak{n}_x\,\}\,.\end{aligned}$$

Indeed, since $M$ is a union of $G_{\mathbb{R}}$-orbits, it suffices to check this description for $\zeta=0$. Recall the definition of the maps $a$ and $q$. Both are submersions, and $a^{-1}(Q)=q^{-1}(M)$; hence

$$(4.17) \qquad \begin{aligned}q^*(T^*_M(G_{\mathbb{R}}/K_{\mathbb{R}}\times X)_{(eK_{\mathbb{R}},x)})\;&=\;(q^*T^*_M(G_{\mathbb{R}}/K_{\mathbb{R}}\times X))_{(e,x)}\\ &=\;(a^*T^*_QX)_{(e,x)}\;=\;a^*((T^*_QX)_x)\,.\end{aligned}$$

We identify the ambient space $T^*_{(e,x)}(G_{\mathbb{R}}\times X)\cong{\mathfrak{g}_{\mathbb{R}}}^*\oplus(\mathfrak{g}/\mathfrak{b}_x)^*$ with $\mathfrak{g}_{\mathbb{R}}\oplus\mathfrak{n}_x$ via $B$ as usual. Then

$$(4.18) \qquad a^*((T^*_QX)_x)\;\cong\;\{\,(-\operatorname{Re}(\eta),\eta))\mid\eta\in\mathfrak{p}\cap\mathfrak{n}_x\,\}\,,$$

since $(T^*_QX)_x\cong\mathfrak{p}\cap\mathfrak{n}_x$, and since $a^*:T^*_xX\to T^*_{(e,x)}(G_{\mathbb{R}}\times X)$ is dual to the map $\mathfrak{g}_{\mathbb{R}}\oplus\mathfrak{g}/\mathfrak{b}_x\to\mathfrak{g}/\mathfrak{b}_x$ which is the identity on the second summand and minus the identity of the first summand, followed by the projection to $\mathfrak{g}/\mathfrak{b}_x$. But $q^*$ is injective on each cotangent space, so (4.17–4.18) imply (4.16) at $\zeta=0$, and therefore in general.

Combining (4.13) and (4.16), we see that the intersection $N_s$ of $G_{\mathbb{R}}/K_{\mathbb{R}} \times T^*X$ with $T^*_M(G_{\mathbb{R}}/K_{\mathbb{R}} \times X) - s d\log f$ is given by the equation $s\zeta = \operatorname{Re}\eta$, which visibly describes a transverse intersection. Explicitly,

(4.19)
$$N_s \;=\; \{\,(\exp\zeta\cdot K_{\mathbb{R}}, \exp\zeta\cdot x, 0, \operatorname{Ad}\exp(\zeta)\eta) \mid x\in Q\,,\ \eta\in\mathfrak{p}\cap\mathfrak{n}_x\,,\ s\zeta = \operatorname{Re}\eta\,\}\,.$$

Note that the first variable can be recovered from the others, so the projection to $T^*X$ maps $N_s$ isomorphically to its image. To identify this image with $F_s(T^*_QX)$, we only need to observe that

$$(4.20)\qquad \begin{array}{rl} & \ell(\exp(-\zeta))^* \;:\; T^*_xX \to T^*_{\exp(\zeta)x}X \\ \text{corresponds to} & \operatorname{Ad}(\exp\zeta) : \mathfrak{n}_x \to \mathfrak{n}_{\exp(\zeta)x} \end{array}$$

via the natural isomorphisms $T^*_xX \cong \mathfrak{n}_x$ and $T^*_{\exp(\zeta)x}X \cong \mathfrak{n}_{\exp(\zeta)x}$. This completes the proof of Lemma 4.10.

Let us summarize what needs to be done to complete the proof of Theorem 3.7. Using the notation established at the beginning of this section, we have $\gamma(\mathcal{F}) = Rp_!\tilde{\mathcal{F}}$; hence

$$(4.21)\qquad \operatorname{CC}(\gamma(\mathcal{F})) \;=\; \operatorname{CC}(Rp_!\tilde{\mathcal{F}})\,.$$

Following the statement of Lemma 4.10, we had argued that we can apply Theorem 6.9 of [SV3], in the context of $\mathcal{C}$-maps and $\mathcal{C}$-functions, and for $Rp_!$ instead of $Rp_*$. We recall the relevant statement. The differential of the projection $p$ induces

$$(4.22)\qquad dp \;:\; G_{\mathbb{R}}/K_{\mathbb{R}} \times T^*X \;\hookrightarrow\; T^*(G_{\mathbb{R}}/K_{\mathbb{R}} \times X)\,,$$

the inclusion of $\{\text{zero section}\}\times T^*X$ into the cotangent bundle of the product. Projecting to the first factor, we get

$$(4.23)\qquad \tau \;:\; G_{\mathbb{R}}/K_{\mathbb{R}} \times T^*X \;\longrightarrow\; T^*X\,.$$

With this notation, the $Rp_!$-version of Theorem 6.9 asserts

$$(4.24)\qquad \operatorname{CC}(Rp_!\tilde{\mathcal{F}}) \;=\; \lim_{s\to 0^+}\ \tau_*\,(dp)^{-1}(\operatorname{CC}(\tilde{\mathcal{F}}) - s d\log f)\,.$$

Thus it suffices to equate $\tau_*\,(dp)^{-1}(\operatorname{CC}(\tilde{\mathcal{F}}) - s d\log f)$ with $F_s(\operatorname{CC}(\mathcal{F}))$, for $s > 0$.

Lemma 4.2 reduces the problem to the case of the conormal bundle of a single $K$-orbit. Let $Q = K\cdot x_\ell \subset X$ be a $K$-orbit and $M = M_\ell$ the corresponding submanifold of $G_{\mathbb{R}}/K_{\mathbb{R}} \times X$. We let $[T^*_QX]$ play the role of $\operatorname{CC}(\mathcal{F})$, and accordingly, $(-1)^{\dim\mathfrak{p}_{\mathbb{R}}}\,[T^*_M(G_{\mathbb{R}}/K_{\mathbb{R}} \times X)]$ the role of $\operatorname{CC}(\tilde{\mathcal{F}})$. We shall argue that

$$(4.25)\quad \tau_*\,(dp)^{-1}((-1)^{\dim\mathfrak{p}_{\mathbb{R}}}[T^*_M(G_{\mathbb{R}}/K_{\mathbb{R}} \times X)] - s d\log f) \;=\; (F_s)_*([T^*_QX])\,,$$

for $s > 0$. Because of (4.2), (4.21), and (4.24), the proof of Theorem 3.7 will be complete once we have established (4.25).

As an operation on cycles, $(dp)^{-1}$ is simply intersection with $G_\mathbb{R}/K_\mathbb{R} \times T^*X$ – in our case, Lemma 4.10 asserts that the intersection is transverse. The operation $\tau$ on the intersected cycle is also particularly simple, since projection from $G_\mathbb{R}/K_\mathbb{R} \times T^*X$ to $T^*X$ induces an isomorphism on the carrier of the cycle. In short, Lemma 4.10 implies (4.25) up to sign.

Let us recall the relevant sign conventions. We choose a specific orientation of $\mathfrak{p}_\mathbb{R} \cong G_\mathbb{R}/K_\mathbb{R}$; the particular choice will not matter. This orients the product of $G_\mathbb{R}/K_\mathbb{R}$ with the complex manifold $X$. Our convention for orienting conormal bundles gives meaning to the cycle $[T^*_M(G_\mathbb{R}/K_\mathbb{R} \times X)]$; the convention [SV3, (2.3)] depends on a choice of orientation of the base manifold $G_\mathbb{R}/K_\mathbb{R} \times X$, which we have made. The operation of subtracting $s d\log f$ defines a diffeomorphism of the ambient manifold $T^*(G_\mathbb{R}/K_\mathbb{R} \times X)$. Via this diffeomorphism, $T^*_M(G_\mathbb{R}/K_\mathbb{R} \times X) - s d\log f$ inherits an orientation from that of $T^*(G_\mathbb{R}/K_\mathbb{R} \times X)$, and

$$(4.26)\qquad \begin{aligned} &(-1)^{\dim \mathfrak{p}_\mathbb{R}}[T^*_M(G_\mathbb{R}/K_\mathbb{R} \times X)] - s d\log f \\ &\quad = (-1)^{\dim \mathfrak{p}_\mathbb{R}}[T^*_M(G_\mathbb{R}/K_\mathbb{R} \times X) - s d\log f]\,. \end{aligned}$$

The ambient manifold $T^*(G_\mathbb{R}/K_\mathbb{R} \times X)$ is canonically oriented by the convention for orienting cotangent bundles – space coordinates first, then the corresponding fiber coordinates in the same order; see [SV3, p. 456]. Thus every oriented submanifold becomes co-oriented by the rule

$$(4.27)\qquad \begin{aligned} &(\text{orientation of the submanifold}) \wedge (\text{co-orientation of the submanifold}) \\ &\quad = (\text{orientation of the ambient manifold})\,, \end{aligned}$$

in symbolic notation. Since $T^*_M(G_\mathbb{R}/K_\mathbb{R} \times X) - s d\log f$ intersects $G_\mathbb{R}/K_\mathbb{R} \times T^*X$ transversely along $N_s$, the normal bundle of $N_s$ in $G_\mathbb{R}/K_\mathbb{R} \times T^*X$ is canonically isomorphic to the normal bundle of $T^*_M(G_\mathbb{R}/K_\mathbb{R} \times X) - s d\log f$ in $T^*(G_\mathbb{R}/K_\mathbb{R} \times X)$ along $N_s$. Thus $N_s$ inherits a co-orientation and, by (4.27), an orientation. At this point the orientation of $G_\mathbb{R}/K_\mathbb{R}$ comes in for the second time. We had remarked earlier that $\tau$ is a diffeomorphism on $N_s$, so $\tau(N_s)$ carries a definite orientation, which gives meaning to the cycle $[\tau(N_s)]$. Except for the factor $(-1)^{\dim \mathfrak{p}_\mathbb{R}}$, this is the cycle on the left-hand side of (4.25):

$$(4.28)\quad \tau_*\,(dp)^{-1}((-1)^{\dim \mathfrak{p}_\mathbb{R}}[T^*_M(G_\mathbb{R}/K_\mathbb{R} \times X)] - s d\log f) \;=\; (-1)^{\dim \mathfrak{p}_\mathbb{R}}[\tau(N_s)]\,.$$

Indeed, the procedure (4.26–4.28) for choosing signs at each step precisely follows the prescription in [SV3].

At this stage, we have two orientations on the connected manifold $\tau(N_s) = F_s(T^*_Q X)$, which we must compare. We shall do so at points of the zero section of $T^*_Q X$, which are fixed by $F_s$. For $x \in Q$,

$$(4.29)\quad N_s \cap \tau^{-1}(x,0) \;=\; (eK_{\mathbb{R}}, x, 0, 0) \;\in\; G_{\mathbb{R}}/K_{\mathbb{R}} \times X \times T^*_{eK_{\mathbb{R}}} G_{\mathbb{R}}/K_{\mathbb{R}} \times T^*_x X \;\cong\; G_{\mathbb{R}}/K_{\mathbb{R}} \times X \times \mathfrak{p}_{\mathbb{R}} \times \mathfrak{n}_x \,.$$

The tangent space of the ambient manifold at the point of intersection (4.29) is

$$(4.30)\quad \begin{aligned} &T_{(eK_{\mathbb{R}},x,0,0)}(T^*(G_{\mathbb{R}}/K_{\mathbb{R}} \times X)) \\ &\quad = T_{(eK_{\mathbb{R}},x)}(G_{\mathbb{R}}/K_{\mathbb{R}} \times X) \oplus T^*_{(eK_{\mathbb{R}},x)}(G_{\mathbb{R}}/K_{\mathbb{R}} \times X) \\ &\quad \cong \mathfrak{p}_{\mathbb{R}} \oplus \mathfrak{g}/\mathfrak{b}_x \oplus \mathfrak{p}_{\mathbb{R}} \oplus \mathfrak{n}_x \,. \end{aligned}$$

It contains the tangent space of submanifold $T^*_M(G_{\mathbb{R}}/K_{\mathbb{R}} \times X) - s d\log f$. From (4.11) and (4.16), we see

$$(4.31)\quad \begin{aligned} &T^*_M(G_{\mathbb{R}}/K_{\mathbb{R}} \times X) - s d\log f \quad \text{is the totality of points} \\ &(\exp\zeta \cdot K_{\mathbb{R}}, \exp\zeta \cdot x, -\operatorname{Ad}(\exp\zeta)(\operatorname{Re}\eta - s\zeta), \operatorname{Ad}(\exp\zeta)\eta)\,, \\ &\text{with } \zeta \text{ ranging over } \mathfrak{p}_{\mathbb{R}},\ x \text{ over } Q, \text{ and } \eta \text{ over } \mathfrak{n}_x \cap \mathfrak{p}\,. \end{aligned}$$

Differentiating this parametrization, we get a description of the tangent space,

$$(4.32)\quad \begin{aligned} &T_{(eK_{\mathbb{R}},x,0,0)}(T^*_M(G_{\mathbb{R}}/K_{\mathbb{R}} \times X) - s d\log f) \\ &\quad = \{\, (\zeta\,,\, (\zeta + \kappa) + \mathfrak{b}_x\,,\, -\operatorname{Re}\eta + s\zeta\,,\, \eta) \mid \zeta \in \mathfrak{p}_{\mathbb{R}}\,,\, \kappa \in \mathfrak{k}\,,\, \eta \in \mathfrak{n}_x \cap \mathfrak{p} \,\}\,; \end{aligned}$$

here $(\dots) + \mathfrak{b}_x$ denotes the image of $\dots$ in $\mathfrak{g}/\mathfrak{b}_x$. We choose $\mathbb{C}$-linear complements $\mathfrak{v}_x$ of $\mathfrak{k}/\mathfrak{k} \cap \mathfrak{b}_x$ in $\mathfrak{g}/\mathfrak{b}_x$ and $\mathfrak{q}_x$ of $\mathfrak{n}_x \cap \mathfrak{p}$ in $\mathfrak{n}_x$. Then, since $s \neq 0$,

$$(4.33)\quad \begin{aligned} &\mathfrak{p}_{\mathbb{R}} \oplus \mathfrak{v}_x \oplus 0 \oplus \mathfrak{q}_x \quad \text{is a linear complement of} \\ &T_{(eK_{\mathbb{R}},x,0,0)}(T^*_M(G_{\mathbb{R}}/K_{\mathbb{R}} \times X) - s d\log f) \quad \text{in} \\ &T_{(eK_{\mathbb{R}},x,0,0)}(T^*(G_{\mathbb{R}}/K_{\mathbb{R}} \times X)) \cong \mathfrak{p}_{\mathbb{R}} \oplus \mathfrak{g}/\mathfrak{b}_x \oplus \mathfrak{p}_{\mathbb{R}} \oplus \mathfrak{n}_x\,; \end{aligned}$$

cf. (4.30).

An orientation of an orientable, connected manifold is described by orienting its tangent space at one point. In the following, we shall think of an orientation for a real vector space as an equivalence class of frames, i.e., of ordered bases. In the case of a complex vector space $V$, the underlying real vector space $V^{\mathbb{R}}$ gets oriented by the complex structure: if $\{v_1, \dots, v_m\}$ is a $\mathbb{C}$-frame of $V$, the corresponding $\mathbb{R}$-frame $\{v_1, iv_1, \dots, v_m, iv_m\}$ orients $V^{\mathbb{R}}$ independently of the order of the $v_j$. The direct sum of two oriented vector spaces is oriented by combining positively oriented frames of the summands in the given order, and the dual space $V^*$ of an oriented vector space $V$ gets oriented by the frame dual to an oriented frame of $V$. The underlying real vector space $(V^*)^{\mathbb{R}}$ of the dual $V^*$ of a complex vector space is canonically isomorphic to $(V^{\mathbb{R}})^*$. It thus has two orientations, related by the rule

$$(4.34)\quad \text{dual orientation of } (V^{\mathbb{R}})^* \;=\; (-1)^{\dim_{\mathbb{C}} V} \text{complex orientation of } (V^*)^{\mathbb{R}}\,;$$

reason: if $\{v_1^*, \dots, v_m^*\}$ is dual to a $\mathbb{C}$-frame $\{v_1, \dots, v_m\}$ of $V$, the real part of $\langle iv_j^*, iv_j\rangle$ equals -1. As a general rule, we orient the conormal bundle of a complex submanifold as the conormal bundle of the underlying real manifold – see the discussion at the beginning of Section 3. In the case of a submanifold $N$ of an oriented manifold $M$, the orientation of the conormal bundle $T_N^*M$ is given by identifying

$$T_{(n,0)}(T_N^*M) \;\cong\; T_nN \oplus (T_N^*M)_n \subset T_nM \oplus T_n^*M \;\cong\; T_{(n,0)}(T^*M)\,,$$

and proceeding as follows: we choose an orientation of $T_nN$, orient the quotient $T_nM/T_nN \cong ((T_N^*M)_n)^*$ consistently with the orientation of $T_nM \cong T_nN \oplus ((T_N^*M)_n)^*$. We then put the dual orientation on $(T_N^*M)_n$; then

(4.35)
$$\text{orientation of } (T_N^*M)_n \;=\; (-1)^{\operatorname{codim}_{\mathbb{R}}(N,M)}\,\text{orientation of } T_nN \oplus (T_N^*M)_n,$$

in accordance with the convention [SV3, (2.3)]. Note that (4.34) and (4.35) are consistent when we identify $V^*$ with the conormal bundle $T_0^*V$.

Our choice of orientation for $G_{\mathbb{R}}/K_{\mathbb{R}}$ orients $\mathfrak{p}_{\mathbb{R}} \cong T_{eK_{\mathbb{R}}}(G_{\mathbb{R}}/K_{\mathbb{R}})$. Recall the choice of complements $\mathfrak{v}_x, \mathfrak{q}_x$ in (4.33). We regard $\mathfrak{n}_x$ , $\mathfrak{g}/\mathfrak{b}_x$ , $\mathfrak{k}/\mathfrak{k}\cap\mathfrak{b}_x$ , $\mathfrak{n}_x \cap \mathfrak{p}$ , $\mathfrak{v}_x$ , $\mathfrak{q}_x$ as real vector spaces, oriented by their complex structure. Then, via the isomorphism

$$T_{(eK_{\mathbb{R}},x,0,0)}(T_M^*(G_{\mathbb{R}}/K_{\mathbb{R}} \times X) - s d\log f) \;\cong\; \mathfrak{p}_{\mathbb{R}} \oplus \mathfrak{k}/\mathfrak{k}\cap\mathfrak{b}_x \oplus \mathfrak{n}_x \cap \mathfrak{p}\,, \tag{4.36}$$

which is implicit in (4.32),

$$\begin{aligned} &\text{orientation of } \; T_{(eK_{\mathbb{R}},x,0,0)}(T_M^*(G_{\mathbb{R}}/K_{\mathbb{R}} \times X) - s d\log f) \\ &\quad = \; (-1)^{\operatorname{codim}_{\mathbb{C}}(Q,X)}\,\text{product orientation of } \; \mathfrak{p}_{\mathbb{R}} \oplus \mathfrak{k}/\mathfrak{k}\cap\mathfrak{b}_x \oplus \mathfrak{n}_x \cap \mathfrak{p}\,. \end{aligned} \tag{4.37}$$

What matters here is the discrepancy (4.34) between the two orientations of $\mathfrak{n}_x \cap \mathfrak{p} \cong \mathfrak{v}_x^*$; the sign in (4.35) does not show up because the real codimension of $M$ is even. The co-orientation of $T_M^*(G_{\mathbb{R}}/K_{\mathbb{R}} \times X)$ orients the linear complement (4.33), and

$$\begin{aligned} &\text{co-orientation on } \; T_{(eK_{\mathbb{R}},x,0,0)}(T_M^*(G_{\mathbb{R}}/K_{\mathbb{R}} \times X) - s d\log f) \\ &\quad = \; (-1)^{\dim_{\mathbb{R}}\mathfrak{p}_{\mathbb{R}} + \dim_{\mathbb{C}} Q}\,\text{product orientation of } \; \mathfrak{p}_{\mathbb{R}} \oplus \mathfrak{v}_x \oplus 0 \oplus \mathfrak{q}_x\,. \end{aligned} \tag{4.38}$$

The sign reflects the sign in (4.37), the sign in

$$\begin{aligned} &(-1)^n \text{complex orientation of } \; \mathfrak{n}_x \cong T_x^*X^{\mathbb{R}} \\ &\quad = \; \text{dual orientation of } \; \mathfrak{g}_x/\mathfrak{b}_x \cong T_xX^{\mathbb{R}}\,, \end{aligned} \tag{4.39}$$

and one other sign. This third sign becomes transparent when we replace the space (4.32) by $0 \oplus \mathfrak{k}/\mathfrak{k}\cap\mathfrak{b}_x \oplus \mathfrak{p}_{\mathbb{R}} \oplus \mathfrak{n}_x \cap \mathfrak{p}$, to which it is congruent modulo the complement (4.33). Then

$$\begin{aligned} &\text{orientation of } \; (0 \oplus \mathfrak{k}/\mathfrak{k}\cap\mathfrak{b}_x \oplus \mathfrak{p}_{\mathbb{R}} \oplus \mathfrak{n}_x \cap \mathfrak{p}) \,\oplus\, (\mathfrak{p}_{\mathbb{R}} \oplus \mathfrak{v}_x \oplus 0 \oplus \mathfrak{q}_x) \\ &\quad = \; (-1)^{\dim_{\mathbb{R}}\mathfrak{p}_{\mathbb{R}}}\,\text{orientation of } \mathfrak{p}_{\mathbb{R}} \oplus \mathfrak{g}/\mathfrak{b}_x \oplus \mathfrak{p}_{\mathbb{R}} \oplus \mathfrak{n}_x\,; \end{aligned} \tag{4.40}$$

in effect, we need to move the second copy of $\mathfrak{p}_{\mathbb{R}}$ past the first copy – the other necessary moves involve even dimensional spaces, and thus do not contribute a sign.

Since $(T^*_M(G_{\mathbb{R}}/K_{\mathbb{R}} \times X) - sd\log f$ intersects $G_{\mathbb{R}}/K_{\mathbb{R}} \times T^*X$ transversely along $N_s$, and since the complement (4.33) lies in $T_{(eK_{\mathbb{R}},x,0)}(G_{\mathbb{R}}/K_{\mathbb{R}} \times T^*X)$,

$$\text{(4.41)} \qquad \begin{aligned} &\mathfrak{p}_{\mathbb{R}} \oplus \mathfrak{v}_x \oplus \mathfrak{q}_x \ \text{ is a linear complement of } \ T_{(eK_{\mathbb{R}},x,0)}N_s \ \text{ in} \\ &\qquad T_{(eK_{\mathbb{R}},x,0)}(G_{\mathbb{R}}/K_{\mathbb{R}} \times T^*X) \cong \ \mathfrak{p}_{\mathbb{R}} \oplus \mathfrak{g}/\mathfrak{b}_x \oplus \mathfrak{n}_x \,. \end{aligned}$$

Because of (4.38),

$$\text{(4.42)} \qquad \begin{aligned} &\text{co-orientation on } \ T_{(eK_{\mathbb{R}},x,0)}N_s \\ &\qquad = \ (-1)^{\dim_{\mathbb{R}} \mathfrak{p}_{\mathbb{R}} + \dim_{\mathbb{C}} Q} \text{product orientation of } \ \mathfrak{p}_{\mathbb{R}} \oplus \mathfrak{v}_x \oplus \mathfrak{q}_x \,. \end{aligned}$$

This co-orientation and the orientation of $G_{\mathbb{R}}/K_{\mathbb{R}} \times T^*X$ induce an orientation on the tangent space of $N_s$ at $(eK_{\mathbb{R}}, x, 0)$. We get a description of this tangent space by differentiating (4.19), at $\eta = 0$:

$$\text{(4.43)} \qquad T_{(eK_{\mathbb{R}},x,0)}N_s \ = \{ \, (s^{-1}\operatorname{Re}\eta \, , \, (s^{-1}\operatorname{Re}\eta + \kappa) + \ fb_x \, , \, \eta) \mid \kappa \in \mathfrak{k} \, , \, \eta \in \mathfrak{n}_x \cap \mathfrak{p} \, \} \, .$$

Modulo the complement (4.41), we get the congruence

$$\text{(4.44)} \qquad T_{(eK_{\mathbb{R}},x,0)}N_s \ \cong \ 0 \oplus \mathfrak{k}/\mathfrak{k} \cap \mathfrak{b}_x \oplus \mathfrak{n}_x \cap \mathfrak{p} \ \subset \ \mathfrak{p}_{\mathbb{R}} \oplus \mathfrak{g}/\mathfrak{b}_x \oplus \mathfrak{n}_x \,.$$

In view of (4.42), (4.39), and the even dimensionality of $\mathfrak{k}/\mathfrak{k} \cap \mathfrak{b}_x$ and $\mathfrak{n}_x \cap \mathfrak{p}$, the orientation of $T_{(eK_{\mathbb{R}},x,0)}N_s$ and the product orientation of $0 \oplus \mathfrak{k}/\mathfrak{k} \cap \mathfrak{b}_x \oplus \mathfrak{n}_x \cap \mathfrak{p}$ are related by the factor $(-1)^{\dim_{\mathbb{R}} \mathfrak{p}_{\mathbb{R}} + \operatorname{codim}_{\mathbb{C}}(Q,X)}$. The projection $\tau$ simply drops the first factor, hence

$$\text{(4.45)} \qquad \begin{aligned} &\text{orientation on } \ T_{(x,0)}\tau(N_s) \\ &\qquad = \ (-1)^{\dim_{\mathbb{R}} \mathfrak{p}_{\mathbb{R}} + \operatorname{codim}_{\mathbb{C}}(Q,X)} \text{product orientation of } \ \mathfrak{k}/\mathfrak{k} \cap \mathfrak{b}_x \oplus \mathfrak{n}_x \cap \mathfrak{p} \,. \end{aligned}$$

The isomorphisms (4.15) induce

$$\text{(4.46)} \qquad T_{(x,0)}(T^*_QX) \ \cong \ \mathfrak{k}/\mathfrak{k} \cap \mathfrak{b}_x \oplus \mathfrak{n}_x \cap \mathfrak{p} \ \subset \ \mathfrak{g}/\mathfrak{b}_x \oplus \mathfrak{n}_x \ \cong \ T_{(x,0)}(T^*X) \,.$$

Our rule for orienting conormal bundles implies

$$\text{(4.47)} \qquad \begin{aligned} &\text{orientation of } \ T_{(x,0)}(T^*_QX) \\ &\qquad = \ (-1)^{\operatorname{codim}_{\mathbb{C}}(Q,X)} \text{ product orientation of } \ \mathfrak{k}/\mathfrak{k} \cap \mathfrak{b}_x \oplus \mathfrak{n}_x \cap \mathfrak{p} \,; \end{aligned}$$

to see this, we argue as in (4.37). The induced map $(F_s)_*$ reduces to the identity along the zero section and $T^*_QX$ is connected; hence

$$\text{(4.48)} \qquad [\tau(N_s)] \ = \ (-1)^{\dim_{\mathbb{R}} \mathfrak{p}_{\mathbb{R}}} (F_s)_*([T^*_QX])$$

by (4.45–4.47). This, in combination with (4.28), proves (4.25).

## 5. Descent to the nilpotent cone

The nilpotent cone $\mathcal{N} \subset \mathfrak{g}$ is a finite union of $G$-orbits. Each orbit $\mathcal{O} \subset \mathcal{N}$ carries a $G$-invariant, nondegenerate, complex algebraic 2-form $\sigma_{\mathcal{O}}$. This form becomes canonical when we identify $\mathcal{O}$ with the corresponding coadjoint orbit via $B : \mathfrak{g} \cong \mathfrak{g}^*$ as in Section 3. In particular, the $G$-orbits $\mathcal{O} \subset \mathcal{N}$ have even complex dimension. We stratify $\mathcal{N}$ by dimension, with

$$\mathcal{N}_k \;=\; \bigcup \{\, \mathcal{O} \mid \dim_{\mathbb{C}} \mathcal{O} = 2k \,\}\;, , \qquad \tilde{\mathcal{N}}_k \;=\; \textstyle\bigcup_{\ell \leq k} \mathcal{N}_\ell\,. \tag{5.1}$$

Then $\tilde{\mathcal{N}}_k$ is closed, and $\mathcal{N}_k$ is open in $\tilde{\mathcal{N}}_k$. For $\zeta \in \mathcal{N}$, the Springer fiber $\mu^{-1}(\zeta)$ is complete, connected, equidimensional, of complex dimension $n - \frac{1}{2}\dim_{\mathbb{C}}(G\cdot\zeta)$ [Spal]. The moment map is $G$-equivariant, so the preceding statement implies:

$$\begin{gathered}\text{for each } G\text{-orbit } \mathcal{O}_k \subset \mathcal{N}_k, \quad \mu : \mu^{-1}(\mathcal{O}_k) \to \mathcal{O}_k \quad \text{is a} \\ G\text{-equivariant fibration, whose typical fiber } \mu^{-1}(\zeta) \text{ is a connected,} \\ \text{complex projective, equidimensional variety of dimension } n-k. \end{gathered} \tag{5.2}$$

We shall use the stratification (5.1) to filter the spaces $T^*_K X$ and $T^*_{G_\mathbb{R}} X$.

In the case of $T^*_K X$, the $k^{\text{th}}$ filtrant is the closed complex algebraic subspace $T^*_K X \cap \mu^{-1}\tilde{\mathcal{N}}_k$. Then

(5.3)

a) $\ldots \subset\; T^*_K X \cap \mu^{-1}\tilde{\mathcal{N}}_k \;\subset\; T^*_K X \cap \mu^{-1}\tilde{\mathcal{N}}_{k+1} \;\subset\; \ldots\;$ ;

b) $\bigcup_k T^*_K X \cap \mu^{-1}\tilde{\mathcal{N}}_k \;=\; T^*_K X\,$;

c) $T^*_K X \cap \mu^{-1}\mathcal{N}_k$ is open in $T^*_K X \cap \mu^{-1}\tilde{\mathcal{N}}_k\,$;

d) the boundary of $T^*_K X \cap \mu^{-1}\mathcal{N}_k$ is contained in $T^*_K X \cap \mu^{-1}\tilde{\mathcal{N}}_{k-1}\,$;

e) $T^*_K X \cap \mu^{-1}\mathcal{N}_k$ is equidimensional, of complex dimension $n$.

The first four statements follow formally from the corresponding statements about the $\tilde{\mathcal{N}}_k$. To see e), note that, $T^*_K X \cap \mu^{-1}\mathcal{N}_k$ fibers over $\mathcal{N}_k \cap \mathfrak{p}$, in other words, over the union of $K$-orbits in $\mathcal{N}_k \cap \mathfrak{p}$. Because of (3.5), the fiber at each point is the full Springer fiber $\mu^{-1}(x)$. Because of (5.2), this fiber is equidimensional of complex dimension $n-k$. The same reasoning gives properties a)–e) for the filtration of $T^*_{G_\mathbb{R}} X$ by the $T^*_{G_\mathbb{R}} X \cap \mu^{-1}\tilde{\mathcal{N}}_k$. In this case, of course, the $T^*_{G_\mathbb{R}} X \cap \mu^{-1}\tilde{\mathcal{N}}_k$ are real algebraic subvarieties, and $T^*_{G_\mathbb{R}} X \cap \mu^{-1}\mathcal{N}_k$ is equidimensional of real dimension $2n$.

Recall that the characteristic cycle maps (3.3) take values in $\mathrm{H}^{\inf}_{2n}(T^*_K X, \mathbb{Z})$ and $\mathrm{H}^{\inf}_{2n}(T^*_{G_\mathbb{R}} X, \mathbb{Z})$. The filtration (5.3) and its real analogue induce filtrations

$$\begin{aligned} \mathrm{H}^{\inf}_{2n}(T^*_K X, \mathbb{Z})_k &= \mathrm{Im}\,\{\mathrm{H}^{\inf}_{2n}(T^*_K X \cap \mu^{-1}\tilde{\mathcal{N}}_k, \mathbb{Z}) \to \mathrm{H}^{\inf}_{2n}(T^*_K X, \mathbb{Z})\}, \\ \mathrm{H}^{\inf}_{2n}(T^*_{G_\mathbb{R}} X, \mathbb{Z})_k &= \mathrm{Im}\,\{\mathrm{H}^{\inf}_{2n}(T^*_{G_\mathbb{R}} X \cap \mu^{-1}\tilde{\mathcal{N}}_k, \mathbb{Z}) \to \mathrm{H}^{\inf}_{2n}(T^*_{G_\mathbb{R}} X, \mathbb{Z})\}\,. \end{aligned} \tag{5.4}$$

Because of (5.3c,d), we get a well-defined restriction map

$$\mathrm{H}^{\inf}_{2n}(T^*_K X, \mathbb{Z})_k \;\longrightarrow\; \mathrm{gr}_k\, \mathrm{H}^{\inf}_{2n}(T^*_K X, \mathbb{Z}) \;\longrightarrow\; \mathrm{H}^{\inf}_{2n}(T^*_K X \cap \mu^{-1}\mathcal{N}_k, \mathbb{Z})\,, \tag{5.5a}$$

and analogously,

$$\text{(5.5b)}\quad \mathrm{H}^{\inf}_{2n}(T^*_{G_\mathbb{R}}X,\mathbb{Z})_k \ \longrightarrow\ \mathrm{gr}_k\ \mathrm{H}^{\inf}_{2n}(T^*_{G_\mathbb{R}}X,\mathbb{Z}) \ \longrightarrow\ \mathrm{H}^{\inf}_{2n}(T^*_{G_\mathbb{R}}X\cap\mu^{-1}\mathcal{N}_k,\mathbb{Z})\,.$$

Integration over the fibers of $\mu$ defines maps from the groups on the right in (5.5a,b) to groups of cycles in $\mathcal{N}_k\cap\mathfrak{p}$ and $\mathcal{N}_k\cap i\mathfrak{g}_\mathbb{R}$, as we shall explain next.

The dual space $\mathfrak{h}^*$ of the universal Cartan contains the universal weight lattice $\Lambda$, i.e., the lattice of differentials of algebraic characters of the universal Cartan group $H$. For $\lambda\in\Lambda$, the character $e^\lambda : H\to\mathbb{C}^*$ determines a $G$-equivariant algebraic line bundle $\mathcal{L}_\lambda\to X$ whose fiber at $x\in X$ is the complex line, on which $B_x$ acts via $e^\lambda$. The map $\lambda\mapsto c_1(\mathcal{L}_\lambda)$ (= first Chern class of $\mathcal{L}_\lambda$) defines a homomorphism from $\Lambda$ to $\mathrm{H}^2(X,\mathbb{Z})$; hence

$$\text{(5.6)}\qquad c_1\,:\,\mathfrak{h}^*\ \longrightarrow\ \mathrm{H}^2(X,\mathbb{C})\,.$$

For simplicity, we write $e^\lambda$ for the cohomology class $1+c_1(\lambda)+\frac{c_1(\lambda)^2}{2}+\dots$ , which we view as a class on $T^*X$. The usual cap product pairing applies also in the setting of homology with infinite support. Thus we can take the cap product of any class in $\mathrm{H}^{\inf}_{2n}(T^*_KX\cap\mu^{-1}\mathcal{N}_k,\mathbb{Z})$ against the component of $e^\lambda$ in degree $2n-2k$; this produces a homology class in degree $2k$, which can then be pushed forward to a class in $\mathrm{H}^{\inf}_{2k}(\mathcal{N}_k\cap\mathfrak{p},\mathbb{C})$. This operation, and its $G_\mathbb{R}$-analogue, is our process of integration over the fibers. In general, the definition of cap product involves certain sign conventions. We do not need to spell them out, since integration over the fiber can be described more concretely in our particular situation.

The moment map defines a fibration $\mu : T^*_KX\cap\mu^{-1}\mathcal{N}_k\to\mathcal{N}_k\cap\mathfrak{p}$. As top dimensional cycle, each $C_{2n}\in\mathrm{H}^{\inf}_{2n}(T^*_KX\cap\mu^{-1}\mathcal{N}_k,\mathbb{Z})$ can be regarded, locally with respect to the base of the fibration, as a product of a top dimensional cycle $C_{2n-2k}(\zeta)$ in the (compact) fiber with a top dimensional cycle in $\mathcal{N}_k\cap\mathfrak{p}$ – here we are using the complex structure of the $\mu^{-1}(\zeta)$ to orient the fiber component of the cycle; the even dimensionality of the fiber makes the order of the product irrelevant. We integrate the component of $e^\lambda$ in degree $2n-2k$ over $C_{2n-2k}(\zeta)$; the resulting function $m(\zeta)$ is locally constant and independent of the particular product decomposition. We multiply the base component of $C_{2n}$ by the multiplicity $m(\zeta)$. This gives us a well-defined class

$$\text{(5.7)}\qquad \int_{C_{2n}} e^\lambda\ \in\ \mathrm{H}^{\inf}_{2k}(\mathcal{N}_k\cap\mathfrak{p},\mathbb{C})\,,$$

the pushforward of the degree $2k$ component of the cap product[1] $e^\lambda\cap C_{2n}$. Analogously, we define integration over the fiber, as a map from $\mathrm{H}^{\inf}_{2n}(T^*_{G_\mathbb{R}}X\,\cap$

[1] To see that the geometric operation of integration over the fiber does agree with cap product, followed by pushforward, one can use the formalism in [SV3], for example: when cycles are viewed as local cohomology classes along their supports, cap product gets converted into cup product; this makes available the full apparatus of cohomology.

$\mu^{-1}\mathcal{N}_k, \mathbb{Z})$ to $\mathrm{H}^{\inf}_{2k}(\mathcal{N}_k \cap i\mathfrak{g}_{\mathbb{R}}, \mathbb{C})$. In general, integration and cap product agree only up to sign. Specifically, if $M$ is a compact oriented $m$-manifold, cap product $\omega \cap [M]$ with the fundamental class, for $\omega \in \mathrm{H}^m(M, \mathbb{C})$ is $(-1)^m$ times the integral of a de Rham representative of $\omega$ over $M$. In our case, the sign is doubly irrelevant: not only are the Springer fibers even dimensional, but if there were a change of sign – as there may be in the more general situation of a semisimple symmetric space – it will occur twice in the statement of Theorem 5.10.

To each $C \in \mathrm{H}^{\inf}_{2n}(T^*_K X \cap \mu^{-1}\mathcal{N}, \mathbb{Z})$, we assign the degree $k = k(C)$, the least integer $k$ such that $C \in \mathrm{H}^{\inf}_{2n}(T^*_K X \cap \mu^{-1}\tilde{\mathcal{N}}_k, \mathbb{Z})$. We then restrict $C$ to $T^*_K X \cap \mu^{-1}\mathcal{N}_k$ and perform the operation (5.7). This gives us

$$\text{(5.8a)} \qquad \begin{aligned} (\mathrm{gr}\,\mu_*)_\lambda \ &: \ \mathrm{H}^{\inf}_{2n}(T^*_K X \cap \mu^{-1}\mathcal{N}, \mathbb{Z}) \ \longrightarrow \ \textstyle\bigoplus_k \mathrm{H}^{\inf}_{2k}(\mathcal{N}_k \cap \mathfrak{p}, \mathbb{C})\,, \\ (\mathrm{gr}\,\mu_*)_\lambda(C) \ &= \ \textstyle\int_{C^0} e^\lambda\,, \quad C^0 = C|_{T^*_K X \cap \mu^{-1}\mathcal{N}_k}\,, \quad k = k(C)\,. \end{aligned}$$

All of this makes sense equally on the $G_{\mathbb{R}}$-side:

$$\text{(5.8b)} \qquad \begin{aligned} (\mathrm{gr}\,\mu_*)_\lambda \ &: \ \mathrm{H}^{\inf}_{2n}(T^*_{G_{\mathbb{R}}} X \cap \mu^{-1}\mathcal{N}, \mathbb{Z}) \ \longrightarrow \ \textstyle\bigoplus_k \mathrm{H}^{\inf}_{2k}(\mathcal{N}_k \cap i\mathfrak{g}_{\mathbb{R}}, \mathbb{C})\,, \\ (\mathrm{gr}\,\mu_*)_\lambda(C) \ &= \ \textstyle\int_{C^0} e^\lambda\,, \quad C^0 = C|_{T^*_{G_{\mathbb{R}}} X \cap \mu^{-1}\mathcal{N}_k}\,, \quad k = k(C)\,. \end{aligned}$$

Note that the two maps $(\mathrm{gr}\,\mu_*)_\lambda$ are not homomorphisms, since the definitions involve going to the leading terms in the graded groups $\mathrm{gr}\,\mathrm{H}^{\inf}_{2n}(T^*_K X \cap \mu^{-1}\mathcal{N}, \mathbb{Z})$ and $\mathrm{gr}\,\mathrm{H}^{\inf}_{2n}(T^*_{G_{\mathbb{R}}} X \cap \mu^{-1}\mathcal{N}, \mathbb{Z})$.

The family $F_s : T^*X \to T^*X$ defined in (3.6) induces a family of bianalytic maps $f_s$ on the nilpotent cone,

$$\text{(5.9)} \qquad \begin{aligned} f_s \ &: \ \mathcal{N} \to \mathcal{N}\,, \qquad s \in \mathbb{R}_{>0}\,, \\ f_s(\eta) \ &= \ \mathrm{Ad}(\exp(s^{-1}\,\mathrm{Re}\,\eta))\eta\,, \qquad (\,\eta \in \mathcal{N}\,)\,. \end{aligned}$$

Because of (4.20), $\mu \circ F_s = f_s \circ \mu$, so $F_s$ does lie over $f_s$.

5.10 Theorem. *For $c \in \mathrm{H}^{\inf}_{2k}(\mathcal{N}_k \cap \mathfrak{p}, \mathbb{C})$, the limit of cycles $\lim_{s\to 0^+}(f_s)_*(c)$ exists as cycle in $\mathcal{N}_k$ and has support in $\mathcal{N}_k \cap i\mathfrak{g}_{\mathbb{R}}$. The resulting homomorphism*

$$\phi \ : \ \textstyle\bigoplus_k \mathrm{H}^{\inf}_{2k}(\mathcal{N}_k \cap \mathfrak{p}, \mathbb{C}) \ \longrightarrow \ \textstyle\bigoplus_k \mathrm{H}^{\inf}_{2k}(\mathcal{N}_K \cap i\mathfrak{g}_{\mathbb{R}}, \mathbb{C})\,,$$

$\phi(c) = \lim_{s\to 0^+}(f_s)_*(c)$ *for* $c \in \mathrm{H}^{\inf}_{2k}(\mathcal{N}_k \cap \mathfrak{p}, \mathbb{C})$, *makes the following diagram commutative*:

$$\mathrm{H}^{\inf}_{2n}(T^*_K X, \mathbb{Z}) @> \ \ \Phi \ \ >> \mathrm{H}^{\inf}_{2n}(T^*_{G_{\mathbb{R}}} X, \mathbb{Z})$$

$$@V(\mathrm{gr}\,\mu_*)_\lambda VV @VV(\mathrm{gr}\,\mu_*)_\lambda V$$

$$\textstyle\bigoplus_k \mathrm{H}^{\inf}_{2k}(\mathcal{N}_k \cap \mathfrak{p}, \mathbb{C}) @> \ \ \phi \ \ >> \textstyle\bigoplus_k \mathrm{H}^{\inf}_{2k}(\mathcal{N}_k \cap i\mathfrak{g}_{\mathbb{R}}, \mathbb{C})\,.$$

*Proof.* Recall the definition of the analytic-geometric category $\mathcal{C}$ in Section 4. To see that the limit of cycles exists, we shall argue that the family of cycles $(f_s)_*(c)$ is a $\mathcal{C}$-family – in other words, that the union of the supports $f_s(|c|)$ is a $\mathcal{C}$-subset of $\mathcal{N}_k \times \mathbb{R}$. We proceed as in the proof of Theorem 3.7. For dimension reasons, $|c|$ is a union of $K$-orbits, hence algebraic in $\mathcal{N}_k$. Thus

(5.11)
$$\{\,(\mathrm{Ad}(g)\eta, gK_{\mathbb{R}}, s) \mid gK_{\mathbb{R}} \in G_{\mathbb{R}}/K_{\mathbb{R}}\,,\ \eta \in |c|\,,\ s \in \mathbb{R}\,\} \ \subset\ \mathcal{N}_k \times G_{\mathbb{R}}/K_{\mathbb{R}} \times \mathbb{R}$$

is a real algebraic subset, and consequently a $\mathcal{C}$-subset of $\mathcal{N}_k \times \overline{G_{\mathbb{R}}/K_{\mathbb{R}}} \times \mathbb{R}$. Since $\overline{\mathfrak{p}_{\mathbb{R}}} \cong \overline{G_{\mathbb{R}}/K_{\mathbb{R}}}$ as $\mathcal{C}$-sets,

$$\{\,(\mathrm{Ad}(\exp\zeta)\eta, \zeta, s) \mid \zeta \in \mathfrak{p}_{\mathbb{R}}, \eta \in |c|, s \in \mathbb{R}\,\} \ \subset\ \mathcal{N}_k \times\ \overline{\mathfrak{p}_{\mathbb{R}}} \times \mathbb{R} \tag{5.12}$$

is a $\mathcal{C}$-subset, as is its intersection with $\{\, s\zeta = \mathrm{Re}\,\eta\,\}$. This intersection is the support of the family of cycles $\{(f_s)_*(c)\}$. In view of [SV3, §3], the limit $\lim_{s\to 0^+}(f_s)_*(c)$ exists as cycle in $\mathcal{N}_k$. More specifically,

$$\lim_{s\to 0^+}(f_s)_*(c) \in \mathrm{H}^{\mathrm{inf}}_{2k}(\mathcal{N}_K \cap i\mathfrak{g}_{\mathbb{R}}, \mathbb{C})\,, \tag{5.13}$$

as follows from the argument below the statement of Theorem 3.7.

We regard $J = [0,\infty]$ as a closed subinterval of the one-point compactification $\mathbb{R}\cup\{\infty\}$ of $\mathbb{R}$. Let us consider a particular $C \in \mathrm{H}^{\mathrm{inf}}_{2n}(T^*_K X, \mathbb{Z})$. The cycles $(F_s)_*(C)$ constitute a family of $2n$-cycles in $T^*X$, parametrized by $I = (0,\infty)$, in the sense of [SV3, §3]. In other words, there exists a $(2n+1)$-chain $C_J$ in $J \times T^*X$ such that

$$\begin{aligned} &\text{a)}\ \ |C_J| \ =\ \text{closure of }\ \{\,(s, F_s(\zeta)) \mid s \in I\,,\ \zeta \in |C|\,\}\,;\\ &\text{b)}\ \ C_I \ =\ C_J|_{(I\times\mu^{-1}(\tilde{\mathcal{N}}_k))}\ \ \text{is a } 2n\text{-cycle in } I \times T^*X\,;\\ &\text{c)}\ \ C_s \ =\ (F_s)_*(C) \quad \text{for } 0 < s < \infty\,;\\ &\text{d)}\ \ \partial C_J \ =\ \{\infty\}\times C \ -\ \{0\}\times \lim_{s\to 0^+} C_s\,. \end{aligned} \tag{5.14}$$

Here $C_s$ denotes the specialization of the family $C_I$ at $s$; i.e.,

$$\{s\}\times C_s \ =\ \partial C_I|_{(0,s]\times T^*X} \quad (\text{boundary in } (0,s]\times T^*X)\,. \tag{5.15}$$

Let $k = k(C)$ be the least integer $k$ such that $C \in \mathrm{H}^{\mathrm{inf}}_{2n}(T^*_K X, \mathbb{Z})_k$ – equivalently, the least integer $k$ such that $C$ is supported in $\mu^{-1}(\tilde{\mathcal{N}}_k)$. Since $\mathcal{N}_k$ is $G$-invariant, this implies that all the cycles $(F_s)_*(C)$ are supported in $\mu^{-1}(\tilde{\mathcal{N}}_k)$. Thus we can regard $C_J$ as a family of cycles in $\mu^{-1}(\tilde{\mathcal{N}}_k)$.

Recall the definition of $\Phi(C)$ as the limit of $(F_s)_*(C)$ as $s \to 0^+$. The preceding remarks imply, in particular, that $\Phi(C)$ is supported in $\mu^{-1}(\tilde{\mathcal{N}}_k)$, hence $k(\Phi(C)) \le k(C)$. In fact,

$$k(\Phi(C)) \ =\ k(C)\,. \tag{5.16}$$

To see this, we observe that $\Phi$ is derived from the functor $\gamma$ of (2.7), which has an inverse, the functor $\Gamma$ of [MUV] with an appropriate shift in degree. Concretely, $\Gamma$ is defined in terms of the diagram (2.8), but with $K$ and $K/K_{\mathbb{R}}$ taking the places of $G_{\mathbb{R}}$ and $G_{\mathbb{R}}/K_{\mathbb{R}}$ and with stars instead of shrieks. Just as $\gamma$ determines $\Phi$, the functor $\Gamma$ determines

$$\Psi \ : \ \mathrm{H}^{\mathrm{inf}}_{2n}(T^*_{G_{\mathbb{R}}}X,\mathbb{Z}) \ \longrightarrow \ \mathrm{H}^{\mathrm{inf}}_{2n}(T^*_K X,\mathbb{Z})\,. \tag{5.17}$$

Theorem (3.7) and its proof carry over to this situation word for word, with one exception: in pinning down the sign in (3.7), it was convenient to use the complex structure of the $K$-orbits. In any case, the same reasoning that gave us the inequality $k(\Phi(C)) \leq k(C)$ gives

$$k(\Psi(C)) \ \leq \ k(C) \qquad (\, C \in \mathrm{H}^{\mathrm{inf}}_{2n}(T^*_{G_{\mathbb{R}}}X,\mathbb{Z})\,)\,. \tag{5.18}$$

Since $\Gamma \circ \gamma$ is equivalent to the identity on $\mathrm{D}_K(X)_{-\lambda}$, up to sign, $\Psi \circ \Phi$ is the identity on the CC-image of $\mathrm{D}_K(X)_{-\lambda}$, up to the same sign. We claim:

$$\mathrm{CC}(\mathrm{D}_K(X)_{-\lambda}) \ = \ \mathrm{H}^{\mathrm{inf}}_{2n}(T^*_K X,\mathbb{Z})\,, \tag{5.19}$$

for integral $\lambda \in \mathfrak{h}^*$, and in particular for $\lambda = \rho$. Assuming this for the moment, we see that $\Psi \circ \Phi$ – which does not depend on the particular choice of $\lambda$ – is the identity, up to sign, on all of $\mathrm{H}^{\mathrm{inf}}_{2n}(T^*_K X,\mathbb{Z})$. Thus (5.16) follows from (5.18) and the earlier inequality $k(\Phi(C)) \leq k(C)$.

We need to establish (5.19). If $\lambda = \rho$, as we may assume, $\mathrm{D}_K(X)_{-\rho} = \mathrm{D}_K(X)$ contains all the direct images $Rj_*\mathbb{C}_S$ of constant sheaves $\mathbb{C}_S$ on $K$-orbits $S \subset X$. The characteristic cycle $\mathrm{CC}(Rj_*\mathbb{C}_S)$ differs from $[T^*_S X]$ by a linear combination of cycles $[T^*_{S'}X]$ with $S' \subset \partial S$. Thus, arguing by induction on the dimension of $S$, we find $[T^*_S X] \in \mathrm{CC}(\mathrm{D}_K(X))$ for all $K$-orbits $S$. These cycles $[T^*_S X]$ span $\mathrm{H}^{\mathrm{inf}}_{2n}(T^*_K X,\mathbb{Z})$, so (5.19) follows.

At this point, we have established the first assertion of Theorem 5.10, and we have shown that $\Phi$ preserves the integer $k(C)$ which enters the definition of integration over the fibers. This operation can be applied to the family $C_I$. The family of diffeomorphisms $F_s$ can be used to trivialize $C_I$; in other words,

$$(0,\infty] \times |C| \; @> \; \sim \; >> \; |C_{(0,\infty]}| \,, \tag{5.20}$$

so that $C_{(0,\infty]}$ becomes the product of the interval $(0,\infty]$ with the cycle $C$. Contrary to appearance, $\infty$ is a generic point of the family, since $F_s$ extends smoothly across $s = \infty$ – recall (3.6). Integrating $e^\lambda$ over the $\mu$-fibers of this family, we obtain a family of $2k$-cycles $c_{(0,\infty]}$ in $\mathcal{N}_k$. It, too, is a product family whose general member is $c_s = (f_s)_*(c)$, with $c = \int_C e^\lambda$. As in (5.14), $c_{(0,\infty]}$ is the restriction to $(0,\infty] \times \mathcal{N}_k$ of a $(2k+1)$-chain $c_J$ in $J \times \mathcal{N}_k$, such that

$$\partial c_J \;=\; \{\infty\} \times c \;-\; \{0\} \times \lim_{s\to 0^+} c_s \,. \tag{5.21}$$

The two families $C_J$, $c_J$ cease to be product families at $s = 0$. To see that $c_s = \int C_s e^\lambda$ even at $s = 0$, we appeal to the formalism of cap product, as follows.

Let $\omega \in \mathrm{H}^{2n-2k}(T^*X, \mathbb{C})$ be the component in degree $2n-2k$ of $e^\lambda$, pulled back from $X$ to $T^*X$. We shall take the cap product of a cochain representative of $\omega$ with the chain $C_J$ – or more precisely, with

$$C'_J \;=\; \text{restriction of } C_J \text{ to } J \times \mu^{-1}(\mathcal{N}_k) \,. \tag{5.22}$$

This can be carried out in several ways, but perhaps most transparently in the simplicial setting. Thus we triangulate, compatibly, the spaces $J \times \mu^{-1}(\mathcal{N}_k)$, $|C'_J|$, $\{0\} \times \mu^{-1}(\mathcal{N}_k)$, and $\{\infty\} \times \mu^{-1}(\mathcal{N}_k)$; we can do so because $|C_J|$ is a $\mathcal{C}$-set. We choose a cochain representative of $\omega$ and, for simplicity, denote it by the same letter. We can then take the cap products $\omega \cap C'_J$, $\omega \cap \partial C'_J$. Since $\omega$ is closed, they are related by the formula

$$\partial(\omega \cap C'_J) \;=\; \omega \cap \partial C'_J \,; \tag{5.23}$$

for the sign conventions, we follow [Span]. Because of the triviality of the family $C_I$, $\mu_*(\omega \cap C'_I)$ is the family in $\mathcal{N}_k$ whose general member is obtained by integrating $\omega$ over the $\mu$-fiber – in other words, the general member is $c_s$; hence

$$\mu_*(\omega \cap C'_I) \;=\; c_I \,. \tag{5.24}$$

All the top dimensional simplices of $C_J$ lie in $C_I$, so (5.24) remains valid with $J$ in place of $I$:

$$\mu_*(\omega \cap C'_J) \;=\; c_J \,. \tag{5.25}$$

Because of (5.14d), (5.23), and (5.25),

$$\mu_*(\omega \cap (\{\infty\} \times C')) - \mu_*(\omega \cap (\{0\} \times \Phi(C)')) = \{\infty\} \times c - \{0\} \times \phi(c) \,, \tag{5.26}$$

where $C'$ and $\Phi(C)'$ denote the restrictions of the two cycles to $\mu^{-1}(\mathcal{N}_k)$.

From our definition of integration over the fiber, we see that $\mu_*(\omega \cap (\{0\} \times \Phi(C)'))$ equals $\{0\} \times \int_{\Phi(C)} e^\lambda$. Hence, because of (5.26),

$$\phi(c) \;=\; \int_{\Phi(C)} e^\lambda\,. \tag{5.27}$$

Taken together, (5.16) and (5.27) give the commutativity of the diagram in the statement of Theorem 5.10.

## 6. Nilpotent orbits

Let $\mathcal{O}_k$ be a $G$-orbit in $\mathcal{N}_k$. Then $\mathcal{O}_k \cap \mathfrak{p}$ is a union of finitely many $K$-orbits, and similarly, $\mathcal{O}_k \cap i\mathfrak{g}_\mathbb{R}$ is a union of finitely many $G_\mathbb{R}$-orbits. These intersections are Lagrangian – in the case of $\mathfrak{p}$, with respect to the holomorphic symplectic structure $\sigma_{\mathcal{O}_k}$ on $\mathcal{O}_k$, and in the case of $i\mathfrak{g}_\mathbb{R}$, with respect to $\operatorname{Re}\sigma_{\mathcal{O}_k}$ [KR]; here, as in Section 5, we identify $\mathfrak{g} \cong \mathfrak{g}^*$ and $\mathcal{N} \cong \mathcal{N}^*$ by means of the bilinear form $B$. In particular, each $K$-orbit in $\mathcal{N}_k \cap \mathfrak{p}$ is a complex manifold of dimension $k$. We use the complex structure to orient the orbits. This allows us to regard them as $K$-invariant, top dimensional cycles in $\mathcal{N}_k \cap \mathfrak{p}$. In fact, these cycles constitute a basis of $\mathrm{H}^{\inf}_{2k}(\mathcal{N}_k \cap \mathfrak{p}, \mathbb{C})^K$, the $K$-invariant part of the top dimensional homology of $\mathcal{N}_k \cap \mathfrak{p}$:

$$\mathrm{H}^{\inf}_{2k}(\mathcal{N}_k \cap \mathfrak{p}, \mathbb{C})^K \;=\; \{\textstyle\sum a_j\,[\mathcal{O}_{\mathfrak{p},j}] \mid a_j \in \mathbb{C}\,,\ \dim_\mathbb{C} \mathcal{O}_{\mathfrak{p},j} = k\,\}\,, \tag{6.1a}$$

with $\mathcal{O}_{\mathfrak{p},j}$ enumerating the $K$-orbits in $\mathcal{N} \cap \mathfrak{p}$. To see this, we note that the connected components of the $K$-orbits provide a basis of $\mathrm{H}^{\inf}_{2k}(\mathcal{N}_k \cap \mathfrak{p}, \mathbb{C})$; the fundamental cycle $[\mathcal{O}_{\mathfrak{p},j}]$ of a $K$-orbit $\mathcal{O}_{\mathfrak{p},j}$ is the sum of the fundamental cycles of the components of $\mathcal{O}_{\mathfrak{p},j}$. Analogously,

$$\mathrm{H}^{\inf}_{2k}(\mathcal{N}_k \cap i\mathfrak{g}_\mathbb{R}, \mathbb{C})^{G_\mathbb{R}} \;=\; \{\textstyle\sum b_j\,[\mathcal{O}_{\mathfrak{g}_\mathbb{R},j}] \mid b_j \in \mathbb{C}\,,\ \dim_\mathbb{R} \mathcal{O}_{\mathfrak{g}_\mathbb{R},j} = 2k\,\} \tag{6.1b}$$

when we enumerate the $G_\mathbb{R}$-orbits in $\mathcal{N} \cap i\mathfrak{g}_\mathbb{R}$ as $\mathcal{O}_{\mathfrak{g}_\mathbb{R},j}$. Each of them lies in a $G$-orbit $\mathcal{O}_k$, from which it inherits the symplectic form $\frac{1}{2\pi i}\sigma_{\mathcal{O}_k}$ – note that the restriction of $\sigma_{\mathcal{O}_k}$ is purely imaginary on $\mathcal{O}_k \cap i\mathfrak{g}_\mathbb{R}$. We use the symplectic structure to orient the $\mathcal{O}_{\mathfrak{g}_\mathbb{R},j}$, to give meaning to the cycles $[\mathcal{O}_{\mathfrak{g}_\mathbb{R},j}]$.

Sekiguchi [Se] and Kostant (unpublished) have described a bijective correspondence between the $K$-orbits in $\mathcal{N} \cap \mathfrak{p}$ on one hand, and the $G_\mathbb{R}$-orbits in $\mathcal{N} \cap i\mathfrak{g}_\mathbb{R}$ on the other. Orbits that correspond to each other lie in the same $G$-orbit, and thus have the same dimension. Recall the definition of

$$\phi \;:\; \textstyle\bigoplus_k \mathrm{H}^{\inf}_{2k}(\mathcal{N}_k \cap \mathfrak{p}, \mathbb{C}) \;\longrightarrow\; \bigoplus_k \mathrm{H}^{\inf}_{2k}(\mathcal{N}_k \cap i\mathfrak{g}_\mathbb{R}, \mathbb{C}) \tag{6.2}$$

in the statement of Theorem 5.10, which was defined in terms of the family of diffeomorphisms $f_s : \mathfrak{g} \to \mathfrak{g}\,,\quad f_s(\eta) = \operatorname{Ad}(\exp(s^{-1}\operatorname{Re}\eta))(\eta)$.

6.3 THEOREM. *The map $\phi$ is an isomorphism. It sends $K$-invariant cycles to $G_{\mathbb{R}}$-invariant cycles. On the invariant part of the homology, $\phi$ coincides with the Kostant-Sekiguchi correspondence via the identifications* (6.1). *Concretely, if $[\mathcal{O}_{\mathfrak{p}}] \in \mathrm{H}^{\mathrm{inf}}_{2k}(\mathcal{N}_k \cap \mathfrak{p}, \mathbb{C})^K$ is the fundamental class of a $K$-orbit $\mathcal{O}_{\mathfrak{p}}$, oriented by its complex structure, then the family of cycles $(f_s)_*[\mathcal{O}_{\mathfrak{p}}]$ has a limit as $s \to 0^+$, and this limit is the fundamental class $[\mathcal{O}_{\mathfrak{g}_{\mathbb{R}}}] \in \mathrm{H}^{\mathrm{inf}}_{2k}(\mathcal{N}_k \cap i\mathfrak{g}_{\mathbb{R}}, \mathbb{C})^{G_{\mathbb{R}}}$ of the Sekiguchi image $\mathcal{O}_{\mathfrak{g}_{\mathbb{R}}}$ of $\mathcal{O}_{\mathfrak{p}}$, oriented by its symplectic form.*

Sekiguchi describes the correspondence between the two types of orbits by reduction to the special case of $G_{\mathbb{R}} = \mathrm{SL}(2,\mathbb{R})$ , $K_{\mathbb{R}} = \mathrm{SO}(2,\mathbb{R})$. Let

$$j \ : \ \mathfrak{sl}(2,\mathbb{C}) \longrightarrow \mathfrak{g} \tag{6.4}$$

be a homomorphism, defined over $\mathbb{R}$ with respect to the real forms $\mathfrak{sl}(2,\mathbb{R})$ and $\mathfrak{g}_{\mathbb{R}}$, and equivariant with respect to the Cartan involutions – in the case of $\mathfrak{sl}(2,\mathbb{C})$, the Cartan involution corresponding to the maximal compact subgroup $\mathrm{SO}(2,\mathbb{R})$ of $\mathrm{SL}(2,\mathbb{R})$. According to Kostant-Rallis [KR],

$$\begin{gathered}\text{every } \zeta \in \mathcal{N} \cap \mathfrak{p} \text{ is } K\text{-conjugate to the } j\text{-image of } \begin{pmatrix} 1 & i \\ i & -1 \end{pmatrix} \\ \text{for some homomorphism } j \text{ as in (6.4).}\end{gathered} \tag{6.5}$$

On the other hand, the Jacobson-Morozov theorem for the Lie algebra $\mathfrak{g}_{\mathbb{R}}$ implies

$$\begin{gathered}\text{every } \eta \in \mathcal{N} \cap i\mathfrak{g}_{\mathbb{R}} \text{ is } G_{\mathbb{R}}\text{-conjugate to the } j\text{-image of } \begin{pmatrix} 0 & i \\ 0 & 0 \end{pmatrix} \\ \text{for some homomorphism } j \text{ as in (6.4);}\end{gathered} \tag{6.6}$$

see, for example, [Ko]. Sekiguchi shows that the $K$-orbit of $\zeta$, and similarly the $G_{\mathbb{R}}$-orbit of $\eta$, determines the homomorphism $j$ up to $K_{\mathbb{R}}$-conjugacy [Se]. Thus

$$K\text{-orbit of } \ j\begin{pmatrix} 1 & i \\ i & -1 \end{pmatrix} \ \longleftrightarrow \ G_{\mathbb{R}}\text{-orbit of } \ j\begin{pmatrix} 0 & i \\ 0 & 0 \end{pmatrix} \ , \tag{6.7}$$

for every homomorphism $j$ as in (6.4), sets up a well-defined correspondence between $K$-orbits in $\mathcal{N} \cap \mathfrak{p}$ and $G_{\mathbb{R}}$-orbits in $\mathcal{N} \cap i\mathfrak{g}_{\mathbb{R}}$.

Every $\zeta \in \mathcal{N} \cap \mathfrak{p}$ is $K$-conjugate to its negative, but $\eta \in \mathcal{N} \cap i\mathfrak{g}_{\mathbb{R}}$ need not be $G_{\mathbb{R}}$-conjugate to $-\eta$. It would be equally natural to let the $K$-orbit of $\zeta$ in (6.5) correspond to the $G_{\mathbb{R}}$-orbit of $-\eta$ in (6.6). From our point of view, the microlocalization of the functor $\gamma$ dictates the choice of $\eta$ over $-\eta$: since

$$\exp\left(s^{-1}\,\mathrm{Re}\begin{pmatrix} 1 & i \\ i & -1 \end{pmatrix}\right)\begin{pmatrix} 1 & i \\ i & -1 \end{pmatrix} \ = \ \begin{pmatrix} 1 & ie^{2s^{-1}} \\ ie^{-2s^{-1}} & -1 \end{pmatrix} \ , \tag{6.8}$$

the definition of $\phi$ forces

$$K \cdot \begin{pmatrix} 1 & i \\ i & -1 \end{pmatrix} \ \longleftrightarrow \ G_{\mathbb{R}} \cdot \begin{pmatrix} 0 & i \\ 0 & 0 \end{pmatrix} \tag{6.9}$$

in the case of $\mathrm{SL}(2,\mathbb{R})$, and correspondingly (6.7) in general.

The proof of Theorem 6.3 is lengthy. It uses methods completely different from those in the rest of the paper. Here we shall reduce the assertion of the theorem to certain technical statements about nilpotent orbits, which are proved in [SV5].

We begin with a simplification of the problem: it suffices to consider the case of a connected semisimple group $G_\mathbb{R}$. Indeed, when $G_\mathbb{R}$ is connected, then so are $K_\mathbb{R}$ and $K$. In that case all homology classes in $\mathrm{H}^{\inf}_{2k}(\mathcal{N}_k \cap \mathfrak{p}, \mathbb{C})$ are $K$-invariant, and similarly all classes in $\mathrm{H}^{\inf}_{2k}(\mathcal{N}_k \cap i\mathfrak{g}_\mathbb{R}, \mathbb{C})$ are $G_\mathbb{R}$-invariant. Now

$$K/K^0 \;\cong\; K_\mathbb{R}/K^0_\mathbb{R} \;\cong\; G_\mathbb{R}/G^0_\mathbb{R}\,, \tag{6.10}$$

so the invariance conditions in (6.1) are equivalent to invariance under the component group $K_\mathbb{R}/K^0_\mathbb{R}$. But $\phi$ commutes with the action of $K_\mathbb{R}$. Thus we may as well assume that $G_\mathbb{R}$ is connected. All nilpotents lie in the derived algebra $[\mathfrak{g},\mathfrak{g}]$. Hence, without changing the problem, we can replace $G_\mathbb{R}$ by its quotient by the connected component of the center. For emphasis,

$$G_\mathbb{R} \;=\; G^0_\mathbb{R}\,, \qquad \mathfrak{g} \;=\; [\mathfrak{g},\mathfrak{g}]\,, \tag{6.11}$$

as will be assumed from now on.

The support of the family of cycles $(f_s)_*[\mathcal{O}_\mathfrak{p}]$, $s > 0$, is contained in a single $G$-orbit $\mathcal{O}$ – the $G$-orbit which contains the $K$-orbit $\mathcal{O}_\mathfrak{p}$. We suppose $\mathcal{O} \neq \{0\}$, since otherwise there is nothing to prove. According to Theorem 5.10, the limit of the family exists as cycle in the union of $G$-orbits having the same dimension as $\mathcal{O}$. This union is disjoint; hence

$$\text{the limit } \lim_{s\to 0^+} (f_s)_*[\mathcal{O}_\mathfrak{p}] \text{ exists as cycle supported on } \mathcal{O}\cap i\mathfrak{g}_\mathbb{R}\;. \tag{6.12}$$

Since $\mathcal{O}_\mathfrak{p}$ has the same dimension as $\mathcal{O}\cap i\mathfrak{g}_\mathbb{R}$, the limit is necessarily an integral linear combination of fundamental classes of the finitely many $G_\mathbb{R}$-orbits in $\mathcal{O}\cap i\mathfrak{g}_\mathbb{R}$. We enumerate these orbits as $\mathcal{O}_{\mathfrak{g}_\mathbb{R},j}$; then

$$\lim_{s\to 0^+} (f_s)_*[\mathcal{O}_\mathfrak{p}] \;=\; \sum b_j\,[\mathcal{O}_{\mathfrak{g}_\mathbb{R},j}]\,, \qquad b_j \in \mathbb{Z}\,. \tag{6.13}$$

Note that both $\mathcal{O}_\mathfrak{p}$ and the $\mathcal{O}_{\mathfrak{g}_\mathbb{R},j}$ are connected because of (6.11). With this notation, Theorem 6.3 amounts to a description of the $b_j$,

$$b_j \;=\; \begin{cases} 1 & \text{if } \mathcal{O}_{\mathfrak{g}_\mathbb{R},j} \text{ is the Sekiguchi image of } \mathcal{O}_\mathfrak{p}\,, \\ 0 & \text{otherwise}\,. \end{cases} \tag{6.14}$$

That is what we must prove.

The multiplicities $b_j$ can be expressed as intersection multiplicities of the cycles $(f_s)_*[\mathcal{O}_\mathfrak{p}]$ with normal slices to the $\mathcal{O}_{\mathfrak{g}_\mathbb{R},j}$. To do this, we fix a particular $\nu \in \mathcal{O}_{\mathfrak{g}_\mathbb{R},j}$ and choose a linear complement $\mathfrak{q}_\mathbb{R}$ to the kernel of $\operatorname{ad}\nu$ in $\mathfrak{g}_\mathbb{R}$. For $a > 0$ sufficiently small,

$$N(\nu,a) \;=\; \{\,\operatorname{Ad}\exp(i\eta)(\nu) \mid \nu \in \mathfrak{q}_\mathbb{R}\,,\ \|\eta\| < a\,\} \tag{6.15}$$

is a real analytic submanifold of $\mathcal{O}$ which meets $\mathcal{O}_{\mathfrak{g}_{\mathbb{R}},j}$ only at $\nu$, and the intersection at $\nu$ is transverse; in other words, $N(\nu,a)$ is a "normal slice" to $\mathcal{O}_{\mathfrak{g}_{\mathbb{R}},j}$, at $\nu$, in $\mathcal{O}$. Then

$$\text{(6.16)}\qquad \begin{array}{c}\text{for generic } \nu \in \mathcal{O}_{\mathfrak{g}_{\mathbb{R}},j}\text{, with } a>0 \text{ sufficiently small, and}\\ s \text{ small in relation to } a\text{, the cycle } (f_s)_*[\mathcal{O}_{\mathfrak{p}}] \text{ intersects } N(\nu,a)\\ \text{transversely, with total intersection multiplicity } b_j.\end{array}$$

In this statement, "generic" means "on an open, dense $\mathcal{C}$-set," where $\mathcal{C}$ refers to the analytic-geometric category introduced and used in Section 4; intersections are to be counted with the same sign convention that makes $\mathcal{O}_{\mathfrak{g}_{\mathbb{R}},j}$ meet $N(\nu,a)$ with multiplicity $+1$. We refer to [SV3, §3] for the notion of limit of a family of cycles, as we have earlier.

It looks prohibitively difficult to compute the intersection multiplicities at a generic point $\nu$ directly. Instead, we shall establish a slightly stronger statement at certain (conceivably) nongeneric points[2], from which we then deduce the needed information about generic points. In preparation for the argument, we introduce the compact real form

$$\text{(6.17a)}\qquad \begin{array}{l}\mathfrak{u}_{\mathbb{R}} \;=\; \mathfrak{k}_{\mathbb{R}} \oplus i\mathfrak{p}_{\mathbb{R}} \;=\; \{\, \zeta \in \mathfrak{g} \mid \theta\bar{\zeta} = \zeta \,\}\\ (\,\bar{\zeta} = \text{complex conjugate of } \zeta \text{ with respect to } \mathfrak{g}_{\mathbb{R}}\,)\end{array}$$

in $\mathfrak{g}$, and the maximal compact subgroup

$$\text{(6.17b)}\qquad U_{\mathbb{R}} \;=\; \text{connected subgroup of } G \text{ with Lie algebra } \mathfrak{u}_{\mathbb{R}}$$

of $G$. Since $G$ is connected,

$$\text{(6.17c)}\qquad G_{\mathbb{R}} \cap U_{\mathbb{R}} \;=\; K_{\mathbb{R}} \;=\; K \cap U_{\mathbb{R}}\,;$$

cf. [He], for example. In previous sections we had chosen a particular Ad-invariant bilinear form $B$ on $\mathfrak{g}$. Now that $\mathfrak{g}$ is semisimple by assumption, we let $B$ denote the Killing form, normalized as follows. By Jacobson-Morozov, any $\zeta \in \mathcal{O}$ can be embedded in an essentially unique $\mathfrak{sl}_2$-triple. In other words, there exist $\tau$, $\zeta_-$ in $\mathfrak{g}$ such that

$$\text{(6.18)}\qquad [\tau,\zeta] = 2\zeta\,, \quad [\tau,\zeta_-] = -2\zeta_-\,, \quad [\zeta,\zeta_-] = \tau\,.$$

$\tau$ is unique up to conjugacy by the centralizer $G_\zeta$ of $\zeta$, and $\zeta_-$ becomes unique once $\tau$ has been chosen. In particular, the orbit $\mathcal{O}$ determines $\tau$ up to $G$-conjugacy. Since $\tau$ is nonzero (recall: we had assumed $\mathcal{O} \neq \{0\}$), semisimple, with integral eigenvalues, we can normalize $B$ by requiring

$$\text{(6.19)}\qquad B(\tau,\tau) \;=\; 2\,.$$

[2] Indeed, we believe that in our particular situation, every point is generic in the sense of (6.16).

Thus $B$ restricts to the linear span of $\zeta, \zeta_-, \tau$ as the trace form of $\mathfrak{sl}(2,\mathbb{C})$, to which this linear span is isomorphic. In terms of $B$, we define

$$(6.20) \qquad (\zeta_1, \zeta_2) \;=\; -B(\zeta_1, \theta\bar{\zeta}_2) \qquad \zeta_1, \zeta_2 \in \mathfrak{g}\,.$$

This is a (positive definite) $U_{\mathbb{R}}$-invariant inner product on $\mathfrak{g}$.

We introduce a moment map for the action of $G$ on $\mathcal{O}$, following Ness [N]. From an intrinsic point of view, we should think of the moment map as taking values in $i\mathfrak{u}_{\mathbb{R}}{}^*$. It will be more convenient, however, to identify $i\mathfrak{u}_{\mathbb{R}} \cong i\mathfrak{u}_{\mathbb{R}}{}^*$, and to define

$$(6.21a) \qquad m \;:\; \mathcal{O} \longrightarrow i\mathfrak{u}_{\mathbb{R}}$$

implicitly, by the equation

$$(6.21b) \qquad 2\,\mathrm{Re}\,(m(\zeta), \eta) \;=\; \frac{1}{\|\zeta\|^2}\left(\frac{d}{dt}\,\|\,\mathrm{Ad}\exp(t\eta)\zeta\|^2\right)\Big|_{t=0}\,.$$

As $\eta$ runs over $\mathfrak{g}$ in this equation, $m(\zeta)$ becomes determined as vector in $\mathfrak{g}$. But the inner product is $U_{\mathbb{R}}$-invariant; hence $m(\zeta)$ does lie in $i\mathfrak{u}_{\mathbb{R}}$. The $U_{\mathbb{R}}$-invariance also implies

$$(6.22) \qquad m(\mathrm{Ad}(u)\zeta) \;=\; \mathrm{Ad}(u)(m(\zeta)) \qquad (\,u \in U_{\mathbb{R}}\,)\,;$$

i.e., the map $m$ is $U_{\mathbb{R}}$-equivariant. To get an explicit formula for $m(\zeta)$, we calculate:

$$\begin{aligned}\left(\frac{d}{dt}\,\|\,\mathrm{Ad}\exp(t\eta)\zeta\|^2\right)\Big|_{t=0} &=\; 2\,\mathrm{Re}\,([\eta,\zeta],\zeta)\\ &=\; -2\,\mathrm{Re}\;B([\eta,\zeta],\theta\bar{\zeta}) \;=\; -2\,\mathrm{Re}\;B(\eta,[\zeta,\theta\bar{\zeta}]) \;=\; 2\,\mathrm{Re}\;B(\eta,\theta\overline{[\zeta,\theta\bar{\zeta}]})\\ &=\; -2\,\mathrm{Re}\,(\eta,[\zeta,\theta\bar{\zeta}]) \;=\; -2\,\mathrm{Re}\,([\zeta,\theta\bar{\zeta}],\eta)\,,\end{aligned}$$

for every test vector $\eta \in \mathfrak{g}$; hence

$$(6.23) \qquad m(\zeta) \;=\; -\,\frac{[\zeta,\theta\bar{\zeta}]}{\|\zeta\|^2}\,.$$

The moment map descends to the image of the orbit $\mathcal{O}$ in the projectivized Lie algebra $\mathbb{P}(\mathfrak{g})$. Viewed as map from $\mathbb{P}(\mathfrak{g})$ to $i\mathfrak{u}_{\mathbb{R}} \cong i\mathfrak{u}_{\mathbb{R}}{}^*$, $m$ coincides with the moment map, in the sense of symplectic geometry, relative to a $U_{\mathbb{R}}$-invariant symplectic structure on $\mathbb{P}(\mathfrak{g})$ [N].

Because of the equivariance (6.22), the square length $\|m(\zeta)\|^2$ is invariant under the $U_{\mathbb{R}}$-action. It is also invariant under scaling by any nonzero complex number. Scaling by positive real numbers plays a special role, since it preserves $G_{\mathbb{R}}$-orbits. We thus regard $\mathcal{O}$ as manifold with $U_{\mathbb{R}} \times \mathbb{R}^+$-action, with $\mathbb{R}^+$, the multiplicative group of positive real numbers, acting by scaling.

6.24 LEMMA. *A point* $\zeta \in \mathcal{O}$ *is a critical points of the function* $\zeta \mapsto \|m(\zeta)\|^2$ *if and only if* $[\zeta, \theta\bar{\zeta}]$ *can be rescaled so that* $\zeta, [\zeta, \theta\bar{\zeta}], \theta\bar{\zeta}$ *becomes an* $\mathfrak{sl}_2$*-triple – in other words, if and only if there exists* $a \in \mathbb{R}$, $a \neq 0$, *such that*

$$[[\zeta, \theta\bar{\zeta}], \zeta] \;=\; a\,\zeta \qquad \text{and} \qquad [[\zeta, \theta\bar{\zeta}], \theta\bar{\zeta}] \;=\; -a\,\theta\bar{\zeta}\,.$$

*The set of critical points is not empty and consists of a single* $U_{\mathbb{R}} \times \mathbb{R}^+$*-orbit. The function* $\|m\|^2$ *on* $\mathcal{O}$ *assumes its minimum values exactly on the critical set. Every* $K$*-orbit in* $\mathcal{O} \cap \mathfrak{p}$ *and every* $G_{\mathbb{R}}$*-orbit in* $\mathcal{O} \cap i\mathfrak{g}_{\mathbb{R}}$ *meets the critical set along exactly one* $K_{\mathbb{R}} \times \mathbb{R}^+$*-orbit.*

This follows from a general property of the moment map [N, Theorem 6.1]; for details, see [SV5]. We shall be able to analyze the intersection of the family of cycles $(f_s)_*[\mathcal{O}_{\mathfrak{p}}]$ with an appropriately chosen normal slice $N(\nu, a)$ for points $\nu \in \mathcal{O}_{G_{\mathbb{R}},j}$ which are critical of the function $\|m\|^2$. As a first step, we show:

6.25 LEMMA. *Let* $\mathcal{O}_{\mathfrak{p}}$ *be a* $K$*-orbit in* $\mathcal{O} \cap \mathfrak{p}$, *and* $\zeta \in \mathcal{O}_{\mathfrak{p}}$ *a critical point for* $\|m\|^2$. *Then* $f_s(\zeta)$ *lies in the critical set for every* $s \in \mathbb{R}_{>0}$. *Moreover, the limit*

$$\lim_{s \to 0^+} \frac{f_s(\zeta)}{\|f_s(\zeta)\|}$$

*exists and lies in the Sekiguchi image* $\mathcal{O}_{\mathfrak{g}_{\mathbb{R}}}$ *of* $\mathcal{O}_{\mathfrak{p}}$.

*Proof.* For $t > 0$, $f_s(t\zeta) = t f_{st^{-1}}(\zeta)$, so we are free to rescale $\zeta$ by a positive real number. Also, $\zeta$ is $K_{\mathbb{R}}$-conjugate to $-\zeta$, because $[[\zeta, \theta\bar{\zeta}], \zeta]$ is a nonzero real multiple of $\zeta$ and $\exp(it[\zeta, \theta\bar{\zeta}]) \in K_{\mathbb{R}}$ for $t \in \mathbb{R}$. Since $f_s$ is $K_{\mathbb{R}}$-equivariant, we can now rescale $\zeta$ by any nonzero real number. In other words, we may assume that $a = -2$, in the notation of the previous lemma. In that case, the linear map $j : \mathfrak{sl}(2, \mathbb{C}) \to \mathfrak{g}$, defined by

$$(6.26)\qquad \begin{aligned} j\begin{pmatrix} 1 & 0 \\ 0 & -1 \end{pmatrix} \;&=\; \zeta - \theta\bar{\zeta}\,, \qquad & j\begin{pmatrix} 0 & 1 \\ 0 & 0 \end{pmatrix} \;&=\; \frac{1}{2i}(\zeta + \theta\bar{\zeta} + [\zeta, \theta\bar{\zeta}])\,, \\ & & j\begin{pmatrix} 0 & 0 \\ 1 & 0 \end{pmatrix} \;&=\; \frac{1}{2i}(\zeta + \theta\bar{\zeta} - [\zeta, \theta\bar{\zeta}])\,, \end{aligned}$$

satisfies the conditions on $j$ in (6.4): it is a homomorphism, defined over $\mathbb{R}$ with respect to the real form $\mathfrak{sl}(2, \mathbb{R}) \subset \mathfrak{sl}(2, \mathbb{C})$, and equivariant with respect to the Cartan involution corresponding to the maximal compact subgroup $\mathrm{SO}(2, \mathbb{R})$ of $\mathrm{SL}(2, \mathbb{R})$. In this way we can reduce the problem to a computation in $\mathfrak{sl}_2$, which we have already done; see (6.8).

For the moment, we keep fixed a particular $\zeta \in \mathcal{O}_{\mathfrak{p}}$. Since $m(f_s(\zeta)) \in i\mathfrak{u}_{\mathbb{R}} = i\mathfrak{k}_{\mathbb{R}} \oplus \mathfrak{p}_{\mathbb{R}}$, we can write

$$(6.27)\qquad \begin{gathered} m(f_s(\zeta)) \;=\; m_1(s, \zeta) + m_2(s, \zeta) + m_3(s, \zeta)\,, \quad \text{with} \\ m_1(s, \zeta) \in \mathbb{R} \cdot \mathrm{Re}\;\zeta\,, \quad m_2(s, \zeta) \in \mathfrak{p}_{\mathbb{R}} \cap (\mathrm{Re}\;\zeta)^{\perp}\,, \quad m_3(s, \zeta) \in i\mathfrak{k}_{\mathbb{R}}\,. \end{gathered}$$

Experimental evidence suggests that $\|m(f_s(\zeta))\|^2$ is decreasing for $s > 0$. The parametric curve $f_s(\zeta)$ would then move away from the critical set of $\|m\|^2$ as $s$ approaches 0. We do not how to prove this; however, the following suffices for our purposes.

6.28 Lemma. *For* $s > 0\,,\quad \|m_1(s,\zeta)\|^2 \;+\; \|m_3(s,\zeta)\|^2 \;\geq\; \|m(\zeta)\|^2\,.$

This lemma plays the crucial role in the proof of the next one; both are established in [SV5].

6.29 Lemma. *Let $\nu_0 \in \mathcal{O}_{\mathfrak{g}_\mathbb{R},j}$ be a critical point for $\|m\|^2$. Then there exists a normal slice $N(\nu_0, a)$ with the following properties. If $\mathcal{O}_{\mathfrak{g}_\mathbb{R},j}$ is the Sekiguchi image of $\mathcal{O}_\mathfrak{p}$, then the submanifolds $f_s(\mathcal{O}_\mathfrak{p})$ and $N(\nu_0, a)$ of $\mathcal{O}$ meet exactly once, for all sufficiently small values of $s$. The intersection is transverse and has multiplicity $+1$, relative to the sign convention which makes $\mathcal{O}_{\mathfrak{g}_\mathbb{R},j}$ meet $N(\nu_0, a)$ with multiplicity $+1$ at $\nu_0$. On the other hand, $f_s(\mathcal{O}_\mathfrak{p}) \cap N(\nu_0, a) = \emptyset$ if $\mathcal{O}_{\mathfrak{g}_\mathbb{R},j}$ is not the Sekiguchi image of $\mathcal{O}_\mathfrak{p}$, again for all sufficiently small $s$.*

We shall deduce Theorem 6.3 from Lemma 6.29. In effect, one can phrase the criterion (6.16) less restrictively in the case of a subanalytic or $\mathcal{C}$-family of cycles: one may use even nongeneric normal slices to calculate the intersection multiplicity, provided they satisfy two conditions. First, the normal slice is normal not only to the carrier of the limit cycle, but also normal to limit of the carriers; secondly, the carrier of the family is transverse to the particular normal slice, except possibly at $s = 0$. We shall not try to establish the more general statement in full generality, but only in our particular situation.

We fix a point $\nu_0 \in \mathcal{O}_{\mathfrak{g}_\mathbb{R},j}$ which is critical for the function $\|m\|^2$. The normal slice $N(\nu_0, a)$ mentioned in the lemma corresponds to a choice of a linear complement $\mathfrak{q}_\mathbb{R}$ to $\operatorname{Ker}(\operatorname{ad}\nu_0)$ in $\mathfrak{g}_\mathbb{R}$. The same $\mathfrak{q}_\mathbb{R}$ will then be a linear complement also to $\operatorname{Ker}(\operatorname{ad}\nu)$ for every $\nu \in \mathcal{O}_{\mathfrak{g}_\mathbb{R},j}$ close to $\nu_0$. We choose a small open neighborhood $V_0$ of $\nu_0$ and a sufficiently small constant $a > 0$ so that

$$B(a) \times V_0 \longrightarrow \mathcal{O} \qquad (\,B(a) \;=\; \{\,\eta \in \mathfrak{q}_\mathbb{R} \mid \|\eta\| < a\,\}\,)\,, \qquad (\eta, \nu) \;\mapsto\; \operatorname{Ad}(\exp\eta)\nu \tag{6.30}$$

is a bianalytic map onto its image. We shrink $a$ further, if necessary, so that the normal slice $N(\nu_0, a)$ satisfies the conclusion of Lemma 6.29, and so that the map (6.30) extends bianalytically to $\overline{B(a)} \times V_0$, the partial closure of $B(a) \times V_0$ in the $\mathfrak{q}_\mathbb{R}$-directions. We use this bianalytic map and $B(a) \cong N(\nu_0, a)$ to identify

$$N(\nu_0, a) \times V_0 \;\cong\; \text{neighborhood of } \nu_0 \text{ in } \mathcal{O}\,. \tag{6.31a}$$

This identification is consistent with the tautological inclusions $N(\nu_0, a) \subset \mathcal{O}$ and $V_0 \subset \mathcal{O}_{\mathfrak{g}_\mathbb{R},j} \subset \mathcal{O}$. The projection

$$\pi \;:\; N(\nu_0, a) \times V_0 \longrightarrow V_0 \tag{6.31b}$$

retracts the neighborhood (6.31a) of $\nu_0$ in $\mathcal{O}$ to the neighborhood $V_0$ of $\nu_0$ in $\mathcal{O}_{\mathfrak{g}_{\mathbb{R}},j}$. The fibers of $\pi$ are normal slices; i.e.,

$$\pi^{-1}(\nu) \;=\; N(\nu, a) \;=\; \text{normal slice to } \nu\,. \tag{6.31c}$$

By construction, the product structure (6.31) extends to the closure in the fiber directions.

The carrier of the family of cycles $\{(f_s)_*[\mathcal{O}_{\mathfrak{p}}]\}_{s>0}$ is a closed real analytic submanifold of $\mathbb{R}_{>0} \times \mathcal{O}$; it is also a $\mathcal{C}$-set in $\mathbb{R} \times \mathcal{O}$. Here, once again, $\mathcal{C}$ refers to the analytic-geometric category used earlier. We denote the submanifold by $M$. The natural bianalytic map

$$\mathbb{R}_{>0} \times \mathcal{O}_{\mathfrak{p}} \; @> \;\sim\; >> \; M\,, \qquad (s, \zeta) \mapsto f_s(\zeta) \tag{6.32}$$

orients $M$. We observe that

$$\begin{aligned} &M \cap (\mathbb{R}_{>0} \times N(\nu_0, a) \times V_0) \;\text{ is closed in }\; \mathbb{R}_{>0} \times N(\nu_0, a) \times V_0\,, \;\text{ and} \\ &M \cap (\mathbb{R}_{>0} \times \overline{N(\nu_0, a)} \times V_0) \;\text{ is closed in }\; \mathbb{R}_{>0} \times \overline{N(\nu_0, a)} \times V_0\,. \end{aligned} \tag{6.33}$$

The projection $M \cap (\mathbb{R}_{>0} \times N(\nu_0, a) \times V_0) \to V_0$ can be partially compactified to a proper $\mathcal{C}$-map

$$\overline{M} \cap (\{0 \le s \le \infty\} \times \overline{N(\nu_0, a)} \times V_0) \;\longrightarrow\; V_0\,. \tag{6.34}$$

We shall show:

$$\begin{gathered} \text{there exists a dense open } \mathcal{C}\text{-set } V_1 \subset V_0 \text{ such that} \\ \text{the projection} \quad M \cap (\mathbb{R}_{>0} \times N(\nu_0, a) \times V_1) \;\longrightarrow\; V_1 \\ \text{is of maximal rank everywhere in the domain.} \end{gathered} \tag{6.35}$$

Indeed, by [DM, statement D.13], the projection (6.34) can be stratified, so that on each stratum in the domain the projection has constant rank. Note that $\dim M = \dim V_0 + 1$. Thus, by dimension count, each stratum in the domain either maps to a lower dimensional stratum in $V_0$, or the projection has maximal rank on it. Our statement (6.35) follows; it may happen, of course, that $\mathbb{R}_{>0} \times N(\nu_0, a) \times V_1$ does not intersect $M$ at all.

The generic triviality statement [DM, 4.11] for compactifiable $\mathcal{C}$-maps allows us to shrink $V_1$ further, so that

$$\begin{gathered} V_1 \subset V_0 \;\text{ is open,} \quad \nu_0 \in \overline{V_1}\,, \quad \text{and} \\ \overline{M} \cap (\mathbb{R}_{\ge 0} \times N(\nu_0, a) \times V_1) \;\longrightarrow\; V_1 \quad \text{is a product} \end{gathered} \tag{6.36}$$

(product in the $\mathcal{C}$-continuous sense). What can we say about the fiber $F_\nu$ of this product over a $\nu \in V_1$? To begin with,

$$F_\nu \;\text{ is a } \mathcal{C}\text{-curve in }\; \mathbb{R}_{\ge 0} \times N(\nu, a)\,, \quad \text{unless } F_\nu \text{ is empty}\,, \tag{6.37a}$$

for dimension reasons. Secondly,

$$F_\nu \cap (\{0\} \times N(\nu,a)) \subset \{(0,\nu)\} \quad \text{and} \qquad (6.37b)$$
$$F_\nu \cap (\mathbb{R}_{>0} \times N(\nu,a)) \ \subset\ \mathbb{R}_{>0} \times (N(\nu,a) - \{\nu\})\,;$$

here we use (5.13), which holds and is proved on the level of supports. Also, for $s \neq 0$ and $\zeta \in \mathcal{O}_{\mathfrak{p}}$, the real part of $f_s(\zeta)$ equals the real part of $\zeta$, hence is nonzero, which prevents $f_s(\zeta)$ from lying in $i\mathfrak{g}_{\mathbb{R}}$. Lastly, because of (6.35),

$$F_\nu \cap (\mathbb{R}_{>0} \times N(\nu,a)) \quad \text{is a closed submanifold of} \quad \mathbb{R}_{>0} \times N(\nu,a)\,. \qquad (6.37c)$$

It is also a $\mathcal{C}$-subset; hence it only has a finite number of connected components.

We enumerate the finitely many connected components of $F_\nu \cap (\mathbb{R}_{>0} \times N(\nu,a))$ as $F_{\nu,\ell}$ , $1 \leq \ell \leq L$; note that

$$F_\nu \cap (\mathbb{R}_{>0} \times N(\nu,a)) \ =\ F_\nu - \{(0,\nu)\} \qquad (6.38)$$

because of (6.37). *A priori* each of the $F_{\nu,\ell}$ can be compact, have zero, one, or two endpoints at $\nu$, with the remaining ends "at infinity," i.e., tending to the boundary of the normal slice. Note that an end at $\nu$ corresponds to $s = 0$, and an end in $\partial N(\nu,a)$ to a strictly positive value of $s$; cf. (6.37b). When we count the net intersection multiplicity of $(f_s)_*[\mathcal{O}_{\mathfrak{p}}]$ with the normal slice $N(\nu,a)$, the component $F_{\nu,\ell}$ contributes only if it has one end at $\nu$ and the other at infinity. After all, we are computing the multiplicity at $\{0\} \times \nu$ of the boundary of the chain $[F_\nu]$. We claim:

(6.39)
a) among the $F_{\nu,\ell}$, either one or none run from zero to infinity, depending on whether or not $\mathcal{O}_{\mathfrak{g}_{\mathbb{R}},j}$ is the Sekiguchi image of $\mathcal{O}_{\mathfrak{p}}$;
b) if $F_{\nu,\ell}$ does run from zero to infinity, this curve can be continued across $\nu = \nu_0$.

According to Lemma 6.29, over $\nu = \nu_0$ we see either one curve or none, again depending on whether or not $\mathcal{O}_{\mathfrak{g}_{\mathbb{R}},j}$ is the Sekiguchi image of $\mathcal{O}_{\mathfrak{p}}$. In the former situation, the curve runs from zero to the boundary of the normal slice and has the same intersection multiplicity with $(f_s)_*[\mathcal{O}_{\mathfrak{p}}]$, $0 < s \ll 1$, as with $[\mathcal{O}_{\mathfrak{g}_{\mathbb{R}}}]$. Thus (6.39) does imply Theorem 6.3.

The verification of (6.39) involves two processes: extending the curve over $\nu_0$ – if there is one – to nearby points $\nu$, and specializing to $\nu = \nu_0$ those $F_{\nu,\ell}$ which run from zero to infinity. For the former, we note that

$$\text{the projection} \quad M \cap (\mathbb{R}_{>0} \times N(\nu_0,a) \times V_0) \ \longrightarrow\ V_0 \ \text{ has maximal} \qquad (6.40)$$
$$\text{rank along} \ \ F_{\nu_0} - \{(0,\nu_0)\} \ =\ M \cap (\mathbb{R}_{>0} \times N(\nu_0,a) \times \{\nu_0\});$$

this follows from the transversality assertion in Lemma 6.29. We conclude that the curve over $\nu_0$ – if it exists – can be continued smoothly to nearby points $\nu$, at least if we stay away from $s = 0$. These nearby curves over $\nu \in V_1$ must run from $\nu$ to $\partial N(\nu, a)$, or from $\nu$ back to $\nu$, or from one point on $\partial N(\nu, a)$ to another, or be compact – these are the only possibilities for curves over points in $V_1$, as was mentioned earlier. All but the first possibility are ruled out by the local smoothness (6.40) of the family across $F_{\nu_0}$, away from $s = 0$. To summarize, when there is a curve $F_{\nu_0}$ over $\nu_0$, it can be continued to a curve over nearby points $\nu \in V_1$ which runs from $\nu$ to $\partial N(\nu, a)$. There is at most one such curve above nearby points $\nu \in V_1$, since otherwise the local smoothness (6.40) would be contradicted.

The reverse process, of specializing from $\nu \in V_1$ to $\nu_0$, depends on a property of our particular setting that we have not used so far. Let us consider a curve $F_{\nu,\ell}$ over a point $\nu \in V_1$ that runs from $\nu$ to $\partial N(\nu, a)$. It was mentioned already that $\nu$ corresponds to $s = 0$ and the boundary point to some strictly positive $s_\infty = s_\infty(\nu)$. We need to know:

(6.41) $\quad s_\infty(\nu)$ is bounded away from $0$ for $\nu \in V_1$ near $\nu_0$.

Assuming this for the moment, we can specialize $F_{\nu,\ell}$ to a curve $F_{\nu_0}$ which runs from $\nu_0$ to $\partial N(\nu_0, a)$. Together with the conclusion of the preceding paragraph, this establishes (6.39).

At this point, only (6.41) needs to be established. If the assertion were false, there would exist sequences $\{s_n\}$ in $\mathbb{R}_{>0}$ and $\{\zeta_n\}$ in $\mathcal{O}_\mathfrak{p}$ such that $s_n \to 0$ and $\lim_{n\to\infty} f_{s_n}(\zeta_n) \in \partial N(\nu_0, a)$. On the other hand, the existence of $\lim_{n\to\infty} f_{s_n}(\zeta_n)$ with $s_n \to 0$ forces $\lim_{n\to\infty} f_{s_n}(\zeta_n) \in i\mathfrak{g}_\mathbb{R}$, as was argued at the end of Section 3. Finally, we can decrease $a$ further, if necessary, to ensure $\partial N(\nu_0, a) \cap i\mathfrak{g}_\mathbb{R} = \emptyset$.

## 7. Completion of the proof of Theorem 1.4

The commutative squares (2.9), (3.7), and (5.10) can be combined into a single commutative diagram, as follows:

(7.1)
$$\begin{CD}
\{\text{virtual H-C-modules}\}_\lambda @>\sim>> \{\text{virtual } G_\mathbb{R}\text{-representations}\}_\lambda \\
@A\alpha AA @AA\beta A \\
\mathrm{D}_K(X)_{-\lambda} @>\gamma>> \mathrm{D}_{G_\mathbb{R}}(X)_{-\lambda} \\
@V\mathrm{CC}VV @VV\mathrm{CC}V \\
\mathrm{H}^{\mathrm{inf}}_{2n}(T^*_K X, \mathbb{Z}) @>\Phi>> \mathrm{H}^{\mathrm{inf}}_{2n}(T^*_{G_\mathbb{R}} X, \mathbb{Z}) \\
@V(\mathrm{gr}\,\mu_*)_\lambda VV @VV(\mathrm{gr}\,\mu_*)_\lambda V \\
\bigoplus_k \mathrm{H}^{\mathrm{inf}}_{2k}(\mathcal{N}_k \cap \mathfrak{p}, \mathbb{C}) @>\phi>> \bigoplus_k \mathrm{H}^{\mathrm{inf}}_{2k}(\mathcal{N}_k \cap i\mathfrak{g}_\mathbb{R}, \mathbb{C})\,.
\end{CD}$$

In Section 6, we saw that $\phi$ induces the Kostant-Sekiguchi correspondence when we make the identifications (6.1).

The parameter $\lambda \in \mathfrak{h}^*$ in (7.1) fixes the infinitesimal character of representations, but not conversely: $\chi_\lambda = \chi_\mu$ when $\lambda$ and $\mu$ are $W$-conjugate. The particular choice of $\lambda$ within its $W$-orbit has not mattered until now. At this point, however, it will become convenient to suppose that $\lambda$ is integrally dominant, in the sense that

$$(7.2) \qquad 2\frac{(\lambda,\alpha)}{(\alpha,\alpha)} \notin \mathbb{Z}_{<0} \quad \text{for every } \alpha \in \Phi^+ \,;$$

here, as before, $\Phi^+$ refers to the universal positive root system. This situation is special for the Beilinson-Bernstein construction [BB1]. First of all, it implies

$$(7.3) \qquad \mathrm{H}^p(X,\mathfrak{M}) \;=\; 0 \quad \text{if } \; p \neq 0\,,$$

for every coherent $\mathcal{D}_\lambda$-module $\mathfrak{M}$; in particular, for every $K$-equivariant, coherent $\mathcal{D}_\lambda$-module. When $\lambda$ is not only integrally dominant, but also regular, the assignment $\mathfrak{M} \mapsto \mathrm{H}^0(X,\mathfrak{M})$ establishes an equivalence of categories between the category of $K$-equivariant, coherent $\mathcal{D}_\lambda$-modules on one hand, and the category of Harish-Chandra modules with infinitesimal character $\chi_\lambda$ on the other. When $\lambda$ is integrally dominant but singular, there exist $\mathcal{D}_\lambda$-modules without sections; however,

(7.4) for each irreducible Harish-Chandra module $M$ with infinitesimal character $\chi_\lambda$, there exists a unique irreducible $K$-equivariant $\mathcal{D}_\lambda$-module $\mathfrak{M}$ such that $\mathrm{H}^0(X,\mathfrak{M}) = M$.

Concretely, $\mathfrak{M}$ is the unique irreducible $K$-equivariant quotient of the Beilinson-Bernstein localization of $M$ which does have nonzero sections.

We use the bilinear form $B$ to identify $\mathfrak{g} \cong \mathfrak{g}^*$, as in Section 5. Correspondingly, we identify $K$-orbits in $\mathcal{N} \cap \mathfrak{p}$ with $K$-orbits in $\mathcal{N}^* \cap \mathfrak{p}^*$ and $G_\mathbb{R}$-orbits in $\mathcal{N} \cap i\mathfrak{g}_\mathbb{R}$ with $G_\mathbb{R}$-orbits in $\mathcal{N}^* \cap i\mathfrak{g}_\mathbb{R}{}^*$. Thus we can think of the associated cycle $\mathrm{Ass}(\pi)$ and the wave front cycle $\mathrm{WF}(\pi)$ as lying in the two groups of nilpotent cycles in the bottom row of (7.1). Recall the definition (2.4) of the $K$-equivariant de Rham functor.

7.5 Proposition. *Let $\lambda$ be integrally dominant, $\pi$ an irreducible representation with infinitesimal character $\chi_\lambda$, and $M$ the Harish-Chandra module of $\pi$. With $\mathfrak{M}$ as in* (7.4), *set $\mathcal{F} = \mathrm{DR}(\mathfrak{M})$. Then* $(\mathrm{gr}\,\mu_*)_\lambda(\mathrm{CC}(\mathcal{F})) = \mathrm{Ass}(\pi)$ *via the identification* (6.1a).

This is essentially a reformulation of a result of J.-T. Chang [C1]. At the end of this section, we shall reduce our statement to Chang's result, and comment on certain aspects of his proof.

7.6 Proposition. *Let $\pi$ be an irreducible representation with infinitesimal character $\chi_\lambda$, and $\mathcal{F} \in \mathrm{D}_{G_\mathbb{R}}(X)_{-\lambda}$ a sheaf such that $\beta(\mathcal{F}) = \pi$ up to infinitesimal equivalence. Then either* $(\mathrm{gr}\,\mu_*)_\lambda(\mathrm{CC}(\mathcal{F}))$ *vanishes or* $(\mathrm{gr}\,\mu_*)_\lambda(\mathrm{CC}(\mathcal{F})) = \mathrm{WF}(\pi)$ *via the identification* (6.1b).

We shall see, after the fact, that the first alternative, i.e., the vanishing of $(\mathrm{gr}\,\mu_*)_\lambda(\mathrm{CC}(\mathcal{F}))$, cannot happen when $\mathcal{F}$ is chosen appropriately. We shall deduce the proposition from our integral formula for characters [SV4]. In the case of complex groups and regular infinitesimal character, (7.6) is due to Rossmann [R2]. Our proof is a generalization of Rossmann's argument.

Before turning to the proof of (7.6), let us argue that the two propositions, together with Theorem 6.3 and the commutativity of the diagram (7.1), do imply Theorem 1.4.

*Proof of Theorem* 1.4. We consider a particular irreducible representation $\pi$, with infinitesimal character $\chi_\lambda$. From the construction of the associated cycle, it is clear that $\mathrm{Ass}(\pi) \neq 0$. With $\lambda$ and $\mathcal{F}$ as in Proposition 7.5, $(\mathrm{gr}\,\mu_*)_\lambda(\mathrm{CC}(\mathcal{F})) = \mathrm{Ass}(\pi) \neq 0$. The commutativity of (7.1) now ensures that $(\mathrm{gr}\,\mu_*)_\lambda(\mathrm{CC}(\gamma\,\mathcal{F})) \neq 0$; hence $(\mathrm{gr}\,\mu_*)_\lambda(\mathrm{CC}(\gamma\,\mathcal{F})) = \mathrm{WF}(\pi)$ by Proposition 7.6. We appeal once more to the commutativity of (7.1) to conclude $\mathrm{WF}(\pi) = \phi(\mathrm{Ass}(\pi))$. Because of Theorem 6.3, this gives the assertion of the theorem.

*Proof of Proposition* 7.6. Let $\Theta_\pi$ denote the character of $\pi$, and $\theta_\pi$ the pullback of $\Theta_\pi$ to the Lie algebra $\mathfrak{g}_\mathbb{R}$,

$$\theta_\pi = \sqrt{\det(\exp_*)}\, \exp^* \Theta_\pi\,. \tag{7.7}$$

Our integral formula for characters[3] ( [SV4]), transferred to $\mathfrak{g}$ via $\mathfrak{g} \cong \mathfrak{g}^*$, asserts

$$\int_{\mathfrak{g}_\mathbb{R}} \theta_\pi\, \phi\, dx \;=\; \frac{1}{(2\pi i)^n n!} \int_{\mathrm{CC}(\mathcal{F})} \mu_\lambda^* \hat{\phi}\; (-\sigma + \pi^* \tau_\lambda)^n\,, \tag{7.8}$$

for every test function $\phi \in C_c^\infty(\mathfrak{g}_\mathbb{R})$. Here $\mu_\lambda : T^*X \to \mathfrak{g}$ denotes Rossmann's twisted moment map, $n$ the complex dimension of $X$, $\sigma$ the canonical holomorphic symplectic form on $T^*X$, $\pi : T^*X \to X$ the natural projection, and finally $\tau_\lambda$ a particular differential form on $X$ such that

$$\frac{\tau_\lambda}{2\pi i} \quad \text{represents the cohomology class } c_1(\lambda) \in \mathrm{H}^2(X, \mathbb{C})\,. \tag{7.9}$$

[3] Our formula is the explicit version of a formula of Rossmann [R1], [R2], who represents invariant eigendistributions on $\mathfrak{g}_\mathbb{R}$, with regular infinitesimal character, as integrals over unspecified cycles.

We define the Fourier transform $\hat{\phi}$ of $\phi$ as a holomorphic function on $\mathfrak{g}$ by the formula

$$\hat{\phi}(\xi) \;=\; \int_{\mathfrak{g}_{\mathbb{R}}} e^{B(\xi,x)}\,\phi(x)\,dx\,. \tag{7.10}$$

Since we have omitted the customary factor of $i$ in the exponent, $\hat{\phi}$ decays rapidly along the imaginary directions. This decay property of $\hat{\phi}$ makes the second integral in (7.8) converge.

The definition of the wave front cycle involves scaling the argument of $\theta_\pi$ by a scaling parameter $t \in \mathbb{R}_{>0}$. Let $m_t : \mathfrak{g} \to \mathfrak{g}$ denote scaling by $t$, i.e., $m_t(\xi) = t\xi$. Since

$$\int_{\mathfrak{g}_{\mathbb{R}}} \theta_\pi(tx)\,\phi(x)\,dx \;=\; t^{-d}\int_{\mathfrak{g}_{\mathbb{R}}} \theta_\pi(x)\,\phi(t^{-1}x)\,dx \qquad (\,d = D\dim \mathfrak{g}_{\mathbb{R}}\,),$$

and since

$$\xi \;\mapsto\; t^d\,\hat{\phi}(t\xi) \quad \text{is the Fourier transform of} \quad x \;\mapsto\; \phi(t^{-1}x)\,,$$

the scaled family of invariant eigendistributions is given by the formula

$$\int_{\mathfrak{g}_{\mathbb{R}}} \theta_\pi(tx)\,\phi(x)\,dx \;=\; \frac{1}{(2\pi i)^n n!}\int_{\mathrm{CC}(\mathcal{F})} \mu_\lambda^* m_t^*\hat{\phi}\;(-\sigma + \pi^*\tau_\lambda)^n\,. \tag{7.11}$$

Scaling of cotangent vectors by $t$ defines a map on $T^*X$; for convenience, we denote this map also by the symbol $m_t$. The definition of the twisted moment map implies

$$m_t \circ \mu_\lambda \;=\; \mu_{t\lambda} \circ m_t \qquad \text{and} \qquad \lim_{t\to 0}\mu_{t\lambda} \;=\; \mu\,. \tag{7.12}$$

In particular,

(7.13)
$$\begin{aligned} \mu_\lambda^* m_t^*\hat{\phi}\;(-\sigma + \pi^*\tau_\lambda)^n \;&=\; m_t^*\mu_{t\lambda}^*\hat{\phi}\;(-\sigma + \pi^*\tau_\lambda)^n \\ &=\; m_t^*\left(\mu_{t\lambda}^*\hat{\phi}\;m_{t^{-1}}^*(-\sigma + \pi^*\tau_\lambda)^n\right) \;=\; m_t^*\left(\mu_{t\lambda}^*\hat{\phi}\;(-t^{-2}\sigma + \pi^*\tau_\lambda)^n\right)\,. \end{aligned}$$

In the last step we have used the identities $m_{t^{-1}}^*\pi^*\tau_\lambda = \pi^*\tau_\lambda$ and $m_{t^{-1}}^*\sigma = t^{-2}\sigma$; the former follows from $\pi \circ m_{t^{-1}} = \pi$, the latter from the definition of $\sigma$. By (7.11) and (7.13),

$$\int_{\mathfrak{g}_{\mathbb{R}}} \theta_\pi(tx)\,\phi(x)\,dx \;=\; \frac{1}{(2\pi i)^n n!}\int_{\mathrm{CC}(\mathcal{F})} \mu_{t\lambda}^*\hat{\phi}\;(-t^{-2}\sigma + \pi^*\tau_\lambda)^n\,, \tag{7.14}$$

since $(m_t)_*\,\mathrm{CC}(\mathcal{F}) = \mathrm{CC}(\mathcal{F})$ – recall: the characteristic cycle is invariant under scaling by a positive factor.

As in Section 5, we let $k = k(\mathrm{CC}(\mathcal{F}))$ denote the least integer such that $\mu^{-1}(\tilde{\mathcal{N}}_k)$ contains the support of $\mathrm{CC}(\mathcal{F})$. On $T_K^*X \cap \mu^{-1}(\mathcal{N}_\ell)$,

$$\frac{1}{n!}(-t^{-2}\sigma + \pi^*\tau_\lambda)^n \;=\; \frac{t^{-2\ell}}{\ell!(n-\ell)!}(-\sigma)^\ell \wedge (\pi^*\tau_\lambda)^{n-\ell}\,.$$

Thus we can rewrite the integral on the right in (7.14) as a sum,

$$\begin{aligned}(7.15)\qquad &\int_{\mathfrak{g}_{\mathbb{R}}} \theta_\pi(tx)\,\phi(x)\,dx \\ &\qquad = \sum_{\ell\le k} \frac{t^{-2\ell}}{(2\pi i)^n \ell!(n-\ell)!} \int_{\mathrm{CC}(\mathcal{F})\cap\mu^{-1}(\mathcal{N}_\ell)} \mu_{t\lambda}^*\hat{\phi}\,(-\sigma)^\ell \wedge (\pi^*\tau_\lambda)^{n-\ell}\,.\end{aligned}$$

To get the full asymptotic expansion of $m_t^*\theta_\pi$, one can expand $\mu_{t\lambda}^*\hat{\phi}$ as a Taylor series in $t$. There are no convergence problems, even though we integrate over cycles with infinite support: the remainder for the truncated Taylor series involves various partial derivatives of $\hat{\phi}$, which satisfy the same kind of bound as $\hat{\phi}$ itself; see [SV4, (3.15–3.16)], where the convergence is deduced from the rapid decay of $\hat{\phi}$. In particular, since $\mu_{t\lambda} \to \mu$ ,

$$\begin{aligned}(7.16)\qquad &\int_{\mathfrak{g}_{\mathbb{R}}} \theta_\pi(tx)\,\phi(x)\,dx \\ &\quad = \frac{t^{-2k}}{(2\pi i)^n k!(n-k)!} \int_{\mathrm{CC}(\mathcal{F})\cap\mu^{-1}(\mathcal{N}_k)} \mu^*\hat{\phi}\,(-\sigma)^k \wedge (\pi^*\tau_\lambda)^{n-k} \;+\; O(t^{1-2k})\end{aligned}$$

as $t \to 0$. The integral on the right is therefore either zero or the leading term in the asymptotic expansion of the left-hand side.

To relate the integral on the right to the wave front cycle of $\pi$, we must express it as linear combination of integrals of $\hat{\phi}$ over $G_{\mathbb{R}}$-orbits in $i\mathfrak{g}_{\mathbb{R}} \cap \mathcal{N}_k$ , in each case with respect to the canonical measure of the orbit in question. We consider a particular $G_{\mathbb{R}}$-orbit $\mathcal{O}_{\mathbb{R}}$ in $i\mathfrak{g}_{\mathbb{R}} \cap \mathcal{N}_k$, and let $\mathcal{O}$ denote the $G$-orbit in which $\mathcal{O}_{\mathbb{R}}$ lies. According to [SV3, Lemma 8.19],

$$\mu^*\sigma_{\mathcal{O}} \;=\; -\sigma|_{\mu^{-1}(\mathcal{O})}\,.$$

Hence, on $\mu^{-1}(\mathcal{O}_{\mathbb{R}})$,

$$(7.17)\qquad \begin{aligned}&\frac{1}{(2\pi i)^k k!}\,\mu^*\hat{\phi}\,(-\sigma)^k \;=\; \mu^*(\hat{\phi}\,dm_{\mathcal{O}_{\mathbb{R}}})\,, \\ &\text{where}\quad dm_{\mathcal{O}_{\mathbb{R}}} \;=\; \frac{(\sigma_{\mathcal{O}})^k}{(2\pi i)^k k!} \quad \text{is the canonical measure on } \mathcal{O}_{\mathbb{R}}\,.\end{aligned}$$

Note that $(2\pi i)^{-1}\sigma_{\mathcal{O}}$ is a real, nondegenerate 2-form on $\mathcal{O}_{\mathbb{R}}$, whose top exterior power orients $\mathcal{O}_{\mathbb{R}}$; this orientation allows us to regard the top exterior power as positive measure. The restriction of $\mathrm{CC}(\mathcal{F})$ to $\mu^{-1}(\mathcal{O}_{\mathbb{R}})$ can be regarded, locally, as a product of $\mathcal{O}_{\mathbb{R}}$, oriented as above, and a top dimensional cycle $\mathrm{CC}(\mathcal{F})(\zeta)$ in the Springer fiber $\mu^{-1}(\zeta)$ over any particular $\zeta \in \mathcal{O}_{\mathbb{R}}$; this fibration was used already in the definition of the map $(\mathrm{gr}\,\mu_*)_\lambda$ in Section 5. We now appeal to (7.9) and the definition of $e^\lambda$ as $1 + c_1(\lambda) + \dots$ , and conclude

$$(7.18)\qquad \frac{1}{(2\pi i)^{n-k}(n-k)!} \int_{\mathrm{CC}(\mathcal{F})(\zeta)} \mu^*\hat{\phi}\,\tau_\lambda^{n-k} \;=\; \hat{\phi}(\zeta) \int_{\mathrm{CC}(\mathcal{F})(\zeta)} e^\lambda\,.$$

The characteristic cycle $\mathrm{CC}(\mathcal{F})$ of the $G_{\mathbb{R}}$-invariant sheaf $\mathcal{F}$ is $G_{\mathbb{R}}$-invariant. It follows that the integral of $e^{\lambda}$ over $\mathrm{CC}(\mathcal{F})(\zeta)$ depends on the orbit $\mathcal{O}_{\mathbb{R}}$, not the particular choice of $\zeta \in \mathcal{O}_{\mathbb{R}}$. Let $b(\mathcal{O}_{\mathbb{R}})$ denote the value of this integral. Combining (7.17–7.18), we find

$$\begin{aligned} &\frac{1}{(2\pi i)^n k!(n-k)!} \int_{\mathrm{CC}(\mathcal{F}) \cap \mu^{-1}(\mathcal{O}_{\mathbb{R}})} \mu^* \hat{\phi}\, (-\sigma)^k \wedge (\pi^* \tau_\lambda)^{n-k} \\ &\qquad = \ b(\mathcal{O}_{\mathbb{R}}) \int_{\mathcal{O}_{\mathbb{R}}} \hat{\phi}\, dm_{\mathcal{O}_{\mathbb{R}}}. \end{aligned} \tag{7.19}$$

The assignment $\phi \mapsto \int_{\mathcal{O}_{\mathbb{R}}} \hat{\phi}\, dm_{\mathcal{O}_{\mathbb{R}}}$ defines a distribution, the Fourier transform of the orbit $\mathcal{O}_{\mathbb{R}} \subset i\mathfrak{g}_{\mathbb{R}}$ – more precisely, of the canonical measure $dm_{\mathcal{O}_{\mathbb{R}}}$ on $\mathcal{O}_{\mathbb{R}}$. We normalize the Fourier transform as in (7.10), without a factor $i$ in the exponent.

Let us summarize what we have established so far. Taking the sum of the expressions (7.19) for all $G_{\mathbb{R}}$-orbits $i\mathfrak{g}_{\mathbb{R}} \cap \mathcal{N}_k$, we find

$$\begin{aligned} \phi \quad \mapsto \quad &\frac{1}{(2\pi i)^n k!(n-k)!} \int_{\mathrm{CC}(\mathcal{F}) \cap \mu^{-1}(\mathcal{N}_k)} \mu^* \hat{\phi}\, (-\sigma)^k \wedge (\pi^* \tau_\lambda)^{n-k} \\ &\text{is the Fourier transform of} \quad \sum_{\mathcal{O}_{\mathbb{R}} \subset i\mathfrak{g}_{\mathbb{R}} \cap \mathcal{N}_k} b(\mathcal{O}_{\mathbb{R}})\ \mathcal{O}_{\mathbb{R}} \end{aligned} \tag{7.20}$$

with $\mathcal{O}_{\mathbb{R}}$ shorthand for the distribution on $i\mathfrak{g}_{\mathbb{R}}$ defined by the measure $dm_{\mathcal{O}_{\mathbb{R}}}$. We do not yet know that

$$\sum_{\mathcal{O}_{\mathbb{R}} \subset i\mathfrak{g}_{\mathbb{R}} \cap \mathcal{N}_k} b(\mathcal{O}_{\mathbb{R}})\ \mathcal{O}_{\mathbb{R}} \ \neq \ 0.$$

However, when this holds, (7.16) shows that the distribution (7.20) is the Fourier transform of the leading term of $m_t^* \theta_\pi$ as $t \to 0$. In other words,

$$\sum_{\mathcal{O}_{\mathbb{R}} \subset i\mathfrak{g}_{\mathbb{R}} \cap \mathcal{N}_k} b(\mathcal{O}_{\mathbb{R}})\ \mathcal{O}_{\mathbb{R}} \ = \ \begin{cases} \mathrm{WF}(\pi) & \text{or} \\ 0 & . \end{cases} \tag{7.21}$$

The definition of the constants $b(\mathcal{O}_{\mathbb{R}})$ tells us that $\sum b(\mathcal{O}_{\mathbb{R}})\ [\mathcal{O}_{\mathbb{R}}]$ is the cycle in $\mathcal{N}_k \cap i\mathfrak{g}_{\mathbb{R}}$ obtained from $\mathrm{CC}(\mathcal{F})$ by integrating $e^{\lambda}$ over the fibers of $\mu$, as defined in Section 5. Thus, in the notation of Section 5,

$$\sum_{\mathcal{O}_{\mathbb{R}} \subset i\mathfrak{g}_{\mathbb{R}} \cap \mathcal{N}_k} b(\mathcal{O}_{\mathbb{R}})\ [\mathcal{O}_{\mathbb{R}}] \ = \ (\mathrm{gr}\, \mu_*)_\lambda(\mathrm{CC}(\mathcal{F})). \tag{7.22}$$

Proposition (7.6) now follows from (7.21-22).

*Remarks on Proposition* 7.5. The characteristic cycle of a holonomic $\mathcal{D}$-module is a local invariant, i.e., local with respect to the base manifold $X$. In particular, the holonomic $\mathcal{D}_\lambda$-module $\mathfrak{M}$ has a well-defined characteristic cycle $\mathrm{CC}(\mathfrak{M})$. By construction, it is a complex algebraic cycle of the same dimension as $X$. We regard it as a geometric cycle, in other words, as cycle in $\mathrm{H}^{\inf}_{2n}(T^*_K X, \mathbb{Z})$, by orienting its components via the complex structure. A result of Kashiwara [K2, §8.2] asserts that $\mathrm{CC}(\mathfrak{M})$ coincides with the characteristic cycle of $\mathcal{F} = \mathrm{DR}(\mathfrak{M})$:

$$\mathrm{CC}(\mathfrak{M}) \ = \ \mathrm{CC}(\mathrm{DR}(\mathfrak{M}))\,. \tag{7.23}$$

We shall give a short proof, based on a result of Ginzburg and on [SV3, Theorem 4.2], which we have used already.

The reason for (7.23) is simple: both notions of characteristic cycle obey the same formalism. We shall establish the equality for every (algebraic) regular holonomic $\mathcal{D}$-module on a complex algebraic manifold $X$. Since the characteristic cycles are local invariants, we may as well assume that $X$ is affine. To begin with, we suppose that $\mathfrak{M}$ is the $\mathcal{D}$-module direct image of a vector bundle $E_S$, with a flat algebraic connection, on a closed irreducible submanifold $S \subset X$. In that particular case,

(7.24) $\mathrm{DR}(\mathfrak{M})$ is the sheaf of flat sections of $E_S$ in degree $\mathrm{codim}_{\mathbb{C}}(S,X)$,

as can be computed directly. Thus, by (4.3a),

$$\mathrm{CC}(\mathrm{DR}(\mathfrak{M})) \ = \ (-1)^{\mathrm{codim}_{\mathbb{C}}(S,X)}\, \mathrm{rk}_{\mathbb{C}}(E_S)[T^*_S X]\,, \tag{7.25}$$

with $T^*_S X$ oriented according to our general convention, as in [SV3, (2.3)]. On the other hand,

(7.26) the $\mathcal{D}$-module characteristic cycle $\mathrm{CC}(\mathfrak{M})$ is the conormal bundle $T^*_S X$, oriented by its complex structure, with multiplicity $\mathrm{rk}_{\mathbb{C}}(E_S)$.

The two orientations of $T^*_S X$ – by our general convention for orienting conormal bundles and via the complex structure – are related by $(-1)^{\mathrm{codim}_{\mathbb{C}}(S,X)}$; see (4.34–4.35). This implies (7.23) in the special case of a flat vector bundle on a closed submanifold.

Next we suppose that $\mathfrak{M}$ is the $\mathcal{D}$-module direct image of a regular holonomic $\mathcal{D}$-module $\mathfrak{N}$ on $U$, the complement of a divisor $\{f = 0\}$ in $X$, and further, that $\mathfrak{N}$ satisfies (7.23) on $U$. Let $j$ denote the embedding $U \hookrightarrow X$. Then

$$\mathrm{DR}(\mathfrak{M}) \ = \ Rj_*\, \mathrm{DR}(\mathfrak{N})\,; \tag{7.27}$$

see, for example, [Bo]. Thus our open embedding theorem [SV3, (4.2)] implies

$$\mathrm{CC}(\mathrm{DR}(\mathfrak{M})) \ = \ \lim_{s \to 0^+} \big(\, \mathrm{CC}(\mathrm{DR}(\mathfrak{N})) \ + \ s\, d \log |f|^2 \,\big)\ . \tag{7.28}$$

Ginsburg's theorem [Gi, Theorem 3.2], which inspired our theorem, asserts

$$\mathrm{CC}(\mathfrak{M}) \;=\; \lim_{s\to 0}\left(\mathrm{CC}(\mathfrak{N}) \;+\; s\,\frac{df}{f}\right). \tag{7.29}$$

The equality (7.28) takes place in the real cotangent bundle of $X$, and (7.27) in the holomorphic cotangent bundle. Our convention (3.1) for identifying the two bundles identifies the two differentials $d\log|f|^2$ and $\frac{df}{f}$. Since $\mathrm{CC}(\mathfrak{N}) = \mathrm{CC}(\mathrm{DR}(\mathfrak{N}))$ by assumption,

$$\mathrm{CC}(\mathrm{DR}(\mathfrak{N})) \;+\; s\,d\log|f|^2 \;\cong\; \mathrm{CC}(\mathfrak{N}) \;+\; s\,\frac{df}{f} \qquad (\, s\in\mathbb{R}\,) \tag{7.30}$$

via $T^*(X^{\mathbb{R}}) \cong T^*X$. The family of cycles whose limit we take in (7.28) is therefore the restriction to $\mathbb{R}_{\geq 0}$ of the complex family appearing in (7.29). The two notions for taking limits are consistent; hence (7.23) holds for the $\mathcal{D}$-module direct image $\mathfrak{M} = j_*\mathfrak{N}$ if it does for $\mathfrak{N}$.

Beginning with flat vector bundles on closed submanifolds, one can generate the $K$-group of holonomic $\mathcal{D}$-modules on smooth affine varieties by a succession of direct images under open affine embeddings. This now gives us the equality (7.23) in general.

The identity we just established reduces the assertion (7.5) to the analogous one about the $\mathcal{D}$-module characteristic cycle $\mathrm{CC}(\mathfrak{M})$ of the $K$-equivariant $\mathcal{D}$-module $\mathfrak{M}$. This cycle is $K$-invariant and thus can be expressed as an integral linear combination of conormal bundles of $K$-orbits in $X$,

$$\mathrm{CC}(\mathfrak{M}) \;=\; \textstyle\sum_j \, m_j \, [T^*_{S_j}X]_{\mathbb{C}}\,. \tag{7.31}$$

Here $[T^*_{S_j}X]_{\mathbb{C}}$ is the fundamental cycle of $T^*_{S_j}X$ oriented by its complex structure, not by our general convention for orienting conormal bundles. The moment map exhibits each of these conormal bundles as a fiber bundle over a $K$-orbit $\mathcal{O}_{\mathfrak{p},j}$ in $\mathcal{N}\cap\mathfrak{p}$,

$$T^*_{S_j} \;@>\; F_j \;>>\; \mathcal{O}_{\mathfrak{p},j}\,, \tag{7.32}$$

with fiber $F_j$, which is a union of irreducible components of the Springer fiber $\mu^{-1}\zeta$ over any particular $\zeta\in\mathcal{O}_{\mathfrak{p},j}$. Note that several conormal bundles may lie over the same $K$-orbit in $\mathcal{N}\cap\mathfrak{p}$; in other words, as the index $j$ enumerates $K$-orbits $S_j$ in $X$, there may be repetition among the $\mathcal{O}_{\mathfrak{p},j}$. J.-T. Chang proves

$$\mathrm{Ass}(\pi) \;=\; \textstyle\sum_j m_j \int_{F_j} e^{\lambda}\, [\mathcal{O}_{\mathfrak{p},j}] \tag{7.33}$$

[C1, (2.5.6)]. Chang's actual statement relates two homogenous polynomials, one of which expresses the multiplicity of a $K$-orbit in $\mathrm{Ass}(\pi)$ as $\pi$ runs over the coherent family generated by $\pi$. Chang expresses this polynomial in terms of the integral of $e^\lambda$ over fibers $F_j$, regarded as a polynomial in the variable $\lambda$, which ranges over $\mathfrak{h}^*$. The polynomial identity, evaluated at the localization parameter $\lambda$, reduces to (7.33). Going back to our definition of integration over the Springer fiber, one finds

$$\sum_j m_j \int_{F_j} e^\lambda \, [\mathcal{O}_{\mathfrak{p},j}] \;=\; (\mathrm{gr}\,\mu_*)_\lambda(\mathrm{CC}(\mathfrak{M}))\,. \tag{7.34}$$

Proposition 7.5 now follows from (7.23), Chang's identity (7.33), and the reinterpretation (7.34) of the right-hand side of (7.33).

Harvard University, Cambridge, MA
*E-mail address*: schmid@math.harvard.edu

Northwestern University, Evanston, IL
*E-mail address*: vilonen@math.nwu.edu

References

[BV] D. Barbasch and D. Vogan, The local structure of characters, *J. Funct. Anal.* **37** (1980), 27–55.

[B] A. Beilinson, Localization of representations of reductive Lie algebras, *Proc. of the Internat. Congr. of Math.*, Vol. 1, 2 (Warsaw, 1983), PWN, Warsaw, 1984, 699–710.

[BB1] A. Beilinson and J. Bernstein, Localisation de $\mathfrak{g}$–modules, *C. R. Acad. Sci. Paris* **292** (1981), 15–18.

[BB2] ———, A generalization of Casselman's submodule theorem, in *Representation Theory of Reductive Groups*, *Progr. in Math.* **40**, Birkhäuser Boston, Boston, MA., 1983, 35–52.

[BL] J. Bernstein and V. Lunts, *Equivariant Sheaves and Functors, Lecture Notes in Math.* **1578**, Springer-Verlag, New York, 1994.

[Bo] A. Borel et al., *Algebraic D-Modules*, *Perspectives in Math.*, Academic Press, New York, 1987.

[C1] J.-T. Chang, Asymptotics and characteristic cycles for representations of complex groups, *Compositio Math.* **88** (1993), 265–283.

[C2] ———, Large components of principal series and characteristic cycles, *Proc. A.M.S.* **127** (1999), 3367–3373.

[DM] L. van den Dries and C. Miller, Geometric categories and 0-minimal structures, *Duke Math. J.* **84** (1996), 497–540.

[DMM] L. van den Dries, A. Macintyre, and D. Marker, The elementary theory of restricted analytic fields with exponentiation, *Ann. of Math.* **140** (1994), 183–205.

[HC] Harish-Chandra, On some applications of the universal enveloping algebra of a semisimple Lie algebra, *Trans. A.M.S.* **70** (1951), 28–96.

[He] S. Helgason, *Differential Geometry, Lie Groups, and Symmetric Spaces*, Academic Press, New York, 1978.

[K1] M. Kashiwara, The Riemann-Hilbert problem for holonomic systems, *Publ. RIMS* **20** (1984), 319–365.

[K2] ———, Index theorem for constructible sheaves, in *Systèmes Differentiels et Singularités* (A. Galligo, M. Maisonobe, and Ph. Granger, eds.), *Astérisque* **130** (1985), 193–209.

[K3] ———, Open problems in group representation theory, *Proc. of Taniguchi Symposium*, 1996, RIMS preprint **569**, Kyoto University, 1987.

[K4] ———, Character, character cycle, fixed point theorem, and group representations, *Adv. Stud. Pure Math.* **14**, Kinokuniya, Tokyo, 1988, 369–378.

[K5] ———, $\mathcal{D}$-modules and representation theory of Lie groups, *Ann. Inst. Fourier* (Grenoble) **443** (1993), 1597–1618.

[KSa] M. Kashiwara and P. Schapira, *Sheaves on Manifolds*, Springer-Verlag, New York, 1990.

[KSd] M. Kashiwara and W. Schmid, Quasi-equivariant $\mathcal{D}$-modules, equivariant derived category, and representations of reductive Lie groups, in *Lie Theory and Geometry, in Honor of Bertram Kostant*, *Progr. in Math.* **123**, Birkhäuser Boston, Boston, MA, 1994, 457–488.

[Ko] B. Kostant, The principal three-dimensional subgroup and the Betti numbers of a complex simple Lie group, *Amer. J. Math.* **81** (1959), 973–1032.

[KR] B. Kostant and S. Rallis, Orbits and representations associated with symmetric spaces, *Amer. J. Math.* **93** (1971), 753–809.

[Ma] T. Matsuki, Orbits on affine symmetric spaces under the action of parabolic subgroups, *Hiroshima Math. J.* **12** (1982), 307–320.

[Me] Z. Mebkhout, Une équivalance de catégories - Une autre équivalence de catégories, *Compositio Math.* **51** (1984 ), 51–62.

[MUV] I. Mirković, T. Uzawa, and K. Vilonen, Matsuki correspondence for sheaves, *Invent. Math.* **109** (1992), 231–245.

[N] L. Ness, A stratification of the null cone via the moment map, *Amer. J. Math.* **106** (1984), 1281–1329.

[R1] W. Rossmann, Characters as contour integrals, *Lie Group Representations*, III, *Lecture Notes in Math.* **1077**, Springer-Verlag, New York, 1984, 375–388.

[R2] ———, Invariant eigendistributions on a semisimple Lie algebra and homology classes on the conormal variety I, II, *J. Funct. Anal.* **96** (1991), 130–193.

[R3] ———, Picard-Lefschetz theory for the coadjoint quotient of a semisimple Lie algebra, *Invent. Math.* **121** (1995), 531–578.

[R4] ———, Picard-Lefschetz theory and characters of a semisimple Lie group, *Invent. Math.* **121** (1995), 579–611.

[S] W. Schmid, Boundary value problems for group invariant differential equations, *The Mathematical Heritage of Élie Cartan*, *Astérisque*, Numero Hors Série (1985), 311–321.

[SV1] W. Schmid and K. Vilonen, Characters, fixed points and Osborne's conjecture, *Representation Theory of Groups and Algebras*, *Contemp. Math.* **145** (1993), 287–303.

[SV2] ———, Characters, characteristic cycles and nilpotent orbits, *Geometry, Topology, and Physics*, Internat. Press, Cambridge, MA, 1995, 329–340.

[SV3] ———, Characteristic cycles of constructible sheaves, *Invent. Math.* **124** (1996), 451–502.

[SV4] ———, Two geometric character formulas for reductive Lie groups, *Jour. A.M.S.* **11** (1998), 799–867.

[SV5] ———, On the geometry of nilpotent orbits, *Asian J. Math.* **3** (1999), 233–274.

[Se] J. Sekiguchi, Remarks on nilpotent orbits of a symmetric pair, *J. Math. Soc. Japan* **39** (1987), 127–138.

[Spal] N. Spaltenstein, On the fixed point set of a unipotent element on the variety of Borel subgroups, *Topology* **16** (1977), 203–204.

[Span] E. Spanier, *Algebraic Topology*, McGraw-Hill, New York, 1966.

[V1] D. Vogan, Gelfand-Kirillov dimension for Harish-Chandra modules, *Invent. Math.* **48** (1978), 75–98.

[V2] ———, Oral lectures at Bowdoin College, August 1989.

[V3] ———, Associated varieties and unipotent representations, in *Harmonic Analysis on Reductive Groups*, *Progr. in Math.* **101**, Birkhäuser Boston, Boston, MA, 1991, 315–388.